\documentclass[PhD] {uclathes}

\usepackage{amsmath,amsthm}

\usepackage{amssymb}

\usepackage{amsfonts}

\usepackage{amscd}

\newtheorem{thm}{Theorem}[section]

\newtheorem{theorem}[thm]{Theorem}

\newtheorem{corollary}[thm]{Corollary}

\newtheorem{lemma}[thm]{Lemma}

\newtheorem{proposition}[thm]{Proposition}

\theoremstyle{definition}

\newtheorem{definition}[thm]{Definition}

\newtheorem{example}[thm]{Example}

\newtheorem{examples}[thm]{Examples}

\newtheorem{remark}[thm]{Remark}

\newtheorem{observation}[thm]{Observation}

\DeclareMathOperator{\End}{End}

\newcommand{\R} {\mathbf{R}}

\newcommand{\Z} {\mathbf{Z}}

\newcommand{\C} {\mathbf{C}}

\newcommand{\A} {\mathbf{A}}

\newcommand{\lie} {\text{Lie}}

\newcommand{\Hom}{\text{Hom}}

\newcommand{\tensor} {\otimes}

\newcommand{\shf}{\mathcal{O}}

\newcommand{\wt}{\widetilde}

\newcommand{\bA} {\mathbf{A}}

\newcommand{\cO}{\mathcal{O}}

\newcommand{\cI}{\mathcal {I}}

\newcommand{\rA}{\mathrm A}

\newcommand{\al}{\alpha}

\newcommand{\be}{\beta}

\newcommand{\s}{\sigma}

\newcommand{\ep}{\epsilon}

\newcommand{\Eta} {\mathrm E}

\newcommand{\GL}{\mathrm{GL}}

\newcommand{\SL}{\mathrm{SL}}

\newcommand{\rSL}{\mathrm{SL}}

\newcommand{\OSp}{\text{OSp}}

\newcommand{\bP}{\mathbf{P}}

\newcommand{\defi} {\text def}

\newcommand{\lra} {\longrightarrow}

\newcommand{\del} {\partial}

\newcommand{\spec}{{\text{Spec}}}

\newcommand{\bm} {\mathbf {m}}

\newcommand{\bp} {\mathbf {p}}

\newcommand{\btm} {\mathbf {M}}

\newcommand{\salg} {\text{ (salg) }}

\newcommand{\sets} {\text{ (sets) }}

\newcommand\sschemes{ \text {(sschemes)} }

\newcommand\sschemesaff{ \text{ (affine sschemes) }}

\newcommand\smod{ \text{(smod) }}

\newcommand{\uspec}{\mathrm{\underline{Spec}}}

\newcommand{\uHom}{\mathrm{\underline{Hom}}}

\newcommand{\rM}{\mathrm{M}}

\newcommand{\id}{\mathrm{id}}

\newcommand{\Der}{\mathrm{Der}}

\newcommand{\der}{\mathrm{Der}}

\newcommand{\Sym}{\mathrm{Sym}}

\newcommand{\sym}{\mathrm{Sym}}

\newcommand{\cF}{\mathcal{F}}

\newcommand{\D}{\Delta}

\newcommand{\rGL}{\mathrm{GL}}

\newcommand{\Ad}{\mathrm{Ad}}

\newcommand{\ad}{\mathrm{ad}}

\newcommand{\Lie}{\mathrm{Lie}}

\newcommand{\Ber}{\mathrm{Ber}}

\newcommand{\fg}{\mathfrak g}

\newcommand{\Span}{\hbox{span}}

\newcommand{\it}{\textit}

\newcommand{\proj}{\text{Proj}}

\begin{document}

  \vspace{3in}

\begin{center}
\textbf{\large{MATHEMATICAL FOUNDATIONS \\ OF SUPERSYMMETRY}}

\vspace{3in}

Lauren Caston and Rita Fioresi \footnote{Investigation supported by the University of Bologna}

\smallskip

\textit{Dipartimento di Matematica, Universit\`a di
Bologna}\\

\textit{Piazza di Porta S. Donato, 5}\\

\textit{40126 Bologna, Italia}\\

{\footnotesize e-mail: caston@dm.UniBo.it, fioresi@dm.UniBo.it}

\end{center}

\newpage

\textbf{Introduction}

\bigskip

Supersymmetry (SUSY) is the machinery mathematicians and physicists have 
developed to treat two types of elementary particles, \textit{bosons} and 
\textit{fermions}, on the same footing.  Supergeometry is the geometric basis 
for supersymmetry; it was first discovered and studied by physicists Wess, 
Zumino \cite{wz}, Salam and Strathde \cite{ss}
(among others) in the early 1970's.  Today 
supergeometry plays an important role in high energy physics.  The objects 
in super geometry generalize the concept of smooth manifolds and algebraic 
schemes to include anticommuting coordinates.  As a result, we employ the 
techniques from algebraic geometry to study such objects, namely A. 
Grothendiek's theory of schemes.

Fermions include all of the material world; they are the building blocks of 
atoms.  Fermions do not like each other. This is in essence the Pauli 
exclusion principle which states that two electrons cannot occupy the same 
quantum mechanical state at the same time. Bosons, on the other hand, can 
occupy the same state at the same time. 

Instead of looking at equations which just describe either bosons or fermions 
separately, supersymmetry seeks out a description of both simultaneously.  
Transitions between fermions and bosons require that we allow transformations 
between the commuting and anticommuting coordinates. Such transitions are 
called supersymmetries.

In classical Minkowski space, physicists classify elementary particles by 
their mass and spin.  Einstein's special theory of relativity requires that 
physical theories must be invariant under the Poincar\'{e} group.  Since 
observable operators (e.g. Hamiltonians) must commute with this action, the 
classification corresponds to finding unitary representations of the 
Poincar\'{e} group.  In the SUSY world, this means that mathematicians are 
interested in unitary representations of the super Poincar\'{e} group.  
A ``super" representation gives a ``multiplet" of ordinary particles which 
include both fermions and bosons.

Up to this point, there have been no colliders that can produce the energy
required to physically expose supersymmetry. However, the Large Hadron
Collider (LHC) in CERN (Geneva, Switzerland) will be operational in 2007.
Physicists are planning proton-proton and proton-antiproton collisions
which will produce energies high enough where it is believed
supersymmetry can be seen.  Such a discovery will solidify supersymmetry as
the most viable path to a unified theory of all known forces.  Even before
the boson-fermion symmetry which SUSY presupposes is made physical fact,
the mathematics behind the theory is quite remarkable.  The concept that
space is an object built out of local pieces with specific local
descriptions has evolved through many centuries of mathematical thought.
Euclidean and non-Euclidean geometry, Riemann surfaces, differentiable
manifolds, complex manifolds, algebraic varieties, and so on represent
various stages of this concept.  In Alexander Grothendieck's theory of
schemes, we find a single structure (a scheme) that encompasses all
previous ideas of space.  However, the fact that conventional descriptions
of space will fail at very small distances (Planck length) has been the
driving force behind the discoveries of unconventional models of space
that are rich enough to portray the quantum fluctuations of space at these
unimaginably small distances.  Supergeometry is perhaps the most highly
developed of these theories; it provides a surprising continuation of the
Grothendieck theory and opens up large vistas.  One should not think of it
as a mere generalization of classical geometry, but as a deep continuation
of the idea of space and its geometric structure.

Out of the first supergeometric objects constructed by the pioneering 
physicists came mathematical models of superanalysis and supermanifolds 
independently by F. A. Berezin \cite{Berezin}, 
B. Kostant \cite{Kostant}, D.A. Leites \cite{Leites}, and 
De Witt \cite{dewitt}.  
The idea to treat a supermanifold as a ringed space with a sheaf of 
$\Z/2\Z$-graded algebras was introduced in these early works.  Later, 
Bernstein \cite{DM}
and Leites made this treatment rigorous and used techniques 
from algebraic geometry to deepen the study of supersymmetry.  
In particular, Bernstein and Leites accented the functor of points approach 
from Grothendieck's theory of schemes.  It is this approach 
(which we call $T$-points) that we present and expand upon in our treatment 
of mathematical supersymmetry.  Interest in SUSY has grown in the past decade, 
and most recently works by V. S. Varadarajan \cite{VSV2}
among others, 
have continued the exploration of the beautiful area of physics and 
mathematics and have inspired this work.  Given the interest and the number 
of people who have contributed greatly to this field from various 
perspectives, it is impossible to give a fair and accurate account of works
related to ours.  

In our exposition of mathematical SUSY, we use the language of $T$-points to 
build supermanifolds up from their foundations in $\Z/2\Z$-graded linear 
algebra (superalgebra).  This treatment is similar to that given by 
Varadarajan in \cite{VSV2}, however we prove some deeper results related 
to the Frobenius theorem as well as give a full treatment of superschemes 
in chapters \ref{sgchap3}-\ref{sgchap4}.  Recently the book by G. Tuynman 
\cite{Tuyn} has been brought to our attention.  The main results from 
chapters \ref{sgchap5}-\ref{sgchap6} can be found in \cite{Tuyn}, however 
we obtained our results independently of this work, moreover, our method of 
$T$-points remains fresh in light of this and other recent works.

Here is a brief description of our work.

In chapter \ref{sgchap1} we begin by studying $\Z/2\Z$-graded linear objects.  
We define super vector spaces and superalgebras, then generalize some 
classical results and ideas from linear algebra to the super setting.  
For example, we define a super Lie algebra, discuss supermatrices, and 
formulate the super trace and determinant (the Berezinian).

In chapter \ref{sgchap2} we introduce the most basic geometric structure: 
a superspace.  We present some general properties of superspaces which leads 
into two key examples of superspaces, supermanifolds and superschemes.  
Here we also introduce the notion of $T$-points which treats our geometric 
objects as functors; it is a fundamental tool to gain geometric intuition 
in supergeometry.

Chapters \ref{sgchap3}-\ref{sgchap4} lay down the full foundations of 
$C^{\infty}$-supermanifolds over $\R$.  We give special attention to super 
Lie groups and their associated Lie algebras, as well as look at how group 
actions translate infinitesimally.  In chapter \ref{sgchap4} we prove the 
local and global Frobenius theorem on supermanifolds, then use the 
infinitesimal actions from chapter \ref{sgchap3} to build the super Lie 
subgroup, subalgebra correspondence.

Chapters \ref{sgchap5}-\ref{sgchap6} expand upon the notion of a superscheme 
which we introduce in chapter \ref{sgchap2}.  We immediately adopt the 
language of $T$-points and give criterion for representability: in 
supersymmetry it is often most convenient to desribe an object functorially, 
then show it is representable.  In chapter \ref{sgchap5}, we explicitly 
construct the Grassmanian superscheme functorially, then use the 
representability criterion to show it is representable.  Chapter \ref{sgchap5} 
concludes with an examination of the infinitesimal theory of superschemes.  
We continue this exploration in chapter \ref{sgchap6} from the point of view 
of algebraic supergroups and their Lie algebras.  We discuss the linear 
representations of affine algebraic supergroups; in particular we show that 
all affine super groups are realized as subgroups of the general linear 
supergroup.

This work is self-contained; we try to keep references to a minimum in the 
body of our work so that the reader can proceed without the aid of other 
texts.  We assume a working knowledge of sheaves, differential geometry, 
and categories and functors.  We suggest that the reader begin with chapters 
\ref{sgchap1} and \ref{sgchap2}, but chapters \ref{sgchap3}-\ref{sgchap4} 
and chapters \ref{sgchap5}-\ref{sgchap6} are somewhat disjoint and may be 
read independently of one another.

We wish to thank professor V. S. Varadarajan for 
introducing us to this beautiful part of mathematics. He has truly inspired 
us through his insight and deep understanding of the subject.  We also wish 
to thank Prof. M. A. Lledo, Prof. A. Vistoli and Prof. M. Duflo
for many helpful remarks.   R. Fioresi thanks the UCLA Department of 
Mathematics for its
kind hospitality during the realization of this work. L. Caston thanks 
the Dipartimento di Matematica, Universita' di Bologna
for support and hospitality during the realization of this work.

\tableofcontents

\chapter{$\Z/2\Z$-Graded Linear Algebra} \label{sgchap1}

We first build the foundations of linear algebra in the super context.  
This is an important starting point as we later build super 
geometric objects from sheaves of super linear spaces.  
Let us fix a ground field $k$, $\text{char}(k) \neq 2,3$.

\section{Super Vector Spaces and Superalgebras}

\begin{definition}
A \textit{super vector space} is a  $\Z/2\Z$-graded vector space
\[ V = V_0 \oplus V_1
\]
where elements of $V_0$ are called ``even" and elements of $V_1$
are called ``odd".
\end{definition}

\begin{definition}
The \textit{parity} of $v \in V$, denoted $p(v)$ or $|v|$, is defined only
on nonzero \textit{homogeneous} elements, that is elements of either $V_0$
or $V_1$:
\[
p(v) = |v| = \left\{ \begin{array}{c} 0 \text{ if $v \in V_0$}\\
1 \text{ if $v \in V_1$} \end{array} \right.
\]
\end{definition}

Since any element may be expressed as the sum of homogeneous
elements, it suffices to only consider homogeneous elements in the statement of definitions, theorems, and proofs.

\begin{definition}
The \textit{super dimension} of a super vector space $V$ is the pair $(p,q)$ where dim($V_0$)=$p$ and dim($V_1$)=$q$ as ordinary vector spaces.  We simply write $\dim(V) =
p|q$.
\end{definition}
From now on we will simply refer to the superdimension as the
dimension when the category is clear.  If $\dim(V)= p|q$, then we
can find a basis $\{ e_1, \ldots, e_p \}$ of $V_0$ and a basis  $\{ \epsilon_1, \ldots, \epsilon_q \}$ of $V_1$ so
that $V$ is canonically isomorphic to the free
$k$-module generated by the $\{ e_1, \ldots, e_p, \epsilon_1,
\ldots, \epsilon_q \}$.  We denote this $k$-module by $k^{p|q}$.

\begin{definition}
A \textit{morphism} from a super vector space $V$ to a super
vector space $W$ is a $\Z/2\Z$-grading preserving linear map from
$V$ to $W$. Let $\text{Hom}(V,W)$ denote the set of morphisms $V
\longrightarrow W$.
\end{definition}

Thus we have formed the abelian category of super vector spaces.
It is important to note that the category of super vector spaces
also admits and ``inner Hom", which we denote
$\text{\underline{Hom}}(V,W)$; it consists of \textit{all} linear
maps from $V$ to $W$:
\[
\begin{array}{lcl}
\text{\underline{Hom}}(V,W)_0 & = & \{ T:V \longrightarrow W
\hspace{.05in} | \hspace{.05in}
\text{$T$ preserves parity} \}\hspace{.1in} (=\text{Hom}(V,W));\\
\text{\underline{Hom}}(V,W)_1 & = & \{ T:V \longrightarrow W
\hspace{.05in} | \hspace{.05in} \text{$T$ reverses parity} \}.
\end{array}
\]
In the category of super vector spaces we have the \textit{parity
reversing functor} $\Pi$ defined by \[(\Pi V)_0 = V_1
\hspace{.3in} (\Pi V)_1 = V_0.
\]

The category of super vector spaces is in fact a tensor category,
where $V \tensor W$ is given the $\Z/2\Z$-grading as follows:
\[
\begin{array}{lcl} (V \tensor W)_0 & = & (V_0 \tensor W_0) \oplus
(V_1 \tensor W_1)\\ (V \tensor W)_1 & = & (V_0 \tensor W_1) \oplus
(V_1 \tensor W_0).
\end{array}
\]
The tensor functor $\otimes$ is additive and exact in each
variable as in the ordinary vector space category; it has a unit
object (i.e. $k$) and is associative.  Moreover, $V\tensor W \cong W \tensor V$ by the
\textit{commutativity map}
\[ c_{V,W}: V \tensor W \longrightarrow
W \tensor V \]
where $v\tensor w \mapsto (-1)^{|v||w|} w  \tensor
v$. This is the so called ``sign rule" that one finds in some physics and math literature. In any tensor
category with an inner Hom, the \textit{dual}, $V^*$, of $V$ is
\[ V^* =_{def} \underline{\text{Hom}}(V,k). \]

\begin{remark}
We understand completely the object $V^{\tensor n} = V \tensor \cdots \tensor V$ ($n$ times) for a super vector space $V$. We can extend this notion to make sense of $V^{\tensor n|m}$ via the parity reversing
functor $\Pi$. Define
\[ V^{n|m} := \underbrace{V \times V \times \ldots \times
V}_{\text{$n$ times}} \times \underbrace{\Pi(V) \times \Pi(V)
\times \ldots \times \Pi(V)}_{\text{$m$ times}},
\]
from which the definition of $V^{\tensor n|m}$ follows by the
universal property.
\end{remark}

Let us now define a super $k$-algebra:
\begin{definition}
We say that a super vector space $A$ is a \textit{superalgebra} if
there is a multiplication morphism $\tau: A \tensor A \longrightarrow A$.
\end{definition}
We then say that a superalgebra $A$ is \textit{commutative} if
\[
\tau \circ c_{A,A} = \tau,
\]
that is, if the product of homogeneous elements obeys the rule
\[
ab = (-1)^{|a||b|}ba.
\]
Similarly we say that $A$ is \textit{associative} if
$\tau \circ \tau \tensor id = \tau \circ id \tensor \tau$ on
$A\tensor A\tensor A$, and that $A$ has a \textit{unit} if there is an even
element $1$ so that $\tau(1\tensor a) = \tau(a\tensor 1) = a$ for
all $a \in A$.

From now on we will assume all superalgebras are
associative and commutative with unit unless specified.

An \textit{even derivation} of a superalgebra $A$ is a super vector space homomorphism $D: A \longrightarrow A$ such that for $a,b \in A$, $D(ab) = D(a)b + aD(b)$.  We may of course extend this definition to include odd linear maps:
\begin{definition}
Let $D \in \underline{\text{Hom}}_k(A, A)$ be a $k$-linear map.  Then $D$ is a \textit{derivation} of the superalgebra $A$ if for $a, b \in A$,
\begin{equation}
\label{der}
D(ab) = D(a)b + (-1)^{|D||a|}aD(b).
\end{equation}
\end{definition}
Those derivations in $\text{Hom}_k(A,A)$ are even (as above) while those in $\underline{\text{Hom}}_k(A,A)_1$ are odd.  The set of all derivations of a superalgebra $A$, denoted $\text{Der}(A)$, is itself a special type of superalgebra called a \textit{super Lie algebra} which we describe in the following section.

\begin{example}
\label{grassex}
\textit{Grassmann coordinates}.
Let
\[ A=k[t_1, \ldots, t_p, \theta_1, \ldots, \theta_q]
\]
where the $t_1, \ldots, t_p$ are ordinary indeterminates and the
$\theta_1, \ldots \theta_q$ are \textit{odd indeterminates}, i.e.
they behave like Grassmannian coordinates:
\[
\theta_i \theta_j = -\theta_j \theta_i.
\]
(This of course implies that $\theta_i^2 = 0$ for all $i$.)  We
claim that $A$ is a supercommutative algebra.  In fact,
\[
A_0 = \{ f_0 + \sum_{|I| \text{ even}} f_I \theta_I | I = \{i_1 < \ldots < i_r\} \}
\]
where $\theta_I = \theta_{i_1}\theta_{i_2}\ldots \theta_{i_r}$ and $f_0, f_I \in k[t_1,
\ldots, t_p]$, and
\[
A_1 = \{\sum_{|J| \text{odd}} f_J \theta_J | J = \{i_1 < \ldots < i_s\} \}
\]
for $s$ odd ($|J| = 2m+1$, $m=1,2,\ldots$) and $f_J \in k[t_1,
\ldots t_q]$. Note that although the $\{\theta_j\} \in A_1$, there
are plenty of nilpotents in $A_0$; take for example
$\theta_1\theta_2 \in A_0$.

Consider the $k$-linear operators $\{ \partial/\partial t_i \}$ and $\{ \partial/\partial \theta_j \}$ of $A$ to itself where
\begin{equation}
\begin{array}{lcl}
\partial/\partial t_i (t_k) = \delta^i_k & \hspace{.3in} & \partial/\partial t_i (\theta_l) = 0;\\
\partial/\partial \theta_j (t_k) = 0 & \hspace{.3in} & \partial/\partial \theta_j (\theta_l) = \delta^j_l.
\end{array}
\end{equation}
It is easy to verify that $\{ \partial/\partial t_i, \partial/\partial \theta_j \} \in \text{Der}(A)$, and we leave it as an exercise to check that
\[
\text{Der}(A) = \text{Span}_A\{\frac{\partial}{\partial t_i}, \frac{\partial}{\partial \theta_j} \}.
\]
\end{example}

\section{Lie Algebras}

An important object in supersymmetry is the super Lie algebra.

\begin{definition}
A \textit{super Lie algebra} $L$ is an object in the category of super
vector spaces together with a morphism $[,]: L \tensor L
\longrightarrow L$ which categorically satisfies the usual conditions.
\end{definition}

It is important to note that in the super category, these
conditions are slightly different to accomodate the odd
variables.  The bracket $[,]$ must satisfy\\
\noindent 1. Anti-symmetry
\[ [,] + [,] \circ c_{L,L} = 0 \]
which may be interpreted as $[x,y] + (-1)^{|x||y|}[y,x] = 0$ for $x, y \in L$ homogeneous.\\
\noindent 2. The Jacobi identity
\[ [,[,]] + [,[,]] \circ \sigma + [,[,]] \circ \sigma^2 = 0 \]
where $\sigma \in S_3$ is a three-cycle, i.e. it takes the first
entry of $[,[,]]$ to the second, the second to the third, and the
third to the first.  So for $x,y,z \in L$ homogeneous, this reads:
\[ [x,[y,z]] + (-1)^{|x||y|+|x||z|}[y,[z,x]] +
(-1)^{|y||z|+|x||z|}[z,[x,y]] = 0. \]

\begin{remark}
We can immediately extend this definition to the case where 
$L$ is an $A$-module. 
\end{remark}

\begin{example}
\label{GrassLie}
In the Grassmannian example above (\ref{grassex}), 
$$
\der(A) =
\text{Span}_A \{ \frac{\partial}{\partial t_i}, 
\frac{\partial}{\partial \theta_j} \}
$$ 
is a super Lie algebra where the bracket is taken for 
$D_1, D_2 \in \text{Der}(A)$ to be $[D_1, D_2] = D_1D_2 - (-1)^{|D_1||D_2|}D_2D_1$.
\end{example}

In fact, we can make any associative algebra $A$ into a 
Lie algebra by taking the bracket to be
\[
[a,b] = ab - (-1)^{|a||b|}ba,
\]
i.e. we take the bracket to be the difference $\tau - \tau \circ c_{A,A}$ where we recall $\tau$ is the multiplication morphism on $A$.  We will discuss other examples of super Lie algebras after the following discussion of superalgebra modules.  In particular we want to examine the SUSY-version 
of a matrix algebra.

\begin{remark}
If the ground field has characteristic 2 o 3 
in addition to the antisymmetry and Jacobi conditions 
one requires that $[x,x]=0$ for $x$ even if the characteristic
is 2 or $[y,[y,y]]=0$ for $y$ odd if the characteristic
is 3. For more details on superalgebras over fields with positive
characteristic see \cite{bmpz}.
\end{remark}

\section{Modules}
\label{Amod}
Let $A$ be a superalgebra, not necessarily commutative in this section.

\begin{definition}
A \textit{left $A$-module} is a super vector space $M$ with a
morphism $A\tensor M \longrightarrow M$ obeying the usual
identities found in the ordinary category.
\end{definition}
A \textit{right $A$-module} is defined similarly.  Note that if
$A$ is commutative, a left $A$-module is also a right
$A$-module using the sign rule
\[
m\cdot a = (-1)^{|m||a|}a\cdot m
\]
for $m \in M$, $a \in A$.  Morphisms of $A$-modules are also
obviously defined, and so we have the category of $A$-modules.  For $A$ commutative, the category of $A$-modules is a
tensor category: for $M_1, M_2$ $A$-modules, $M_1 \tensor M_2$ is
taken as the tensor of $M_1$ as a right module with $M_2$ as a
left module.

Let us now turn our attention to \textit{free} $A$-modules.  We already have the notion of the vector space $k^{p|q}$ over $k$, and so we define $A^{p|q}:= A \tensor k^{p|q}$ where
\[
\begin{array}{lcl}
(A^{p|q})_0 & = & A_0 \tensor (k^{p|q})_0 \oplus A_1
\tensor(k^{p|q})_1\\
(A^{p|q})_1 & = & A_1 \tensor (k^{p|q})_0 \oplus A_0
\tensor(k^{p|q})_1.
\end{array}
\]

\begin{definition}
We say that an $A$-module $M$ is \textit{free} if it is isomorphic (in the category of $A$-modules) to $A^{p|q}$ for some $(p,q)$.
\end{definition}

This definition is equivalent to saying that there are \textit{even} elements $\{e_1, \ldots, e_p
\}$ and \textit{odd} elements $\{\epsilon_1, \ldots,\epsilon_q
\}$ which generate $M$ over $A$.

Let $T:A^{p|q} \longrightarrow A^{r|s}$ be a morphism
of free $A$-modules and write $e_{p+1}, \ldots, e_{p+q}$ for the odd variables $\epsilon_1, \ldots, \epsilon_q$.  Then $T$ is defined on the basis elements
$\{e_1, \ldots e_{p+q} \}$ by
\begin{equation}
\label{Teqn}
T(e_j) = \sum_{i=1}^{p+q} e_i t^i_j.
\end{equation}
Hence $T$ can be represented as a matrix of size $(r+s) \times
(p+q)$:
\begin{equation}
\label{matrixT}
T = \left( \begin{array}{cc}T_1 & T_2\\ T_3 & T_4
\end{array} \right)
\end{equation}
where $T_1$ is an $r \times p$ matrix
consisting of even elements of $A$, $T_2$ is an $r \times q$
matrix of odd elements, $T_3$ is an $s \times p$ matrix of even
elements, and $T_4$ is an $s \times q$ matrix of odd elements.
We say that $T_1$ and $T_4$ are \textit{even blocks} and that $T_2$
and $T_3$ are \textit{odd blocks}.  Because $T$ is a morphism of super $A$-modules, it
must preserve parity, and therefore the parity of the blocks is
determined. Note that when we define $T$ on the basis
elements, in the expression (\ref{Teqn}) the basis element \textit{preceeds} the
coordinates $t^i_j$. This is important to keep the signs in
order and comes naturally from composing morphisms.  For any
$x \in A^{p|q}$, we can express $x$ as the column vector $x = \sum
e_i x^i$ and so $T(x)$ is given by the
matrix product $Tx$. Similarly the composition of morphisms is
given by a matrix product.

\section{Matrices}

Let us now consider all endomorphisms of $M = A^{p|q}$, i.e.
Hom($M,M$).  This is an ordinary algebra (i.e. \textit{not} super)
of matrices of the same type as $T$ above.  Even though in matrix
form each morphism contains blocks of odd elements of $A$, each
morphism is an even linear map from $M$ to itself since a morphism
in the super category must preserve parity.  In order to get a
truly SUSY-version of the ordinary matrix algebra, we must
consider \textit{all} linear maps $M$ to $M$, i.e. we are
interested in \underline{Hom}($M,M$).  Now we can talk about even
and odd matrices.  An even matrix $T$ takes on the block form from
above.  But the parity of the blocks is reversed for an odd matrix $S$; we get
\[
S = \left( \begin{array}{cc} S_1 & S_2\\ S_3 & S_4 \end{array}
\right)
\]
where $S_1$ is a $p \times p$ odd block, $S_4$ is a $q \times q$
odd block, $S_2$ is a $ p\times q$ even block, and $S_3$ is a
$q\times p$ even block.  Note that in the case where $M=k^{p|q}$,
the odd blocks are just zero blocks.  We will denote this super
algebra of even and odd $(p+q) \times (p+q) = p|q \times p|q$
matrices by $\text{Mat}(A^{p|q})$.  This super algebra is in fact
a super Lie algebra where we define the bracket $[ , ]$ as in Example \ref{GrassLie}:
\[
[T,S] = TS - (-1)^{|T||S|}ST
\]
for $S,T \in\text{Mat}(A^{p|q})$.

\begin{remark}
\label{MatNotM}
Note that $\text{Mat}(A^{p|q}) = \underline{\text{Hom}}(A^{p|q}, A^{p|q})$.  We do not want to confuse this with what we will later denote as $M_{p|q}(A)$, which will functorially only include the \textit{even part} of $\text{Mat}(A^{p|q})$, i.e.
\[
\text{Mat}(A^{p|q})_0 = M_{p|q}(A) = \text{Hom}(A^{p|q}, A^{p|q})
\]
(see chapter \ref{sgchap3}).
\end{remark}

We now turn to the SUSY-extensions of the trace and
determinant. Let $T:A^{p|q} \longrightarrow A^{p|q}$ be a morphism
(i.e. $T\in(\text{Mat}(A^{p|q}))_0$) with block form (\ref{matrixT}).
\begin{definition}
We define the \textit{super trace} of $T$ to be:
\begin{equation}
\label{trace}
\text{Tr}(T) := \text{tr}(T_1) - \text{tr}(T_4)
\end{equation}
where ``tr" denotes the ordinary trace.
\end{definition}

This negative sign is
actually forced upon us when we take a categorical view of the
trace.  We will not discuss this here, but we later give
motivation to this definition when we explore the SUSY-extension
of the determinant.

\begin{remark} The trace is actually defined for
\textit{all} linear maps.  For $S \in \text{Mat}(A^{p|q})_1$ an odd matrix,
\[ \text{Tr}(S) =
\text{tr}S_1 + \text{tr}S_4.
\]
Note the sign change.  Note also that the trace is commutative, meaning that for even matrices $A,B \in \text{Mat}(A^{p|q})_0$, we have the familiar formula
\[
\text{Tr}(AB) = \text{Tr}(BA).
\]
\end{remark}

\begin{definition}
\label{gldef}
Again let $M = A^{p|q}$, the free $A$-module generated by $p$ even
and $q$ odd variables.  Then $\text{GL}(A^{p|q})$ denotes the \textit{super general linear
group of automorphisms of $M$}; we may also use the notation
$\text{GL}_{p|q}(A) = \text{GL}(A^{p|q})$.
\end{definition}

\begin{remark}
If $M$ is an $A$-module, then $\GL(M)$ is defined as the group of automorphisms of $M$.  If $M = A^{p|q}$, then we write $\GL(M) = \GL_{p|q}(A)$ as above.
\end{remark}

Next we define the generalization of the determinant, called the \textit{Berezinian}, on elements of $\text{GL}(A^{p|q})$.

\begin{definition}  
Let $T \in \text{GL}(A^{p|q})$ have the standard block form (\ref{matrixT}) from
above. Then we formulate Ber:
\begin{equation}
\label{bereq}
\text{Ber}(T) = \det(T_1 - T_2T_4^{-1}T_3)\det(T_4)^{-1}
\end{equation}
where ``det" is the usual determinant.
\end{definition}

\begin{remark}
The first thing we
notice is that in the super category, we only define the Berezinian 
for \textit{invertible} transformations.  We immediately see 
that it is necessary
that the block $T_4$ be invertible for the formula (\ref{bereq}) to
make sense, however one can actually define the Berezinian on 
all matrices with \textit{only} the $T_4$ block invertible 
(i.e. the matrix itself may not be invertible, but the $T_4$ block is).  
There is a similar formulation of the Berezinian which requires that 
only the $T_1$ block be invertible:
$$
\text{Ber}(T) = \det(T_4 - T_3T_1^{-1}T_2)\det(T_1)^{-1}
$$
So we can actually define the Berezinian on all matrices with 
\textit{either} the $T_1$ \textit{or} the $T_4$ block invertible.  
Note that in the case where both blocks are invertible 
(i.e. when the matrix $T$ is invertible), both formulae of 
the Berezinian give the same answer.
\end{remark}

We leave the following proposition as an exercise.

\begin{proposition}
Let $T: A^{p|q} \longrightarrow A^{p|q}$ be a morphism with the
usual block form (\ref{matrixT}).  Then $T$ is invertible if and
only if $T_1$ and $T_4$ are invertible.
\end{proposition}

\begin{proposition}
The Berezinian is multiplicative: For $S,T \in \text{GL}(A^{p|q})$,
\[
\text{Ber}(ST) = \text{Ber}(S)\text{Ber}(T).
\]
\end{proposition}

\begin{proof}
We will only briefly sketch the proof here and leave the details
to the reader.  First note that any $T \in \text{GL}(A^{p|q})$ with block form (\ref{matrixT}) may be
written as the product of the following ``elementary matrices":

\begin{equation}
T_+ = \left( \begin{matrix} 1 & X\\ 0 & 1 \end{matrix} \right), \hspace{.1in} T_0 =  \left( \begin{matrix} Y_1 & 0\\ 0 & Y_2 \end{matrix} \right), \hspace{.1in} T_- =  \left( \begin{matrix} 1 & 0\\ Z & 1 \end{matrix} \right).
\end{equation}
If we equate $T=T_+T_0T_-$, we get a system of equations which lead to the solution
\begin{equation*}
\begin{array}{lcl}
X & = & T_2T_4^{-1},\\
Y_1 & = & T_1 - T_2T_4^{-1}T_3,\\
Y_2 & = & T_4,\\
Z & = & T_4^{-1}T_3.
\end{array}
\end{equation*}

It is also easy to verify that $\text{Ber}(ST)=\text{Ber}(S)\text{Ber}(T)$ for $S$ of type $\{T_+, T_0\}$ \textit{or} $T$ of type $\{T_-, T_0\}$.  The last case to verify is for
\[
S = \left( \begin{matrix} 1 & 0\\ Z & 1 \end{matrix} \right) \hspace{.1in} \text{and} \hspace{.1in} T = \left( \begin{matrix} 1 & X\\ 0 & 1 \end{matrix} \right).
\]
We may assume that both $X$ and $Z$ each have only one non-zero entry since the product of two matrices of type $T_+$ results in the sum of the upper right blocks, and likewise with the product of two type $T_-$ matrices.  Let $x_{ij}, z_{kl} \neq 0$.  Then
\[
ST = \left( \begin{matrix} 1 & X\\ Z & 1+ZX \end{matrix} \right)
\]
and $\text{Ber}(ST) = \det(1-X(1+ZX)^{-1}Z)\det(1+ZX)^{-1}$.  Because all the values within the determinants are either upper triangular or contain an entire column of zeros ($X,Z$ have at most one non-zero entry), the values $x_{ij},z_{kl}$ contribute to the determinant only when the product $ZX$ has its non-zero term on the diagonal, i.e. only when $i=j=k=l$.  But then $ \det(1-X(1+ZX)^{-1}Z) = 1+x_{ii}z_{ii}$, and it is clear that $\text{Ber}(ST) = 1$.  A direct calculation shows that $\text{Ber}(S) = \text{Ber}(T) = 1$.
\end{proof}

\begin{corollary}
The Berezinian is a
homomorphism
\[
\text{Ber}: \text{GL}(A^{p|q}) \longrightarrow \text{GL}_{1|0}(A) =
A_0^{\times}
\]
into the invertible elements of $A$.
\end{corollary}

\begin{proof}
This follows immediately from above proposition.
\end{proof}

The usual determinant on the general linear group $\text{GL}_n$ induces the trace on its
Lie algebra, namely the matrices $\text{M}_n$ (see Remark \ref{MatNotM}).  The extension to the
Berezinian gives
\[
\text{Ber}(I + \epsilon T) = 1 + \epsilon \text{Tr}(T)
\]
where $I$ is the $p|q \times p|q$ identity matrix (ones down the
diagonal, zeros elsewhere) and $\epsilon^2 = 0$.  An easy
calculation then exposes the super trace formula with the negative
sign.  This of course leads to the question of how the formula for
the Berezinian arises.  The answer lies in the SUSY-version of
integral forms on supermanifolds called
\textit{densities}.  In F.A. Berezin's pioneering work in
superanalysis, Berezin calculated the change of variables formula
for densities on isomorphic open submanifolds of $\R^{p|q}$ (\cite{Berezin}). This
lead to an extension of the Jacobian in ordinary differential
geometry; the Berezinian is so named after him.

We finish our summary of superlinear algebra by giving meaning to the
\textit{rank} of a endomorphism of $A^{p|q}$.
\begin{definition}
\label{rankDef}
Let $T \in \text{End}(A^{p|q})$. Then the \textit{rank} of
$T$, rank($T$), is the superdimension of the largest submatrix of
$T$ (obtained by removing columns and rows).
\end{definition}

\begin{proposition}
Again, let $T \in \text{End}(A^{p|q})$ with block form
(\ref{matrixT}).  Then $\text{rank}(T) =
\text{rank}(T_1)|\text{rank}(T_4)$.
\end{proposition}

\begin{proof}
Assume that $\text{rank}(T) = r|s$.  Then there is an invertible $r|s \times r|s$ submatrix of $T$ and it is clear that $r\leq \text{rank}(T_1), s\leq \text{rank}(T_4)$.  Conversely, if $\text{rank}(T_1) = r', \text{rank}(T_4) = s'$ it is also clear that there exists an invertible $r'|s' \times r'|s'$ submatrix of $T$.  Therefore we must have $r=r', s=s'$.
\end{proof}


\chapter{Supergeometry} \label{sgchap2}

In this chapter we discuss the foundations of supergeometric objects.  We
begin by defining 
the most basic object, the \textit{super ringed space} and build some basic
concepts from 
this definition.

\section{Superspaces}

\begin{definition}

As in ordinary algebraic geometry, a \textit{super ringed space} is a
topological 
space $|S|$ endowed with a sheaf of supercommuting rings which we denote by
$\shf_S$.  
Let $S$ denote the super ringed space $(|S|, \shf_S)$.

\end{definition}

\begin{definition}

A \textit{superspace} is a super ringed space $S$ with the property that the
stalk $\shf_{S,x}$ 
is a local ring for all $x \in |S|$.

\end{definition}

Given an open subset $U\subset |S|$, we get an induced \textit{subsuperspace}
given by 
restriction: $(U, \shf_S|_U)$.  For a closed superspace, we make the following definition:

\begin{definition}

Let $S$ be a superspace.  Then we say that $S'$ is a \textit{closed subsuperspace} of $S$ if\\

\noindent(i) $|S'| \subset |S|$ is a closed subset;\\

\noindent(ii) The structure sheaf on $S'$ is obtained by taking the quotient
of $\shf_S$ 
by a quasi-coherent sheaf of ideals $\mathcal{I}$ in $\shf_S$:
\[
\shf_{S'}(U) = \shf_S(U)/\mathcal{I}(U)
\]
for all open subsets $U$.

\end{definition}

Next we define a morphism of superspaces so that we can talk about the category of superspaces.

\begin{definition}

Let $S$ and $T$ be superspaces.  Then a morphism $S \longrightarrow T$ is a
continuous map 
$|\varphi|: |S| \longrightarrow |T|$ together with a sheaf map
$\varphi^*:\shf_T \longrightarrow 
\varphi_*\shf_S$ so that $\varphi^*_x(\mathfrak{m}_{|\varphi|(x)}) \subset
\mathfrak{m}_x$ where 
$\mathfrak{m}_x$ is the maximal ideal in $\shf_{S,x}$ and $\varphi^*_x$ is the
stalk map.  We denote the pair $(|\varphi|, \varphi^*)$ by $\varphi: S \longrightarrow T$.

\end{definition}

\begin{remark}

The sheaf map $\varphi^*:\shf_T \longrightarrow \varphi_*\shf_S$ corresponds
to the system 
of maps $\varphi^*|_U: \shf_T(U) \longrightarrow \shf_S(|\varphi|^{-1}(U))$
for all open sets
 $U \subset T$.  To ease notation, we also refer to the maps $\varphi^*|_U$ as $\varphi^*$.

\end{remark}

Essentially the condition $\varphi^*_x(\mathfrak{m}_{|\varphi|(x)}) \subset \mathfrak{m}_x$ means that the sheaf homomorphism is local.  Note also that $\varphi^*$ is a morphism of supersheaves, so it preserves parity.  The main point to make here is that the sheaf map must be specified along with the continuous topological map since sections are not necessarily genuine functions on the topological space as in ordinary differential geometry.  An arbitrary section cannot be viewed as a function because supercommutative rings have many nilpotent elements, and nilpotent sections are identically zero as functions on the underlying topological space.  Therefore we employ the methods of algebraic geometry to study such objects.  We will address this in more detail later.  Now we introduce two types of superspaces that we examine in detail in the forthcoming chapters: supermanifolds and superschemes.

\section{Supermanifolds}

A supermanifold is a specific type of ``smooth" superspace which we describe via a local model.  Because we always keep an eye on the physics literature from which supersymmetry springs, the supermanifolds of interest to us are the $C^{\infty}$-supermanifolds over $\R$.

Let $C_U^{\infty}$ be the sheaf of $C^{\infty}$-functions on the domain $U \subset \R^p$.  We define the \textit{superdomain} $U^{p|q}$ to be the super ringed space $(U, C_U^{\infty}[\theta^1, \ldots, \theta^q])$ where $C_U^{\infty}[\theta^1, \ldots, \theta^q]$ is the sheaf of supercommutative $\R$-algebras given by (for $V \subset U$ open)

\[
V \mapsto C_U^{\infty}(V)[\theta^1, \ldots, \theta^q].
\]

The $\theta^j$ are odd (anti-commuting) global sections which we restrict to $V$ .  Most immediately, the superspaces $\R^{p|q}$ are superdomains with sheaf $C_{\R^p}^{\infty}[\theta^1, \ldots, \theta^q]$.

\begin{definition}

A \textit{supermanifold} of dimension $p|q$ is a superspace which is locally isomorphic to $\R^{p|q}$.  Given any point $x \in |M|$, there exists a neighborhood $V\subset |M|$ of $x$ with $p$ even functions $(t^i)$ and $q$ odd functions $\theta^j$ on $V$ so that

\begin{equation}
\label{localM}
\shf_M|_V = \underbrace{C^{\infty}(t^1,\ldots, t^p)}_{C_M^{\infty}(V)}[\theta^1, \ldots, \theta^q].
\end{equation}

\end{definition}

Morphisms of supermanifolds are morphisms of the underlying superspaces.  For $M,N$ supermanifolds, a morphism $\varphi:M\longrightarrow N$ is a continuous map $|\varphi|:|M| \longrightarrow |N|$ together with a (local) morphism of sheaves of superalgebras $\varphi^*: \shf_N \longrightarrow \varphi_* \shf_M$.  Note that in the purely even case of ordinary $C^{\infty}$-manifolds, the above notion of a morphism agrees with the ordinary one.  We may now talk about the category of supermanifolds.  The difficulty in dealing with $C^{\infty}$-supermanifolds arises when one tries to think of ``points" or ``functions" in the traditional sense.  The ordinary points only account for the topological space and the underlying sheaf of ordinary $C^{\infty}$-functions, and one may truly only talk about the ``value" of a section $f \in \shf_M(U)$ for $U \subset |M|$ an open subset; the value of $f$ at $x\in U$ is the unique real number $c$ so that $f-c$ is not invertible in any neighborhood of $x$.  What this says is that we cannot reconstruct a section by knowing only its values at topological points.  Such sections are then not truly functions in the ordinary sence, however, now that we have clarified this point, we may adhere to the established notation and call such sections$f$ \textit{``functions on $U$"}.

\begin{remark}

Let $M$ be a supermanifold, $U$ an open subset in $|M|$, and $f$ a function on $U$.  If $\shf_M(U) = C^{\infty}(t^1, \ldots, t^p)[\theta^1, \ldots, \theta^q]$ as in (\ref{localM}), there exist even functions $f_I \in C^{\infty}(t)$ ($t=t^1, \ldots t^p)$) so that

\begin{equation}
\label{functionExpress}
f(t, \theta) = f_0(t) + \sum_i f_i(t) \theta^i + 
\sum_{i < j} f_{ij}(t)\theta^i \theta^j + \ldots = f_0(t) + \sum_{|I|=1}^q f_I(t)\theta^I
\end{equation}
where $I = \{i_1 < i_2 < \ldots < i_r \}_{r=1}^q$.

\end{remark}

Let us establish the following notation.  Let $M$ be a supermanifold, then we write the nilpotent sections as

\begin{equation}
J_M = \shf_{M,1}+\shf_{M,1}^2 = \langle \shf_{M, 1} \rangle_{\shf_M}.
\end{equation}

This is an ideal sheaf in $\shf_M$ and thus defines a natural subspace of $M$ we shall call $M_{\text{red}}$, or $\widetilde{M}$, where

\begin{equation}
\widetilde{M} = (|M|, \shf_M/J_M).
\end{equation}

Note that $\widetilde{M}$ is a completely even superspace, and hence lies in the ordinary category of ordinary $C^{\infty}$-manifolds,  i.e. it is locally isomorphic to $\R^p$.  The quotient map from $\shf_M \longrightarrow \shf_M/J_M$ defines the inclusion morphism $\widetilde{M} \hookrightarrow M$.  The subspace $\widetilde{M}$ should not be confused with the purely even superspace $(|M|, \shf_{M,0})$ which is \textit{not} a $C^{\infty}$-manifold since the structure sheaf still contains nilpotents.

\begin{observation}
\label{closedsubmflds}
Here we examine closed submanifolds in the super category.  Let $M$ be a supermanifold.  Then a submanifold of $M$ is a supermanifold $N$ together with a immersion, that is a map $i:N \longrightarrow M$ so that $i$ induces an imbedding of $\widetilde{N}$ onto a closed (locally closed) ordinary submanifold of $\widetilde{M}$ and $i^*_U : \shf_M(U) \longrightarrow \shf_N(i^{-1}(U))$ is surjective for all open $U \subset |M|$.

Locally, this means that we can find a system of coordinates $(t^1, \ldots, t^p, \theta^1, \ldots, \theta^q)$ in any open neighborhood of $M$ so that $N$ restricted to this neighborhood is described by the vanishing of some of the coordinates:
\[
t^1 = \ldots = t^r = \theta^1 = \ldots = \theta^s = 0.
\]

One can check that this explanation of submanifolds agrees with the definition of a closed sub superspace given earlier.

\end{observation}

\section{Superschemes}

A \textit{superscheme} is an object in the category of superspaces which generalizes the notion of a scheme.

\begin{definition}

A superspace $S = (|S|, \shf_S)$ is a superscheme if $(|S|, \shf_{S,0})$ is an ordinary scheme and $\shf_{S,1}$ is a quasi-coherent sheaf of $\shf_{S,0}$-modules.

\end{definition}

Because any non-trivial supercommutative ring has non-zero nilpotents, we need to redefine what we mean by a reduced superscheme.

\begin{definition}

We say that a superscheme $S$ is \textit{super reduced} if $\shf_S/J_S$ is reduced.  In other words, in a super reduced superscheme, we want that the odd sections generate all the nilpotents.

\end{definition}

\begin{example}

Let $\A^m$ be the ordinary affine space of dimension $m$ over $\C$ given with the Zariski topology.  On $\A^m$ we define the following sheaf $\shf_{\A^{m|n}}$ of superalgebras.  Given $U \subset \A^m$ open,

\begin{equation}
\label{Ascheme}
\shf(U) = \shf_{\A^n}(U)[\xi_1, \ldots, \xi_n]
\end{equation}

where $\shf_{\A^m}$ is the ordinary sheaf on $\A^m$ and the $\xi_1, \ldots, \xi_n$ are anti-commuting (or \textit{odd}) variables.  One may readily check that $(\A^m, \shf_{\A^{m|n}})$ is a reduced supercheme which we hereon denote by $\A^{m|n}$.

\end{example}

\begin{remark}

The superscheme $\A^{m|n}$ is more than reduced; it is a smooth superscheme.  The difference being the local splitting in (\ref{Ascheme}).  We do not further explore the notion of smoothness in these notes.

\end{remark}

Morphisms of superschemes are just morphisms of superspaces, so we may talk about the subcategory of superschemes.  The category of superschemes is larger than the category of schemes; any scheme is a superscheme if we take a trivial odd component in the structure sheaf.  We will complete our exposition of the category of superschemes in chapters \ref{sgchap5}-\ref{sgchap6}.

\section{T-Points}

The presence of odd coordinates steals some of the geometric intuition away from the language of supergeometry.  For instance, we cannot see an ``odd point" -- they are invisible both topologically and as classical functions on the underlying topological space.  We see the odd points only as sections of the structure sheaf.  To bring some of the intuition back, we turn to the functor of points approach from algebraic geometry.

\begin{definition}

Let $S$ and $T$ be superspaces.  Then a \textit{$T$-point} of $S$ is a morphism $T \longrightarrow S$.  We denote the set of all $T$-points by $S(T)$.  Equivalently, \[S(T) = \text{Hom}(T,S).\]

\end{definition}

Let us recall an important lemma.

\begin{lemma}

(Yoneda's Lemma)  There is a bijection from the set of morphisms $\varphi:M \longrightarrow N$ to the set of maps $\varphi_T: M(T) \longrightarrow N(T)$, functorial in $T$.

\end{lemma}

\begin{proof}

Given a map $\varphi: M \longrightarrow N$, for any morphism $t:T\longrightarrow M$, $\varphi \circ t$ is a morphism $T \longrightarrow N$.  Conversely, we attach to the system $(\varphi_T)$ the image of the identity map from $\varphi_M: M(M) \longrightarrow N(M)$.

\end{proof}

Yoneda's lemma allows us to replace a superspace $S$ with its set of $T$-points, $S(T)$.  We can now think of a superspace $S$ as a representable functor from the category of superspaces to the category of sets.  In fact, when constructing a superspace, it is often most convenient to construct the functor of points, then prove that the functor is \textit{representable}.  Let us give a couple examples of $T$-points from the category of supermanifolds.

\begin{example}

\noindent (i) If $T$ is just an ordinary topological point (i.e. $T = (\R^{0|0}, \R)$), then a $T$-point of $M$ is an ordinary topological point of $|M|$.\\

\noindent (ii) If $M=\R^{p|q}$, then a $T$-point of $M$ is a system of $p$ even and $q$ odd functions on $T$ by definition of a superspace morphism.  This is made more clear in chapter \ref{sgchap5} by Proposition \ref{chartTheorem}.  Thus $\R^{p|q}(T) = \shf_{T,0}^p \oplus \shf_{T,1}^q = (\shf_T^{p|q})_0$.\\

\end{example}

We already see the power of $T$-points in these two examples.  The first example ($T=\R^{0|0}$) gives us complete topological information while the second ($M=\R^{p|q}$) will allow us to talk about coordinates on supermanifolds.  We fully explore these topics in the next chapter.

In chapter 5 we give a criterion for the representability of functors from the category of superschemes to the category of sets.  As in the ordinary case, it turns out that representable functors must be \textit{local}, i.e. they should admit a cover by open affine subfunctors which glue together in some sense.


\chapter{$C^{\infty}$-supermanifolds} \label{sgchap3}

We have already described a supermanifold in chapter 2 as a superspace 
which is locally isomorphic to $\R^{p|q}$.
Recall also that given a supermanifold $M$, there is a surjection 
$\shf_M \longrightarrow \shf_M/J_M$ which corresponds to the natural 
inclusion $\widetilde{M} \hookrightarrow M$.  For local functions $f$ 
on $\shf_M$, this means $f \mapsto \wt{f} = f_0$ from the decomposition 
(\ref{functionExpress}).

\section{Charts}

Let us begin studying supermanifold morphisms in detail through the following example.

\begin{example}
\label{chartEx}
Consider the supermanifold $\R^{1|2}$ with a morphism $\varphi: \R^{1|2} \longrightarrow \R^{1|2}$.  
On $\R^{1|1}$ we have global coordinates $t, \theta^1, \theta^2$ and so we may express any function $f$ as in (\ref{functionExpress}):
\[
f=f(t, \theta^1, \theta^2) = f_0(t) + f_1(t)\theta^1 + f_2(t)\theta^2 + f_{12}(t)\theta^1\theta^2.
\]
Then $\wt{f} = f_0(t) \in C^{\infty}(\R)$ which sits as a function on the reduced $C^{\infty}$-manifold $\wt{\R^{1|2}} = \R$.  The morphism $\varphi$ is described by a continuous map $|\varphi|$ and a sheaf map $\varphi^*$. 

We first prescribe the global coordinates under $\varphi^*$:
\begin{equation}
\label{coorStar}
\begin{array}{lclcl}
t & \mapsto & t^* & = & t + \theta^1\theta^2\\
\theta^1 & \mapsto & \theta^{1*} & = & \theta^1\\
\theta^2 & \mapsto & \theta^{2*} & = & \theta^2.
\end{array}
\end{equation}

We claim that knowing $\varphi^*$ on only these global coordinates is enough to completely describe $\varphi$.  Indeed, we first see that $t \mapsto t + \theta^1\theta^2$ tells us that $|\varphi|$ is just the identity map.  Next, let $f \in C^{\infty}(t)[\theta^1, \theta^2]$ be as above.  Then $f \mapsto f^*$;
\begin{equation}
f^* = f(t^*, \theta^{1*}, \theta^{2*}) = f_0(t + \theta^1\theta^2) + f_1(t + \theta^1\theta^2)\theta^1 + f_2(t + \theta^1\theta^2)\theta^2 + f_{12}(t + \theta^1\theta^2)\theta^1\theta^2.
\end{equation}
And so we must only make sense of $f_I(t + \theta^1\theta^2)$.  The key is that we take a \textit{Taylor series expansion}; the series of course terminates thanks to the nilpotence of the the odd coordinates:
\begin{equation}
f_I(t + \theta^1\theta^2) = f_I(t) + \theta^1\theta^2 f_I'(t).
\end{equation}
It is easy to check that this in fact gives a homomorphism of superalgebras.  For $g,h \in C^{\infty}(\R)$, $(gh)^* = gh + \theta^1\theta^2 (gh)' = g^*h^*$.  The global sections are enough since the full sheaf map is given by restriction.  So in this example, it is enough to know $\varphi$ on only the coordinates.  In fact, the morphism induced by the equations (\ref{coorStar}) is unique via the above construction.
\end{example}

That a morphism is determined by local coordinates is true in general; we summarize this fact in the following \textit{Chart Theorem}.

\begin{theorem}
\label{chartTheorem}
(Chart) Let $U\subset \R^{p|q}$ be an open submanifold of $\R^{p|q}$ ($U = U^{p|q}$ is a superdomain).  There is a bijection between\\
\noindent (i) the set of morphisms $\varphi:M \longrightarrow U$ and\\
\noindent (ii) the set of systems of $p$ even functions $t^{i*}$ and $q$ odd functions $\theta^{j*}$ on $M$ so that $\wt{t^{i*}}(m) \in |U|$ for all $m \in |M|$.
\end{theorem}

\begin{proof} We sketch the proof of this well-known result here (for more details, see for example \cite{Leites}).  The key point is that given a system of $p$ even functions $t^{i*}$ and $q$ odd functions $\theta^{j*}$ on $M$, we can define a sheaf map.  As in the example (\ref{chartEx}) it is enough to define the sheaf map for $f \in C^{\infty}(U)$ since the expansion of an arbitrary section is linear in the odd coordinates over $C^{\infty}(U)$ and since we can restrict to an open $V \subset U$.

We define $\varphi$ formally by
\[
f = f(t^1, \ldots, t^p) \mapsto f^* = f(t^{1*}, \ldots, t^{p*}).
\]
We can write $t^{i*} = \wt{t^{i*}} + n^i$ where the $n^i$ are nilpotent, and we are set up to take a Taylor series expansion, just as in the above example:
\begin{equation}
\label{taylorSeries}
f(\wt{t^{1*}} + n^1, \ldots, \wt{t^{p*}} + n^p) := \sum_k \frac{\partial^k}{\partial t^k} f(\wt{t^{1*}}, \ldots, \wt{t^{p*}}) \frac{n^k}{k!}
\end{equation}
where $n^k$ is $k$-tuples of $\{n^i\}$.  This series terminates again thanks to the nilpotent $n_i$.  The $C^{\infty}$-functions $\wt{t^{i*}}$ completely determine the topological map $|\varphi|$.
\end{proof}

\begin{remark}
Note that because the expansion (\ref{taylorSeries}) involves an arbitrary number of derivatives, there is no way to make sense of $C^k$-supermanifolds.  We may, however, talk about the category of analytic (over $\R$ or $\C$) supermanifolds.  We refer to theorem (\ref{chartTheorem}) as the \textit{Chart Theorem} because it equates the definition of a supermanifold to giving an atlas of local charts.  These local charts glue together isomorphic copies of open subsets of $\R^{p|q}$.
\end{remark}

As in the category of superschemes, products exist in the category.  Let $M$ be a dimension $p|q$ supermanifold and $N$ be a dimension $r|s$ supermanifold, then we describe $M \times N$ by
\[
M \times N = (|M|\times|N|, \shf_{M\times N}).
\]
We define the sheaf $\shf_{M \times N}$ as follows.  For coordinate neighborhoods $U (= (x,\theta)) \subset |M|, V (= (t, \eta)) \subset |N|$, $\shf_{M\times N}(U \times V) = C^{\infty}(x,t)[\theta, \eta]$.  One must show that gluing conditions are satisfied, but this calculation mimics that in the ordinary category and is left for the reader.  So $M \times N$ is a $(p+r)|(q+s)$-dimensional supermanifold with $\widetilde{M \times N} = \wt{M} \times \wt{N}$.  As in the ordinary category, $\shf_{M \times N} \neq \shf_M \times \shf_N$; instead the we must take the completion of the tensor product to get an equality.

\begin{remark}
We cannot think of a supermanifold simply as a fiber space over an ordinary manifold.  Morphisms between supermanifolds mix both even and odd coordinates and thus for an open neighborhood $U$ of a supermanifold $M$, $C^{\infty}$ cannot be realized as a subsheaf of $\shf_M$; it follows that there is no natural morphism $M \longrightarrow \wt{M}$.  The symmetries of interest in these extensions of classical manifolds are those which place even and odd on the same level.  Such symmetries are called \textit{supersymmetries} and are at the foundation of the physical supersymmetry theory which aims to treat bosons and fermions on the same footing.
\end{remark}

\section{Vector Fields}

Many concepts and results from ordinary differential geometry extend naturally to the category of supermanifolds.  If we keep the categorical language we have developed, there is hardly any difference in fundamental differential geometry between the ordinary and the super categories.  For example, a vector bundle on a supermanifold $M$ is a locally free sheaf of (super)modules over $\shf_M$.  This leads to the notion of a tangent bundle on $M$, where we find super extensions of the inverse and implicit function theorem (see \cite{Leites}), and the local and global Frobenius theorem which we will prove in the next chapter.

\begin{definition}
\label{vfdef}
A \textit{vector field} $V$ on a supermanifold $M$ is an $\R$-linear derivation of $\shf_M$, i.e. it is a family of of derivations $\shf_M(U) \longrightarrow \shf_M(U)$ that is compatible with restrictions.
\end{definition}

The vector fields form a sheaf of modules over $\shf_M$, the \textit{tangent sheaf}  which we denote by $\text{Vec}_M$.  The sheaf $\text{Vec}_M$ is actually locally free over $\shf_M$ which we establish with the following lemma.  The lemma also helps us understand the local structure of a vector field.

\begin{lemma}
\label{distlemma}
Let $(t,\theta)$ be coordinates on some open subsupermanifold $U \subset \R^{p|q}$.  Then the $\shf_U$-module of $\R$-linear derivations of $\shf_U$ is a rank $p|q$ free sheaf over $\shf_U$ with basis $\{ \partial/\partial t^i, \partial/\partial \theta^j \}$ where
\begin{equation}
\frac{\partial}{\partial t^i} (f_I(t)\theta^I) = \frac{\partial f_I(t)}{\partial t^i}\theta^I, \hspace{.5in}
\frac{\partial}{\partial \theta^j} (f_I(t)\theta^j\theta^I) = f_I(t)\theta^I
\end{equation}
where $j \notin I$.
\end{lemma}

\begin{proof}
The proof is the same as in the classical case since the $\theta$-variables are polynomial (in fact, they are linear).
\end{proof}

Since $U \subset \R^{p|q}$ is the local model for any dimension $p|q$ supermanifold $M$, the lemma implies that $\text{Vec}_M$ is a vector bundle of rank $p|q$.  If $V$ is a vector field on $U$, then in a coordinate neighborhood $U' \subset U$ with coordinates $(t, \theta)$, there exist functions $f_i, g_j$ on $U'$, so that $V$ has the unique expression
\begin{equation}
\label{vfexpress}
V|_{U'} = \sum_{i=1}^p f_i(t, \theta)\frac{\partial}{\partial t^i} + \sum_{j=1}^q g_j(t, \theta) \frac{\partial}{\partial \theta^j}.
\end{equation}

We similarly define the tangent space at a point $m$ of $M$, which we denote $T_m(M)$.  We think of tangent vectors as $\R$-linear derivations $\shf_m \longrightarrow \R$ of the stalk at $m$; we may think of a tangent vector $v\in T_m(M)$ as a vector field on $U$, a neighborhood of $m$, composed with evaluation at $m$.  If the open subset $U$ from definition \ref{vfdef} is a coordinate neighborhood around $m$, the vector $v$ takes the expression
\begin{equation}
v= \sum a_i \frac{\partial}{\partial x^i}|_m + b_j \frac{\partial}{\partial \theta^j}|_m
\end{equation}
for $a_i, b_j \in \R$.

For $M$ and $N$ supermanifolds, we can extend a vector on $M$ to a $\shf_N$-linear derivation on $M \times N$, and likewise we may trivially treat any vector field on $M$ as a vector field on $M \times N$.  We will call these extensions \textit{extended vectors} and \textit{extended vector fields} respectively.

\begin{definition}
Let $v$ be a tangent vector of $M$ at $m$, $U_m \subset |M|$ an open neighborhood of $m$.  We view $v$ as a derivation
$\shf_M(U_m) \longrightarrow \R$ and identify $\shf_N$ with $\shf_{\R \times N}$. Then $v$ extends uniquely to a \textit{$\shf_N$-linear derivation}:
\[
\begin{array}{rrcl}
v_N : & \shf_{M\times N}(U_m \times V) & \longrightarrow & \shf_N(V)\\
& \searrow & & \nearrow\\
& & \shf_{\R \times N}(\R \times V) & \\
\end{array}
\]
for any open $V \subset |N|$ (this is easily seen locally by using coordinates, and then by
patching using local uniqueness) so that
\begin{equation}
v_N (a \otimes b)=v(a)b
\end{equation}
where $a$ and $b$ are local functions of $M$ and $N$ respectively.
\end{definition}




One may similarly ``extend" vector fields: let $V$ be a vector field on $M$.  Then we extend $V$ to a derivation $(V \tensor id)$ on $M \times N$ by forcing $V$ to act trivially on $N$.  If $(t, \theta)$ and $(x, \xi)$ are local coordinates on $M$ and $N$ respectively, $V$ has the coordinate expression as in (\ref{vfexpress}).  Then the extension $(V \tensor id)$ has the same coordinate expression on $M \times N$ described by coordinates $(t,x,\theta,\xi)$, i.e. it is identically zero on $(x, \xi)$.  Again the extension is unique by patching using local uniqueness.

\section{Differential Calculus}

In this section we discuss the notion of differential
of a morphism of supermanifolds. In this context
the theory of supermanifolds resembles very closely the classical theory.
For completeness, we give a summary of the well known results, 
sketching only briefly the proofs or leaving them as exercises.

\medskip

\begin{definition}
Let $\al:M \lra N$ be a morphism of supermanifolds. We define
\textit{differential} of $\al$ at a topological point $m \in |M|$ the
map $(d \al)_m:T_mM \lra T_{|\al|(m)}N$ given by:
$$
(d \al)_m(X)(f)=X(\al^*_m(f)), \qquad \al^*_m:\cO_{N,|\al|(m)} 
\lra \cO_{M,m},
$$
where $X \in T_{|\al|(m)}N \cong \Der(\cO_{M,m}, \R)$, as
we have seen in the previous section.
\end{definition}

In local coordinates one can readily check that $(d \al)_m$
has the usual jacobian expression. In fact, let's choose suitable
open submanifolds $U \subset M$ and $V \subset N$ 
such that  $m \in |U|$ and $|\al|(m) \in |V|$, 
homeomorphic respectively to open submanifolds in
$\R^{r|s}$ and $\R^{u|v}$, and let $(t^i,\theta^j)$ be local
coordinates in $U$.
We have that:
$$
\al(t^i,\theta^j)=(f^k, \phi^l) \subset V. 
$$
Then 
$$
(d \al)_m= {\begin{pmatrix}
{\partial f^k \over \partial t^i} & {\partial \phi^l \over \partial t^i} \\
{\partial f^k \over \partial \theta^j} & {\partial \phi^l \over \partial \theta^j}
\end{pmatrix}}_m,
$$
where the subscript $m$ means evaluation at $m$.

\medskip

As an example let's compute the differential of the morphism 
described in example \ref{chartEx}.

\begin{example} Let $\al:\R^{1|2} \lra \R^{1|2}$ be the 
morphism given locally (and globally) by:
$$
\al(t^1, \theta^1, \theta^2)=(t^1+\theta^1\theta^2, \theta^1, \theta^2)
$$ 
Then the differential at a generic topological point $m=(t_0, 0,0)$
is:
$$
(d\al)_m={\begin{pmatrix}
1 & 0 & 0 \\
\theta^2 & 1 & 0 \\
-\theta^1 & 0 & 1 
\end{pmatrix}}_m=\begin{pmatrix} 1 & 0 & 0 \\
0 & 1 & 0 \\ 0 & 0 & 1 \end{pmatrix}
$$
\end{example}

As in the ordinary theory there are classes of morphisms
that play a key role: immersions, submersions and diffeomorphisms.

\begin{definition} 
Let $\al:M \lra N$ be a supermanifold morphism and 
$\widetilde{\al}: \widetilde{M} \lra \widetilde{N}$ the underlying 
classical morphism on the
reduced spaces. $\al$ is an \textit{immersion} at $m \in |M|$ if 
$\widetilde{\al}$ is an
ordinary immersion at $m$ and $(d\al)_m$ is injective.
Likewise $\al$ is a \textit{submersion} at $m$ if 
$\widetilde{\al}$ is a
submersion at $m$ and $(d\al)_m$ is surjective. 
Finally $\al$ is a \textit{diffeomorphism}
at $m$ if it is a submersion and an immersion. When we say
$\al$ is a immersion (resp. submersion or diffeomorphism) we mean
$\al$ is such at all points of $|M|$.
\end{definition}

As in the classical setting submersions and immersions have the
usual local models.

\medskip

\begin{proposition} Let $\al:M \lra N$ be a supermanifold morphism.
Let $m \in |M|$ and let $U \subset M$ and $V \subset N$ two suitable open sets
homeomorphic respectively to open sets in $\R^{r|s}$ and $\R^{u|v}$,
$m \in |U|$, $\al(m) \in |V|$.

\noindent
1. If $\al$ is an immersion at $m$, after suitable
changes of coordinate in $U$ and $V$, locally we have that
$$
\al(t^i,\theta^k)=(t^i,0,\theta^k,0)
$$

\noindent
2. If $\al$ is an submersion at $m$, after suitable
changes of coordinate in $U$ and $V$, locally we have that
$$
\al(t^i,s^j,\theta^k,\sigma^l)=(t^i,\theta^k)
$$
\end{proposition}

The proof of this result is essentially the same as in the
ordinary setting and can be found for example in \cite{VSV2} chapter 4.

\begin{remark}
A \textit{ closed } sub supermanifold $N$ of $M$ can be
equivalently defined as a supermanifold
such that $\widetilde{N}$ is a closed supermanifold of 
$\widetilde{M}$ and
$N \subset M$ is an immersion. We leave as an exercise to the
reader the check that this definition is equivalent to the one
given in observation \ref{closedsubmflds} in chapter \ref{sgchap2}.
\end{remark}

\medskip

The Submersion Theorem in the supercontext plays an important role
in proving that certain closed subset of a supermanifold
admit a supermanifold structure.
 
\begin{theorem} Submersion Theorem. \label{submersionthm}
Let $f:M \lra N$ be a submersion at $n \in |N|$, and let $|P|=|g|^{-1}(n)$. 
Then $|P|$ admits a supermanifold structure, i. e. 
there exists a supermanifold $P=(|P|, \cO_P)$ where $\cO_P=\cO_M|_{|P|}$.
Moreover:
$$
\dim P= \dim M - \dim N
$$
\end{theorem}

\begin{proof} (Sketch).
Locally we can define 
for $p \in |P|$,
$\cO_{P, p}=\cO_{M,p}/f^*(I_n)$, where $I_n$ is the ideal
in $\cO_{N,n}$ of elements vanishing at $n$.
For $W \in |P|$ open, sections in $\cO_P(W)$ are defined as maps
$$
\begin{array}{cccc}
u:W & \lra & \coprod_{p \in W} \cO_{P,p} \\
q & \mapsto & u(q) \in \cO_{P,q}
\end{array}
$$
This gives $P$ the structure of a superspace. By the previous
proposition we know that locally $f(t,s,\theta, \sigma)=(t,\theta)$. Hence 
locally we have coordinates $(s, \sigma)$ and
$P$ is a supermanifold of the prescribed dimension.
\end{proof}

\begin{example} \label{sl11}
Let $X$ be the open submanifold of $\R^{2|2}$ whose topological
space consists of the points of plane $\R^2$ with the coordinate axis, 
$x=0$, $y=0$ removed. 
Let $f:X \lra \R^{1|0}$ be the morphism
$\al(x,y,\xi, \eta)=y^{-1}(x-\xi y^{-1} \eta)$.
One can check that at the topological point $1 \in \R$, $f$ 
is a submersion, hence $|P|=|f|^{-1}(1)$ admits a sub supermanifold 
structure.
If we identify $\R^{2|2}$ with the $1|1 \times 1|1$ supermatrices 
$P$ corresponds to the supermatrices with Berezinian equal to $1$ and
we denote $P$ with $SL(1|1)$, the \textit{special linear supergroup}.
\end{example}

\medskip

We next turn our attention to a special type of supermanifold, 
a \textit{super Lie group}.

\section{Super Lie Groups}

A Lie group is a group object in the category of manifolds. Likewise a super Lie group is a group object in the category of supermanifolds. This means that there are appropriate morphisms which correspond to the group operations: product $\mu: G \times G 
\longrightarrow G$, unit $e: 1 \longrightarrow G$ (where $1 \in |G|$ may be equated to $\R^{0|0}$, a single topological point), and inverse $i: G 
\longrightarrow G$ so that the necessary diagrams commute (these are, in fact, the same diagrams as in the ordinary setting).  We may of course interpret all these maps and diagrams in the language of $T$-points, which gives us (for any supermanifold $T$) morphisms $\mu_T: G(T) \times G(T) \longrightarrow G(T)$, etc. that obey the same commutative diagrams.  In other words, Yoneda's Lemma says that the set $G(T)$ is in fact a group for all $T$.  This leads us to our working definition of a super Lie group.

\begin{definition}
A supermanifold $G$ is a \textit{super Lie group} if for any supermanifold $T$, $G(T)$ is a group, and for any supermanifold $S$ and morphism $T\longrightarrow S$, the corresponding map $G(S) \longrightarrow G(T)$ is a group homomorphism.
\end{definition}
In other words, $T \mapsto G(T)$ is a functor into the category of groups.
\begin{example}
\label{r11ex}
Let us consider the super Lie group $\R^{1|1}$ through the symbolic language of $T$-points.
The product morphism $\mu: \R^{1|1} \times \R^{1|1} \longrightarrow \R^{1|1}$ is given by
\begin{equation}
\label{grouplaw}
(t, \theta) \cdot (t', \theta') = (t+t' + \theta \theta', \theta + \theta')
\end{equation}
where the coordinates $(t, \theta)$ and $(t', \theta')$ represent two distinct $T$-points for some supermanifold $T$.  It is then clear by the formula (\ref{grouplaw}) that the group axioms inverse, identity, and associativity are satisfied.
\end{example}

Also in the language of $T$-points, the definition given above is equivalent to saying that a super Lie group is a functor from the category of supermanifolds to the category of groups which is representable.  In this vein, let us further examine the $\GL_{p|q}$ example.

\begin{example}
Let's first construct the supermanifold $\GL_{p|q}$.  The reduced space $\widetilde{\GL}_{p|q} = \GL_p \times \GL_q$ is an open subset of $\R^{p^2} \times \R^{q^2}$.  We build the sheaf on $\GL_{p|q}$ from the smooth functions on $\GL_p \times \GL_q$ and the restriction of the global odd coordinates $\theta^1, \ldots, \theta^{2pq}$ on $\R^{p^2 + q^2|2pq}$, i.e. for an open $U \subset \GL_p \times \GL_q$,
\[
\shf_{\GL_{p|q}}(U) = (\GL_p \times \GL_q)(U) \tensor [\theta^1, \ldots, \theta^{2pq}]|_U.
\]
Now we can examine the $T$-points of $\GL_{p|q}$.
Let $T$ be a supermanifold, $t \in \GL_{p|q}(T)$ a $T$-point, then $t: T \longrightarrow \GL_{p|q}$ is a morphism.  By proposition \ref{chartTheorem}, $t$ is equivalent to giving the pullbacks of coordinates.  By taking into account the determinant identities which must be satisfied, we see that $t$ is then equivalent to an invertible matrix with coefficients in $\shf_T$, and so $t$ corresponds to an automorphism of $\shf_T^{p|q}$.  Thus $\GL_{p|q}(T)$ is the group of automorphisms of $\shf_T^{p|q}$.
\end{example}

\begin{example}
\label{sl}
Let us consider another example of a super Lie group, $\SL_{p|q}$.
We define $\SL_{p|q}$ in a way which mimics the classical construction.  For each supermanifold $T$, the Berezinian gives a morphism from the $T$-points of $\GL_{p|q}$ to the $T$-points of $\GL_{1|0}$:
\[
\text{Ber}_T: \GL_{p|q}(T) \longrightarrow \GL_{1|0}(T).
\]
The super special linear group $\SL_{p|q}$ is the kernel of $\text{Ber}_T$.

Using a similar argument as in example \ref{sl11} one can show that
the functor $\SL_{p|q}$ is the functor
of points of a super Lie group, closed sub supermanifold of $\GL_{p|q}$.
In fact $|\SL_{p|q}|=|\text{Ber}|^{-1}(1)$, where $\text{Ber}$ is the
map $\text{Ber}: \GL_{p|q} \lra \GL_{1|0}$ between supermanifold corresponding
by Yoneda's lemma to the family of maps $\text{Ber}_T$ given above.
\end{example}

\begin{example}
\label{osp}
In our last example we extend the classical orthogonal group to the super category.
Let $\Phi$ be an even nondegenerate bilinear form on $\R^{p|2q}$ with values in $\R^{1|0}$.  The form $\Phi$ is equivalent to giving nondegenerate symmetric bilinear form on $\R^p$ and a nondegenerate alternating form on $\R^{2q}$.  Then for any supermanifold $T$, define $\OSp_{p|2q}(T)$ as the subgroup of $\GL_{p|2q}(T)$ which preserves $\Phi$.
\end{example}

\begin{remark}
A word of caution.  In the above two examples, we only give 
$\SL_{p|q}$ and $\OSp_{p|2q}$ in terms of their $T$-points.  
It is clear that each is a functor from supermanifolds to groups.  
However, it is not clear without a further argument,
that the functors defined above are \textit{representable}.  
\end{remark}

\section{Left Invariant Vector Fields}

In the remainder of this chapter, we discuss left invariant vector fields 
on a super Lie group, then examine the infinitesimal interpretation of a 
super Lie group acting on a supermanifold, which will be most relevant when 
we examine the super Lie group/algebra, super Lie subgroup/subalgebra pairing.

\medskip

Let $G$ be a super Lie group with group law $\mu: G \times G \rightarrow G$.  Via $T$-points, we
can symbolically understand this group law as $(x, \xi)
\cdot (x', \xi') = (t, \theta)$ where
$t = t(x, x', \xi, \xi')$   
are even functions and $\theta = \theta(x, x', \xi, \xi')$
are odd functions.  Again, all this really says is that $\mu^*(t^i) = t^{i*}
= t_i(x, x', \xi, \xi')$ for some even section
$t_i$ of $\shf_{G \times G}$ and $\mu^*(\theta^j) = \theta^{j*} = \theta_j(x, x', \xi, \xi')$ for some odd section.

Recall classically that for an ordinary Lie group $H$, we could define a map $\ell_h$, ``left multiplication by $h$"($h \in H$):
\begin{equation}
H \stackrel{\ell_h}{\longrightarrow} H;\hspace{.5in}
a \mapsto ha
\end{equation}
(for $a \in H$).  The differential of this map gives
\begin{equation}
T_e(H) \stackrel{d\ell_h}{\longrightarrow} T_h(H)
\end{equation}
and for a vector field $X$ on $H$, we say that $X$ is \textit{left invariant} if
\begin{equation}
\label{classicLeftInvariance}
d\ell_h \cdot X = X \cdot \ell_h.
\end{equation}
We interpret this in the super category by saying that a left invariant vector field on $G$ is invariant with respect to the group law $\mu^*$ ``on the left".  What this amounts to in making a formal definition is that we replace the ordinary group law $\mu$ with the anti-group law $\iota$ given (via $T$-points) by:
\[
\iota(g,g') = \mu(g',g) = g' \cdot g.
\]


Since $V$ is a vector field, $V|_U :\shf_G(U) \rightarrow \shf_G(U)$ is a derivation for all open $U \subset |G|$, the expression $\iota^* \circ V$ makes sense.  Now in the spirit of (\ref{classicLeftInvariance}) we need to understand
``$V \circ \iota^*$.''  We trivially extend the derivation $V$ to $\shf_{G \times G}$, and the expression $(V \tensor id) \circ \iota^*$ is formal.  We can now make a definition.

\begin{definition}
If $V$ is a vector field on the super Lie group $G$, we say that
$V$ is \text{left invariant} if $(V \tensor id)\iota^* = \iota^* V$.
\end{definition}

As in the classical theory, we have the following theorem.
\begin{theorem}
\label{bijection}
There is a bijection between\\
i. left-invariant vector fields on $G$ and\\
ii. $T_e(G)$.
\end{theorem}

Before we prove theorem (\ref{bijection}), let us first establish some useful notation and a technical lemma.  Recall that for a supermanifold $X$ and any $v \in T_xX$, we have the extended $\shf_T$-linear derivation $v_T$ for any supermanifold $T$.  Moreover, let $\varphi:X \longrightarrow Y$ be a morphism of supermanifolds.  Then $\varphi$ induces the morphism $\varphi \times id_T: X \times T \longrightarrow Y \times T$ and we denote the pullback by
\begin{equation}
(\varphi \times id_T)^* = \varphi^* \tensor id_T.
\end{equation}
Similarly we can define $id_T \tensor \varphi^*$.


\begin{proof} (Theorem \ref{bijection})\\
\noindent Since $G$ is a super Lie group, there is a map $\text{id}:\{e\} \rightarrow
G$ which gives $\epsilon: \shf_G \rightarrow k$, ``evaluation at $e$".
If $V$ is a left invariant vector field on $G$, then $\epsilon V = v$ is
a tangent vector at the origin of $G$.  We claim that this $v$ in
fact determines $V$:
\[
V = (v_G)\iota^*.
\]
Let us first show that given any tangent vector $v$,
$(v_G)\iota^*$ is a left invariant vector field on $G$.

Heuristically we are doing the same thing as in the classical setting;
we are infinitesimally pushing the vector $v$ with the group law.  It is clear that $(v_G)\iota^*$ is
locally a derivation on $\shf_G$; we next show it is
left invariant, i.e. we must
show that
\begin{equation}
\label{v is l-i}
((v_G)\iota^* \tensor id_G)\iota^* = \iota^* (v_G)\iota^*.
\end{equation}
A direct check on local coordinates (one can always choose coordinates of the form $a \tensor b$ on $G \times G$) shows that
\[
((v_G)\iota^* \tensor id_G) = v_{G \times G}(\iota^* \tensor id_G).
\]
But $v_{G \times G}(\iota^* \tensor id_G) = v_{G \times G}(id_G \tensor \iota^*)$ by the coassociativity of $G$, and another direct check shows that
\[
v_{G \times G}(id_G \tensor \iota^*) = \iota^*(v_G)
\]
Hence the claimed equality (\ref{v is l-i}).
  
The only item left to show is that $V = (v_G)\iota^*$.  Note that we have the equality $id_G = (\epsilon \tensor id_G)\iota^*$ from the ``identity" group axiom.  Then $V = (\epsilon \tensor id_G)\iota^*V = (\epsilon \tensor id_G)(V \tensor id_G)\iota^*$ by left-invariance of $V$, but the right hand side of the last equality is precisely $(v_G)\iota^*$ by evaluation on local coordinates.
\end{proof}

\begin{remark}
A \textit{right invariant} vector field is similarly defined; we need only replace $\iota$ by $\mu$ in the above definitions and theorems.  There is a natural anti-homomorphism from left invariant vector fields to right invariant vector fields induced by the inverse map $i: G \longrightarrow G$.
\end{remark}

The left invariant vector fields are a subsuper Lie algebra of $\text{Vec}_G$ which we denote by
\[
\mathfrak{g} = \{ V \in \text{Vec}_G \hspace{.1in} | \hspace{.1in} (V \tensor id)\iota^* = \iota^*V \}.
\]
Since the bracket of left invariant vector fields is left invariant, in fact, $\mathfrak{g}$ is the super Lie algebra associated to the super Lie group $G$, and we write $\mathfrak{g} = \text{Lie}(G)$ as usual.

\begin{example}
\label{why not iota}
We will calculate the left invariant vector fields on $\R^{1|1}$ with the group law from example (\ref{r11ex})
\begin{equation}
(t,\theta)\cdot_{\mu} (t', \theta') = (t+t'+\theta\theta', \theta + \theta').
\end{equation}
From theorem (\ref{bijection}), we know that the Lie algebra of left invariant vector fields can be extracted from $T_eG = \text{span}\{\frac{\partial}{\partial t}|_e, \frac{\partial}{\partial \theta}|_e \}$.  We use the identity $V = (v_G)\iota^*$ from the proof of Theorem \ref{bijection} to calculate the corresponding left invariant vector fields:
\begin{equation}
\label{r11livf}
(\frac{\partial}{\partial t}|_e)_G \circ \iota^*, \hspace{.3in} (\frac{\partial}{\partial \theta}|_e)_G \circ \iota^*.
\end{equation}
To get coordinate representations of (\ref{r11livf}), we apply them to coordinates $(t, \theta)$:
\begin{equation}
\label{iota 1}
\begin{array}{lcl}
(\frac{\partial}{\partial t}|_e)_G \circ \iota^*(t) & = & (\frac{\partial}{\partial t}|_e)_G (t' + t +\theta'\theta) = 1\\
(\frac{\partial}{\partial t}|_e)_G \circ \iota^*(\theta) & = & (\frac{\partial}{\partial t}|_e)_G (\theta' + \theta) = 0;
\end{array}
\end{equation}
\begin{equation}
\label{iota 2}
\begin{array}{lcl}
(\frac{\partial}{\partial \theta}|_e)_G \circ \iota^*(t) & = & (\frac{\partial}{\partial \theta}|_e)_G (t' + t +\theta'\theta) = -\theta'\\
(\frac{\partial}{\partial \theta}|_e)_G \circ \iota^*(\theta) & = & (\frac{\partial}{\partial \theta}|_e)_G (\theta' + \theta) = 1.
\end{array}
\end{equation}
Thus the left invariant vector fields on $(\R^{1|1}, \mu)$ are
\begin{equation}
\label{thelivf}
\frac{\partial}{\partial t}, \hspace{.3in} -\theta\frac{\partial}{\partial t} + \frac{\partial}{\partial \theta}.
\end{equation}
A quick check using the definition shows that (\ref{thelivf}) are in fact left invariant.
\end{example}

\section{Infinitesimal Action}

In this section we discuss the infinitesimal interpretation of a super Lie group acting on a supermanifold.  We later use the results from this section in chapter 4 to build the super Lie group/algebra super Lie subgroup/subalgebra correspondence.

Let $G$ be a super Lie group, $M$ a supermanifold, and
\[
\varphi : G\times M\longrightarrow M
\]
be a morphism.  If $v\in T_eG$ ($e \in G$ the identity), then the composition
\begin{equation}
\label{newVF}
\shf_M(U) \stackrel{\varphi^*}{\longrightarrow} \shf_{G\times M}(U_e \times U)\stackrel{v_M}{\longrightarrow} \shf_M(U)
\end{equation}
is then a derivation of $\shf_M(U)$ for any open $U \subset |M|$ ($U_e \subset |G|$ is some open neighborhood of the identity $e$).  The Leibniz property can be verified directly by calculating on local coordinates.  Then the composition $v_M \circ \varphi^*$ defines a vector
field on $M$ which we denote by $V_M(v,\varphi)$:
\begin{equation}
\label{two}
V_M(v,\varphi)(f)=v_M(\varphi^*(f))
\end{equation}
for $f$ local function on $M$.  It is clear that as $v$ varies we get a map of super vector
spaces from $T_x(X)$ into the super vector space of vector fields
on $M$.  Let $S$ be another supermanifold and consider the morphism
\[
\varphi \times id_S : G\times M\times S\longrightarrow M\times S
\]
and we see that
\begin{equation}
\label{three} v_{M\times S}(\varphi^* (f)\otimes s)=v_M(\varphi^*(f))\otimes s
\end{equation}
for $f$ again a function on $M$ and $s$ a function on $S$.  We thus obtain the equality:
\begin{equation}
\label{four} V_{M\times S}(v,\varphi \times id_S)= V_M(v,\varphi)\otimes id_S.
\end{equation}
Note that it is enough to verify (\ref{four}) on sections of the form $f \tensor s$.  The equation establishes an equality of vector fields, and so it is enough to check it on coordinates.  We can always write coordinates in the form $f \tensor s$, and so the calculation (\ref{three}) is enough.

\begin{example}
\label{iota and mu}
Let $M=G$ and let the map
$\varphi= \iota$ be the anti-group law.
Then 
\[
V_G(v, \iota)=:V^\ell
\]
is the unique left invariant vector field on $G$ which defines the
tangent vector $v$ at $e$.  If we take
$\varphi=\mu$ where $\mu(gg')=gg'$ is the ordinary group law, then
\[
V_G(v, \mu)=V^r
\]
where $V^r$ is the unique right invariant vector field defining
the tangent vector $v$ at $e$. We know that $v\longmapsto V^\ell$
is a linear isomorphism of $T_e(G)$ with $\lie(G)$ and one can check that
$V^\ell \longmapsto V^r$ is an anti-isomorphism of super Lie
algebras.
\end{example}

\begin{definition}
\label{rhoex}
Let $\varphi=\sigma$ be an action of $G$ on $M$
\[
\sigma  : G\times M\longrightarrow M.
\]
We define a linear map $\rho$ by
\begin{equation}
\label{five}
\rho (v)=_{def}V_M(v,\sigma), \quad \rho(v)(f) = v_M(\sigma^*(f))\qquad
\end{equation}
for $f$ a section on $\shf_M$.
\end{definition}
In fact, the next theorem asserts that the association
\[
V^\ell\longmapsto \rho (v)
\]
is a linear map from $\lie(G)$ to $\text{Vec}_M$.

The definition of action of $G$ on $M$ gives rise to a commutative
diagram
\begin{equation}
\label{action}
\begin{matrix} G\times G\times M & \stackrel{\mu \times id_M}{\longrightarrow} & G\times M\\
_{id_G \times \sigma} \downarrow & & \downarrow_{\sigma}\\
G\times M & \stackrel{\sigma}{\longrightarrow} & M.
\end{matrix}
\end{equation}

\begin{theorem}
The map $\rho$ (\ref{rhoex}) is an antimorphism of super Lie algebras $\mathrm{Lie}(G) \longrightarrow \mathrm{Vec}_M$. It
moreover satisfies the property
\begin{equation}
\label{six}
(V^r\otimes id_M)(\sigma^* f)=\sigma^* (\rho(v)f)
\end{equation}
for $v \in T_e(G)$, $V^r$ its corresponding right-invariant vector field,  and $f$ a function on $M$.
\end{theorem}

\begin{proof}
It is enough to prove (\ref{six}) to prove both assertions.
Indeed, suppose we have proved (\ref{six}). Then we see that the
image of $\shf_M$ under $\sigma^*$ is stable under all the
vector fields $V^r\otimes id_M$, and that $V^r\otimes id_M$ and
$\rho(v)$ are $\sigma$-related. Moreover, as $f$ varies, $\sigma^*(f)$ is surjective onto all sections of $\shf_M$ by restriction since $G$ is a super Lie group (i.e. $G$ contains an identity element which acts trivially on $M$).  It is then immediate that
\[
V^r\longmapsto \rho(v)
\]
is a morphism of super Lie algebras. Hence $\rho$ is an
antimorphism of $\lie(G)$ into $\text{Vec}_M$.

It thus remains to prove (\ref{six}). It will come as a
consequence of the relation from the commutative diagram (\ref{action}) that
\begin{equation}
\label{seven} (\mu \times id_M)^* (\sigma^* (f)) =( id_G \times
\sigma)^* (\sigma^*(f)),
\end{equation}
and the equality we seek will come by evaluating $v_{G\times M}$ on both
sides of (\ref{seven}).

By (\ref{four}), $V_{G\times M}(v,\mu \times
id_M)=V_G(v,\mu)\otimes id_M$. Hence
\begin{equation}
\label{eight} \begin{array}{rcl} v_{G\times M}((\mu \times id_M)^*
(\sigma^* (f))) & = & V_{G\times M}(v,\mu \times
id_M)(\sigma^*(f))\\ & = & (V_G(v,\mu)\otimes id_M)(\sigma^*(f))\\
& = & (V^r\otimes id_M)(\sigma^*(f)).
\end{array}\end{equation}

We shall next evaluate $v_{G\times M}$ on the right side of
(\ref{seven}). Now for $u\in \shf_{G \times M}$,
\[
v_{G\times M}((id_G \times \sigma )^*(u))= V_{G\times M}(v,
id_G\times \sigma)(u).
\]
Let $Z$ denote the vector field $V_{G\times M}(v,id_G \times
\sigma)$ on $G\times M$ for brevity. Let $a\in \shf_G$ and $b\in \shf_M$; we get
\[
Z(a\otimes b)=v_{G\times M}(a\otimes \sigma^*(b))=
v(a)\sigma^*(b)=\sigma^*(v(a)b).
\]
On the other hand, $v_M(a\otimes b)= v(a)b$, so we may rewrite
the last equation as
\[
Z(a\otimes b)=\sigma^*(v_M(a\otimes b)).
\]
But $Z$ and $\sigma^* \circ v_M$ are both derivations of $\shf_{
G\times M}$, and hence are vector fields on $G\times M$.  Then
the above relation shows that they must be identical as it is enough to check the equality of vector fields on coordinates, and we may always find coordinates of the form $a \tensor b$. Hence
\[
Z(u)=\sigma^*(v_M(u))
\]
for $u\in \shf_{G\times M}$.  If we take $u=\sigma^*(f)$ again for
$f$ a local function on $M$, by definition we have $v_M(\sigma^*(f))=\rho(v)(f)$ and so
the right side is equal to $\sigma^*(\rho(v)f)$. The left side is
equal to $v_{G\times M}((id_G\times \sigma)^*(f))$. Hence
\begin{equation}
\label{nine} v_{G\times M}((id_G\times \sigma)^*(f))=
\sigma^*(\rho(v)f).
\end{equation}
The equations (\ref{eight}) and (\ref{nine}) give us our result.
\end{proof}

\begin{corollary}
The anti-morphism $\rho$ extends to an associative algebra
anti-morphism (which we also call $\rho$),
\[
\rho: \mathcal{U}(\mathrm{Lie}(G)) \longrightarrow
\mathcal{U}(\mathrm{Vec}_M).
\]
\end{corollary}

\begin{proof}
We use the universal property of the universal enveloping algebra
and extend the anti-morphism by mapping basis to basis.  We can
characterize the extension also by the relation (\ref{six}): for
$v_1,v_2,\ldots, v_k \in T_e(G)$ and $f$ a local section of $\shf_M$,
\[
(V_1^rV_2^r\ldots V_k^r \otimes id_G)(\sigma^*f) = \sigma^*(\rho
(v_1v_2\ldots v_k)f) = \sigma^*(\rho(v_1)\rho(v_2)\ldots
\rho(v_k)f).
\]
\end{proof}


\chapter{The Frobenius Theorem} \label{sgchap4}

\section{The Local Frobenius Theorem}

We want a mechanism by which we can construct a subsupermanifold of a given supermanifold $M$.  In this chapter we present a construction from the tangent bundle of $M$.  We first prove the super extension of the Frobenius theorem on manifolds, then prove a global result.

Let $M$ be a supermanifold with tangent bundle $\text{Vec}_M$.

\begin{definition} A \textit{distribution} on $M$
is an $\shf_M$-submodule $\mathcal{D}$ of $\text{Vec}_M$ which is locally a direct factor.
\end{definition}

\begin{definition}
We say that a distribution $\mathcal{D}$ is \textit{integrable} if it is stable under the bracket on $\text{Vec}_M$, i.e. for $D_1, D_2 \in \mathcal{D}$, $[D_1, D_2] \in \mathcal{D}$.
\end{definition}

\begin{lemma}
\label{distislocfree}
Any distribution $\mathcal{D}$ is locally free.
\end{lemma}

\begin{proof}
By definition, a distribution is a locally direct subsheaf of the
tangent sheaf $\text{Vec}_M$.  Let $x \in |M|$, then
\[
T_x(M) = \mathcal{D}_x \oplus D'
\]
where $\mathcal{D}_x$ is a subsuper vector space of
$T_x(M)$ and we may say that $\mathcal{D}_x$ has basis
$\{s_1, s_2, \ldots, s_k\}$. Then by Nakayama's Lemma (\ref{snakayama}), the $\{s_i\}$
correspond to vector fields which span $\mathcal{D}$ in a
neighborhood  of $x$, and by the locally direct property of a
distribution, these vector fields are linearly independent in this
neighborhood.  Hence $\mathcal{D}$ is actually locally free.
\end{proof}

We can then define the \textit{rank} of a distribution.

\begin{definition}

Let $\mathcal{D}$ be a distribution as above.  Then
$\text{rank}(\mathcal{D})$ is the dimension of $\mathcal{D}_x$ for
$x \in |M|$.  This definition is well-defined thanks to Lemma \ref{distislocfree}.
\end{definition}

Now we prove a series of lemmas before we prove the local Frobenius theorem on supermanifolds.

\begin{remark}
Note that all the following lemmas which pertain to the local Frobenius theorem are local results.  Thus it suffices to consider the case $M = \R^{p|q}$ in a coordinate neighborhood of the origin.
\end{remark}

\begin{lemma}
\label{Dsupercommutes}
Let $\mathcal{D}$ be an integrable distribution.  Then there exist linearly independent
supercommuting vector fields which span $\mathcal{D}$.
\end{lemma}

\begin{proof}
Let $X_1, \ldots, X_r, \chi_1, \ldots, \chi_s$ be a basis for $\mathcal{D}$ and
let $(t, \theta) = (t^1, \ldots, t^p$, $\theta^1, \ldots,
\theta^q)$ be a local set of coordinates.  Then we can express the
vector fields:
\begin{equation}
\label{Dexpressed}
\begin{array}{lcl}
X_j & = & \sum_i a_{ij}\frac{\partial}{\partial t^i} +
\sum_l \alpha_{lj} \frac{\partial}{\partial \theta^l}\\
\chi_k & = & \sum_i \beta_{ik}\frac{\partial}{\partial t^i} +
\sum_l b_{lk}\frac{\partial}{\partial \theta^l}.
\end{array}
\end{equation}

The coefficients form an $r|s \times p|q$ matrix $T$;
\[
T = \left( \begin{array}{cc} a & \alpha\\ \beta & b \end{array}
\right)
\]
of rank $r|s$ since the $\{X_i, \chi_j\}$ are linearly independent.  This is to say that the submatrix $(a)$ has rank $r$

and rank$(b)$ = $s$.  Then by renumerating coordinates $(t, \theta)$, we may assume that
\[
T = \left( T_0 | * \right)
\]
where $T_0$ is an invertible $r|s \times r|s$ matrix.  Multiplying
$T$ by any invertible matrix on the left does not change the row
space of $T$ (i.e. the distribution $\mathcal{D}$), so we can
multply by $T_0^{-1}$ and assume that
\[
T = \left( \begin{array}{ccc} I_r & 0 & *\\ 0 & I_s & * \end{array}
\right),
\]
which is to say that we may assume that
\begin{equation}
\begin{array}{lcl}
X_j  & = & \frac{\partial}{\partial t^j} + \sum_{i=r+1}^p a_{ij} \frac{\partial}{\partial t^i} + \sum_{l = s+1}^q \alpha_{lj} \frac{\partial}{\partial \theta^l}\\
\chi_k  & = & \frac{\partial}{\partial \theta^k} + \sum_{l = s+1}^q b_{lk}\frac{\partial}{\partial \theta^l} + \sum_{i = r+1}^p \beta_{ij} \frac{\partial}{\partial t^i}.
\end{array}
\label{triangle}
\end{equation}

We then claim that $[X_j,X_k] = 0$.  By the involutive property of
$\mathcal{D}$, we know that
\[
[X_j,X_k] = \sum_{i=1}^r f_i X_i + \sum_{l = 1}^s \varphi_l
\chi_l
\]
where the $f_i$ are even functions and the $\varphi_l$ are
odd functions.  Then by (\ref{triangle}), $f_i$ is the coefficient
of the $\frac{\partial}{\partial t^i}$ term in the vector field
$[X_j,X_k]$.  However, again by (\ref{triangle}), it is clear that
$[X_j,X_k]$ has only $\frac{\partial}{\partial t^i}$ terms for
$i>r$, and so we have that $f_i = 0$ for all $i$.  Similarly,
$[X_j,X_k]$ has only $\frac{\partial}{\partial \theta^l}$ terms
for $l > s$, hence also $\varphi_l = 0$ for all $l$.

The cases $[X_j, \chi_k] = 0$ and $[\chi_l, \chi_k] = 0$
follow by using the same argument above.
\end{proof}

\begin{lemma}
\label{1|0-lemma}
Let $X$ be an even vector field.  There exist local coordinates so that
\[
X = \frac{\partial}{\partial t^1}.
\]
\end{lemma}

\begin{proof}
Let $r = 1$, i.e. we begin with a single even vector field $X$, and we want to show that we may express $X=\frac{\partial}{\partial t^1}$ in some coordinate system.  Let $\mathcal{J}$ be the ideal generated by the odd functions on $\R^{p|q}$.  Then since $X$ is even, $X$ maps $\mathcal{J}$ to itself.  Thus $X$ induces a vector field, and hence an integrable distribution, on the reduced space $\R^p$.  Then we may apply the classical Frobenius theorem to get a coordinate system where $X = \frac{\partial}{\partial t^1}$ (mod $\mathcal{J}$).  So we may assume

\[
X = \frac{\partial}{\partial t^1} + \sum_{i\geq 2} a_i \frac{\partial}{\partial t^i} + \sum_j \alpha_j \frac{\partial}{\partial \theta^j}
\]
where the $a_i$ are even, $\alpha_j$ are odd, and $a_i, \alpha_j \in \mathcal{J}$.  That the $a_i$ are even implies that $a_i \in \mathcal{J}^2$.  Moreover, we can find an even matrix $(b_{jk})$ so that $\alpha_j = \sum_k b_{jk}\theta^k$ (mod $\mathcal{J}^2$), and so modulo $\mathcal{J}^2$ we have that
\[
X = \frac{\partial}{\partial t^1} + \sum_{j,k} b_{jk}\theta^k
\frac{\partial}{\partial \theta^j}.
\]
Let $(t, \theta) \mapsto (y, \eta)$ be a change of coordinates where $y = t$ and $\eta = g(t)\theta$ for $g(t) = g_{ij}(t)$ a suitable invertible matrix of smooth functions, that is, $\eta_j = \sum_i g_{ij}(t)\theta^i$.  Then
\begin{equation}
\label{1|0 case}
X = \frac{\partial}{\partial y^1} + \sum_{jk} \theta^k(\frac{\partial g_{jk}}{\partial t^1} + \sum_l g_{jl}b_{lk}) \frac{\partial}{\partial \eta^j},
\end{equation}
and we choose $g(t)$ so that it satisfies the matrix differential equation and initial condition
\[
\frac{\partial g}{\partial t^1} = -gb, \hspace{.3in} g(0) = I.
\]
Then from (\ref{1|0 case}) we may then assume that modulo $\mathcal{J}^2$,
\[
X = \frac{\partial}{\partial y^1}.
\]

Next we claim that  if $X = \frac{\partial}{\partial t^1}$ (mod $\mathcal{J}^k$), then $X = \frac{\partial}{\partial t^1}$ (mod $\mathcal{J}^{k+1}$).  Since $\mathcal{J}$ is nilpotent, this claim will imply the result for the $1|0$-case.

Again, let $(t, \theta) \mapsto (y, \eta)$ be a change of coordinates so that $y^i = t^i +c_i$ and $\eta^j = \theta^j + \gamma_j$ for $c_i, \gamma_j \in \mathcal{J}^k$ suitably chosen.  In the $(t, \theta)$ coordinate system, let
\[
X = \frac{\partial}{\partial t^1} + \sum_{i\geq 2} h_i \frac{\partial}{\partial t^i} + \sum_u \varphi_u \frac{\partial}{\partial \theta^u}
\]
for $h_i, \varphi_u \in \mathcal{J}^k$.  In the new coordinate system, this becomes
\[
X = \frac{\partial}{\partial y^1} + \sum_i (h_i + \frac{\partial c_i}{\partial t^1}) \frac{\partial}{\partial y^i} + \sum_l (\varphi_l + \frac{\partial \gamma_l}{\partial t^1}) \frac{\partial}{\partial \eta^l} + Y
\]
for some $Y = 0$ (mod $\mathcal{J}^{k+1}$) since $2k-1 \geq k+1$ for $k \geq 2$.  So choose the $c_i$ and $\gamma_l$ so that they satisfy the differential equations
\[
\frac{\partial c_i}{\partial t^1} = -h_i, \hspace{.3in} \frac{\partial \gamma_l}{\partial t^1} = -\varphi_l,
\]
and we get that $X = \frac{\partial}{\partial y^1}$ (mod $\mathcal{J}^{k+1}$) as we wanted.
\end{proof}

The above Lemma \ref{1|0-lemma} sets us up to prove the following.

\begin{lemma}
\label{Deven}
Let $\{ X_j \}$ be a set of supercommuting even vector fields.  Then
there exist local coordinates $(t,\theta)$ so that
\[
X_j = \frac{\partial}{\partial t^j} + \sum_{i=1}^{j-1} a_{ij}
\frac{\partial}{\partial t^i}
\]
for some even functions $a_{ij}$.
\end{lemma}

\begin{proof}
Notice that since the $\{ X_j \}$ supercommute, they in fact form a
distribution.  Now we proceed by induction.  The $r=1$ case is presented obove.

We may now assume that we can find coordinates which work for $r-1$
supercommuting vector fields, and we want to prove the lemma for $r$.
Again, assume there are coordinates so that
$X_j = \frac{\partial}{\partial t^j} +\sum_{i=1}^{j-1} a_{ij}
\frac{\partial}{\partial t^i}$ for $j<r$.
Then
\[
X_r = \sum_{i=1}^p f_i \frac{\partial}{\partial t^i} +
\sum_{k=1}^q \varphi_k \frac{\partial}{\partial
\theta^k}
\]
for some even functions $f_i$ and odd functions $\varphi_k$.  The assumption $[X_r, X_j] = 0$ gives
\[
\sum f_i [\frac{\partial}{\partial t^i}, X_j] + \sum \varphi_k
[\frac{\partial}{\partial \theta^k}, X_j] -
\sum (X_jf_i)\frac{\partial}{\partial t^i} -
\sum (X_j\varphi_k)\frac{\partial}{\partial \theta^k} = 0.
\]
We know that $[\frac{\partial}{\partial t^i}, X_j]$ is a linear combination
of $\frac{\partial}{\partial t^l}$ for $l<r$, which means that
$X_jf_i = 0$ for all $j \geq r-1$.  Because the coefficients of the $X_j$ are ``upper triangular" for $j \leq r-1$, we see that $f_i$ depends only on
$(t^r, \ldots, t^p,\theta^1, \ldots, \theta^q)$ for $i \geq r$.
We also have that $[\frac{\partial}{\partial \theta^k}, X_j] = 0$
for all $k$, and so $X_j\varphi_k = 0$ for all $j$ as well.
We can then again conclude that the $\varphi_k$ depend only on
$(t^r, \ldots, t^p,\theta^1, \ldots, \theta^q)$ as well.

Now we can rewrite $X_r$ as follows:
\[
X_r = \left( \sum_{i=1}^{r-1} f_i \frac{\partial}{\partial t^i} \right)
+ \underbrace{\sum_{l=r}^p f_l\frac{\partial}{\partial t^l} +
\sum_{k = 1}^q \varphi_k \frac{\partial}{\partial \theta^k}}
_{\begin{array}{c} ||\\ X_r'\end{array}}.
\]
Here the $X_r'$ depends only on $(t^r, \ldots, \theta^q)$, and so by
an application of the \textit{$1|0$-lemma} on $X_r'$, we may change the
coordinates $(t^r, \ldots, \theta^q)$ so that $X_r' = \frac{\partial}
{\partial t^r}$, and so
\[
X_r = \frac{\partial}{\partial t^r} + \sum_{i=1}^{r-1}f_i'
\frac{\partial}{\partial t^i}
\]
(where the $f_i'$ are the $f_i$ above under the
change of coordinates prescribed by Lemma \ref{1|0-lemma}).

\end{proof}

In fact, the above lemma proves the local Frobenius theorem in the case when
$\mathcal{D}$ is a purely even distribution (i.e. of rank $r|0$).
For the most general case we need one more lemma.

\begin{lemma}
\label{Dodd}
Say $\chi$ is an odd vector field so that $\chi^2 = 0$ and that
$\text{Span}\{\chi\}$ is a distribution.  Then there exist coordinates so
that locally $\chi=\frac{\partial}{\partial \theta^1}$.
\end{lemma}

\begin{proof}
As we have previously remarked, since we want a local result, it suffices to assume that $\chi$ is a
vector field on $\R^{p|q}$ near the origin.
Let us say $(y, \eta)$ are coordinates on $\R^{p|q}$.  Then
\[
\chi = \sum_i \alpha_i(y, \eta)\frac{\partial}{\partial y^i}
+ \sum_j a_j(y, \eta)\frac{\partial}{\partial
\eta^j}
\]
where the $\alpha_i$ are odd, the $a_{\sigma}$ are even, and we
may assume that $a_1(0) \neq 0$.

Now consider the map
\[
\pi: \R^{0|1} \times \R^{p|q-1} \longrightarrow \R^{p|q}
\]
given by
\[
\begin{array}{lcl}
y^i & = & t^i + \epsilon \alpha_i(t,0,\hat{\theta}),\\
\eta^1 & = & \epsilon a_1(t,0,\hat{\theta}),\\
\eta^{j\geq 2} & = & \theta^j + \epsilon a_j(t,0,\hat{\theta})
\end{array}
\]
where $\epsilon$ is the coordinate on $\R^{0|1}$ and $(t^1, \ldots, t^p, \theta^2, \ldots, \theta^q)$ are the coordinates on $\R^{p|q-1}$, and
$\hat{\theta}$ denotes the $\theta$-indices $2, \ldots, q$.  The $\alpha(t, 0 , \hat{\theta})$ and $a(t, 0 , \hat{\theta})$
are the functions $\alpha_i$ and $a_{\sigma}$ where
we substitute $t$ for $y$, let $\theta^1 = 0$, and
substitute $\hat{\theta}$ for $\eta^2, \ldots, \eta^q$.  We claim that the map
$\pi$ is a diffeomorphism in a neighborhood of the origin.  Indeed,
the Jacobian of $\pi$ at $0$ is
\[
J = \text{Ber} \left( \begin{array}{ccc} I_p & * & 0\\
0 & a_1(0) & 0\\ 0 & * & I_{q-1} \end{array} \right) = a_1^{-1}(0) \neq 0.
\]
So we may think of $(t,\epsilon, \hat{\theta})$
as coordinates on $\R^{p|q}$ with $\pi$ being a change of coordinates.
Under this change of coordinates, we have
\[
\frac{\partial}{\partial \epsilon} = \sum_i
\frac{\partial y^i}{\partial \epsilon}\frac{\partial}{\partial y^i} +
\sum_j \frac{\partial \eta^j}{\partial \epsilon}\frac{\partial}
{\partial \eta^j},
\]
which amounts to
\[
\frac{\partial}{\partial \epsilon} = \sum_i \alpha_i
(t, 0, \hat{\theta})
\frac{\partial}{\partial y^i} + \sum_j a_j(t, 0, \hat{\theta}) \frac{\partial}{\partial
\eta^j}.
\]
The $\alpha_i(t, 0, \hat{\theta})$ and $a_j(t, 0, \hat{\theta})$ terms must be expressed as functions of $(y, \eta)$.  Notice that by a simple Taylor series expansion, $\alpha_i(y, \eta) = \alpha_i(t^i + \epsilon
\alpha_i, \epsilon a_1, \theta^{k \geq 2}+ \epsilon a_k)
= \alpha_i(t^i,0,\hat{\theta}) + \epsilon \beta_i$ for some
odd function $\beta_i$.  Similarly we get
$a_j(y, \eta) = a_j(t, 0, \hat{\theta}) + \epsilon b_j$
for some even function $b_j$.  Thus we can write
\[
\frac{\partial}{\partial \epsilon} = \chi + \epsilon Z
\]
for some even vector field $Z$.  Recall that $\eta^1 = \epsilon \hat{a}_1$ where $\hat{a}_1$ is an even
invertible section.  Hence $\epsilon = \eta^1 A$ from some invertible
even section $A$.

Then we see that under the change of coordinates given by $\pi$,
\[
\frac{\partial}{\partial \epsilon} - \eta^1
\underbrace{A \cdot Z^*}_{=Z'} = \chi
\]
where $Z^*$ denotes the pullback of $Z$ by $\pi$ and $Z'$ is some even
vector field (since both $A$ and $Z$ are even).  Now,
\[
\begin{array}{lcl}
\chi^2 = 0 & \implies & (\frac{\partial}{\partial \epsilon}
-\eta^1Z')^2 = 0\\
& \implies & \underbrace{(\frac{\partial}{\partial \epsilon})^2}_{=0}

- \frac{\partial}{\partial \epsilon} (\eta^1Z') - (\eta^1Z')
\frac{\partial}{\partial \epsilon} +
\underbrace{(\eta^1Z')^2}_{=0} = 0\\
& \implies & -\hat{a}_1Z' + \eta^1\frac{\partial}{\partial \epsilon}Z'
-\eta^1Z'\frac{\partial}{\partial \epsilon} = 0\\
& \implies & \hat{a}_1Z' = 0\\
& \implies & Z' = 0,
\end{array}
\]
so we really have $\frac{\partial}{\partial \epsilon} = \chi$ under the change of coordinates.
\end{proof}

Now we can prove of the full local Frobenius theorem.

\begin{theorem} \text{(Local Frobenius Theorem)}
Let $\mathcal{D}$ be an integrable (involutive) distribution of
rank $r|s$.  Then there exist local coordinates so that
$\mathcal{D}$ is spanned by
\[
\frac{\partial}{\partial t^1},
\ldots, \frac{\partial} {\partial t^r}, \frac{\partial}{\partial
\theta^1}, \ldots, \frac{\partial}{\partial \theta^s}.
\]
\end{theorem}

\begin{proof}

Let $\{ X_1, \ldots, X_r,\chi_1, \ldots, \chi_s \}$ be a basis of vector
fields for the distribution $\mathcal{D}$.  By Lemma \ref{Dsupercommutes} we may
assume that these basis elements supercommute, so then
$\mathcal{D}' = \text{span}\{X_1, \ldots, X_r\}$ is a
subdistribution, and by lemma \ref{Deven} we get that there exist
coordinates so that $X_i = \frac{\partial} {\partial t^i}$.

We then use the fact that $[\chi_1, X_i] = 0$ for all $i$, to see that
$\chi_1$ depends only on coordinates $(t^{r+1}, \ldots, \theta^q)$ (as in the proof of 
Lemma \ref{Deven}).  In fact, this is not completely accurate.  If we express $\chi_1$ as in (\ref{Dexpressed}), we see that it is only the $\beta_{ik}$ and $b_{lk}$ which depend only on the coordinates $(t^{r+1}, \ldots, \theta^q)$.  However, we can always kill off the first $r$ $\partial/\partial t^i$ terms by subtracting off appropriate linear combinations of the $\{X_1 = \partial/\partial t^1, \ldots, X_r = \partial/\partial t^r \}$.

Then since $\chi_1^2 = 0$, by Lemma \ref{Dodd} we may change only the coordinates $(t^{r+1}, \ldots, \theta^q)$ and express $\chi_1 =
\frac{\partial}{\partial \theta^1}$.  For $\chi_2$ we apply the same
idea: that $[\chi_2, X_i] = 0$ and $[\chi_2,\chi_1] = 0$ again shows that
$\chi_2$ depends only on coordinates $(t^{r+1}, \ldots, t^p, \theta^2,
\ldots, \theta^q)$, and again applying Lemma \ref{Dodd} gives $\chi_2 =
\frac{\partial}{\partial \theta^2}$.  And so on with $\chi_3, \ldots,
\chi_s$.
\end{proof}

We are now in a position to state and prove the global Frobenius theorem on supermanifolds.

\section{The Global Frobenius Theorem on Supermanifolds}

\begin{theorem} \text{(Global Frobenius Theorem)}
Let $M$ be a $C^{\infty}$-supermanifold, and let $\mathcal{D}$ be
an integrable distribution on $M$.  Then given any point of $M$
there is a unique maximal supermanifold corresponding to
$\mathcal{D}$ which contains that point.
\end{theorem}

\begin{proof}
Let $\mathcal{D} = \text{span} \{ X_1, \ldots, X_r, \chi_1, \ldots, \chi_s \}$ as in the previous section (again the $X_i$ are even and the $\chi_j$ are odd).  Then let $\mathcal{D}_0 = \text{span} \{X_1, \ldots X_r \}$; this subdistribution maps
odd sections to odd sections, and so descends to an integral
distribution $\widetilde{\mathcal{D}}_0$ on $\widetilde{M}$.  Let $x \in |M|$. Then by the
classical global Frobenius Theorem, there is a unique maximal
integral manifold $\widetilde{M}_x \subset \widetilde{M}$ of
$\widetilde{\mathcal{D}}_0$ containing $x$.  We want to build a sheaf of commutative superalgebras on $\widetilde{M}_x$.

By the local Frobenius
theorem, given any point $y \in |M|$, there exists an open
coordinate neighborhood around $y$, $U_y \subset |M|$, so that $U_y$ is characterized by coordinates $(t, z, \theta,
\eta)$ (i.e. $\shf_M(U_y)= C^{\infty}(t, z)[\theta, \eta]$) where $\mathcal{D}$ is given by the

$\mathcal{T}_M$-span of $\{ \frac{\partial}{\partial t},
\frac{\partial}{\partial \theta} \}$.  Now let $U\subset |M|$
and define the following presheaf $\mathcal{I}$:
\[
\mathcal{I}(U) = \langle \{ f \in \shf_M(U) \hspace{.1in} |
\hspace{.1in} \text{$\forall y \in \widetilde{M}_x \cap U$,
$\exists$ $V_y \subset U$ so that $f|_{V_y}\in C^{\infty}(z)[\eta]$} \} \rangle.
\]
We claim that $\mathcal{I}$ is a subsheaf of $\shf_M$.  Again let
$U\subset |M|$ be an open subset and let $\{U_{\alpha} \}$ be an
open covering of $U$ so that for $s_{\alpha} \in U_{\alpha}$ we
have $s_{\alpha}|_{U_{\alpha} \cap U_{\beta}} =
s_{\beta}|_{U_{\alpha} \cap U_{\beta}}$.  We know that there
exists a unique $s\in \shf_M(U)$ so that $s|_{U_\alpha} =
s_{\alpha}$.  Let $y\in M_x \cap U$, then $y\in U_{\alpha}$ for
some $\alpha$.  Then there exists $V_y \subset U_{\alpha}$ where $s|_{V_y} \in C^{\infty}(z)[ \eta]$ since $s|_{V_y} = s_{\alpha}|_{V_y}$.
Hence $\mathcal{I}$ is a subsheaf of $\shf_M$.  Moreover
$\mathcal{I}$ is an ideal sheaf by construction.

It is clear that if $p \notin \widetilde{M}_x$, then
$\mathcal{I}_p = \shf_{M,p}$ since we can find some neighborhood
of $p$, $W_p\cap \widetilde{M}_x = \emptyset$, where
$\mathcal{I}(W_p) = \shf_M(W_p)$.  Thus the support of
$\mathcal{I}$ is $\widetilde{M}_x$, and we have defined a
quasi-coherent sheaf of ideals with support $\widetilde{M}_x$
which defines a unique closed subspace of $M$.  By going to
coordinate neighborhoods it is clear that this closed subspace is
in fact a closed subsupermanifold which we shall now call
$M_x$.

The maximality condition is clear.  From the classical theory we have that the reduced space is maximal, and locally we can verify that we have the maximal number of odd coordinates that $\mathcal{D}$ allows.
\end{proof}

\section{Lie Subalgebra, Subgroup Correspondence}

From a slice of the tangent bundle of a given supermanifold, the global Frobenius theorem allows us to build a subsupermanifold.  We now use this construction to make the super Lie group/algebra super Lie subgroup/subalgebra correspondence.  We begin with a technical lemma we will need later.

\begin{lemma}
\label{tech-sub}
Let $M$ be a supermanifold, $N\subset M$ a subsupermanifold, and
let $\varphi:M \longrightarrow M$ be a diffeomorphism so that
$\varphi(N)\subset N$ and
$\widetilde{\varphi}(\widetilde{N})=\widetilde{N}$. Then
$\varphi(N)=N$.
\end{lemma}

\begin{proof}
Because $\varphi$ is a diffeomorphism, $\varphi(N)$ is a super
submanifold of $M$ with the same super dimension as $N$.  We
assume that they both have the same underlying space.  Then since
they have the same odd dimension and $\varphi(N)$ sits inside $N$,
they must be the same space.  This can be checked
at the coordinate neighborhood level.
\end{proof}

Let $G$ be a super Lie group with super Lie algebra
$\lie(G)$.  Let $(H',\mathfrak{h})$ be a pair consisting of an ordinary
Lie subgroup and a Lie superalgebra, so that
\[
\begin{array}{ll}
1. & H' \subset \widetilde{G} \text{ is a Lie subgroup};\\
2. & \mathfrak{h}\subset \lie(G) \text{ is a super Lie subalgebra}.
\end{array}
\]

\begin{theorem}

There is a super Lie subgroup $H$ of $G$ so that
\[
\begin{array}{ll}
1. & \widetilde{H} = H';\\
2. & \lie(H) = \mathfrak{h}.
\end{array}
\]
\end{theorem}

\begin{proof}
Let $\mathcal{D}$ be the distribution generated by $\mathfrak{h}$
on $G$, i.e. $\mathcal{D} = \langle \mathfrak{h}
\rangle_{\mathcal{T}_G}$.

Then we can use the Global Frobenius Theorem to get maximal
integral super submanifolds $G_p$ through each point $p \in G$
which correspond to $\mathcal{D}$.  Since $H'$ is second
countable, it is the union of a countable number of connected
components.  Take a collection of points $p_i \in H'$, each of
which corresponds to exactly one of the $G_{p_i}$ and let $H$ be
the supermanifold of the union of these maximal integral super
submanifolds, i.e. $H= \cup G_{p_i}$.  This construction makes it
clear that $\widetilde{H}= H'$.

All that is left to show is that $H$ is in fact a super Lie group.
We already have a morphism $\mu_H:H \times H \longrightarrow G$
which comes from restricting the multiplication morphism $\mu:G \times G
\longrightarrow G$. This gives the sheaf map
\[ \mu_H^*:\shf_G \longrightarrow \shf_{H \times H}.\]
Let us restrict our view to some coordinate neighborhood of $G$,
and let $(t,z,\theta,\eta)$ represent the coordinates on $G$ so that $H$ is
described locally by the vanishing of the coordinates $(z, \eta)$.  This is equivalent to saying that

$\mathfrak{h}$ kills these coordinates; for $h\in \mathfrak{h}$,
$h(z)=h(\eta)=0$. To show that $\mu_H^*$ is actually a map from
$\shf_H \longrightarrow \shf_{H \times H}$ we have to show that
$\mu_H^*$ vanishes on $(z,\eta)$.

Let $h\in\mathfrak{h}$.  Because $\mathfrak{h} \subset \lie(G)$,
$h$ is left invariant, and hence commutes with $\mu$ (more precisely, with
$\mu_H$).  Then the following diagram commutes:
\[
\begin{array}{rcl}
\shf_G & \stackrel{\mu_H^*}{\longrightarrow} & \shf_{H\times H}\\
_h \downarrow & & \downarrow_{id_H \tensor h}\\
\shf_G & \stackrel{\mu_H^*}{\longrightarrow} & \shf_{H\times H}.
\end{array}
\]
We already know $\mu_H^*\circ h (z) = 0$, from which
commutativity gives that $id_H \tensor h \circ \mu_H^*
(z) = 0$. Thus $\mu_H^*(z)$ is killed by $id_H \tensor
h$ which implies that $\mu_H^*(z) = 0$ since $h$ ranges over all left invariant vector fields.  Similarly,
$\mu_H^*(\eta) = 0$. Thus we have a product structure on $H$.

Last we show that there is an inverse map $\iota:H \longrightarrow
H$. Consider the morphism
\[ \nu: G \times G \longrightarrow G \times G \]
given by $\nu = (id, \mu)$ which we define as follows via $T$-points.  Let $T$ be any supermanifold.  Then
\[ \nu(T): G(T) \times G(T) \longrightarrow G(T) \times G(T) \]
so that for $g,h \in G(T)$,
\[ (g,h) \mapsto (g,gh). \]
If is then clear that $\nu(T)$ is bijective for all $T$, and so $\nu$ is thus a diffeomorphism of $G \times G$ to itself by

Yoneda's Lemma.  Recall that $H\times H \subset G \times G$ is a closed
subsupermanifold.  From the arguments above, $\nu$ maps $H\times
H$ into itself.  Moreover, we claim that
$\widetilde{\nu}(\widetilde{H \times H})=\widetilde{H \times H}$.
First note that $\widetilde{H \times H} = \widetilde{H} \times
\widetilde{H}$ and that for any ordinary manifold $S$, we have
that $\widetilde{\nu}(S):\widetilde{H}(S) \times \widetilde{H}(S)
\longrightarrow \widetilde{H}(S) \times \widetilde{H}(S)$ is a
bijection because $\widetilde{H}$ is a Lie group.  Thus
\[
\widetilde{\nu}(\widetilde{H}\times \widetilde{H}) = \widetilde{H}
\times \widetilde{H}
\]
and we can use the general Lemma (\ref{tech-sub}) to see that $\nu(H \times
H)=H \times H$ from which it follows that the inverse map of $G$
descends to $H$.

The necessary diagrams commute (associativity, inverses, etc.)
because they do for $G$ and all the maps for $H$ are derived from
those of $G$.  We have thus produced a super Lie subgroup $H$ of
$G$ with the additional property that $\lie(H)=\mathfrak{h}$.
\end{proof}


\chapter{Supervarieties and Superschemes} \label{sgchap5}

\section{ Basic definitions}

\medskip

In this section we give the basic definitions of
algebraic supergeometry. Because we are in need of a more
general setting in the next two chapters
we no longer assume the ground field to
be $\R$.

\medskip

Let $k$ be a commutative ring.

\medskip

Assume all superalgebras are associative, commutative (i.e.
$xy=(-1)^{p(x)p(y)} yx$) with unit and over $k$.  
We denote their category with $\salg$.
For a superalgebra $A$ let
$J_A$ denote the ideal generated by the odd elements i. e. $J_A=<A_1>_A$.
Denote the quotient $A/J_A$ by $A^r$.

\medskip

In chapter \ref{sgchap2} we have introduced
the notion of {\textit {superspace}} and of
{\textit {superscheme}}. Recall that a superspace
$X=(|X|,\cO_X)$ is a topological space
$|X|$ together with a sheaf of superalgebras $\cO_X$ such that
$\cO_{X,x}$ is a local superalgebra, i.e. it has a unique two sided
maximal homogeneous ideal.

\medskip

The sheaf of superalgebras $\cO_X$ is a
sheaf of $\cO_{X,0}$-modules,
where $\cO_{X,0}(U)$ $=_{\defi}$ $\cO_{X}(U)_0$, $\forall U$ open in $|X|$.

\medskip

Let $\cO_X^r$ denote the sheaf of algebras:
$$
\cO_X^r(U)=\cO_X(U)/J_{\cO_X(U)}
$$
We will call $X^r=(|X|, \cO_X^r)$ the \textit{ reduced space} associated
to the superspace $X=(|X|,\cO_X)$. 
This is a locally ringed space in the classical sense.

\medskip

Recall that given two superspaces $X=(|X|,\cO_X)$ and $Y=(|Y|,\cO_Y)$
a \textit{ morphism} $f : X \lra Y$ of superspaces
is given by a pair $f =(|f|,f^*)$ such that

\noindent
1. $|f|:X \rightarrow Y$ is a continuous map.

\noindent
2. $f^*:\cO_Y \rightarrow f_*\cO_X$
is a map of sheaves of superalgebras on $|Y|$, 
that is for all
$U$ open in $|Y|$ there exists
a morphism $f^*_U:\cO_Y(U) \rightarrow \cO_X(|f|^{-1}(U))$.

\noindent
3. The map of local superalgebras
$f^*_p:\cO_{Y,|f|(p)} \rightarrow \cO_{X,p}$ is a local morphisms i.e.
sends the maximal ideal of $\cO_{Y,|f|(p)}$ into the maximal
ideal of  $\cO_{X,p}$.

\medskip


\medskip

Recall that a \textit{ superscheme} $S$ is a superspace $(|S|, \cO_S)$
such that $(|S|, \cO_{S,0})$ is a quasi coherent sheaf of
$\cO_{x,1}$-modules. 
A \textit{ morphism} of superschemes 
is a morphisms of the corresponding superspaces.

\medskip

For any open $U \subset X$ we have the superscheme
$U=(|U|,\cO_X|_U)$,  called 
\textit{ open subscheme} in the superscheme $X$.

\medskip

One of the most important examples of superscheme is given by the spectrum
of the even part of a given superalgebra (the topological structure)
together with
a certain sheaf of superalgebras on it that plays the role of
the structural sheaf in the classical theory.
Let's see this construction in detail.

\medskip

\begin{definition} $\uspec A$.

Let $A$ be an object of $\salg$. We have that
$\spec(A_0)=\spec(A^r)$, since the
algebras $A^r$ and $A_0$ differ only by nilpotent elements. 

Let's consider $\cO_{A_0}$ the structural sheaf of $\spec(A_0)$.
The stalk of the sheaf at the prime $\bp\in
\spec(A_0)$ is the localization of $A_0$ at $\bp$.
As for any superalgebra, $A$ is a
module over $A_0$. We have indeed a sheaf $\cO_A$ of
$\cO_{A_0}$-modules over $\spec A_0$ with stalk $A_{\bp}$, the localization
of the $A_0$-module $A$ over each prime
$\bp \in \spec(A_0)$.
$$
A_{\bp}=\{ {f \over g} \quad | \quad f\in A, g \in A_0-\bp\}
$$
The localization $A_\bp$ has a unique two-sided maximal ideal 
which consists of
the maximal ideal in the local ring
$(A_\bp)_0$ and the generators of $(A_\bp)_1$ as $A_0$-module.
For more details on this construction see \cite{Harts} II \S 5.

$\cO_A$  is a sheaf of superalgebras  and
$(\spec A_0,\cO_A)$ is a superscheme that we will
denote with $\uspec A$. 
Notice that on the open sets:
$$
U_f=\{ \bp \in \spec A_0 | (f) \not\subset \bp \}, \qquad f \in A_0
$$
we have that $\cO_A(U_f) = A_f=\{a/f^n \quad | \quad a \in A\}$.
\end{definition}

\medskip

\begin{definition}
An \textit{ affine superscheme} is a superspace
that is isomorphic to $\uspec A$ for some
superalgebra $A$ in $\salg$.
An \textit{ affine  algebraic supervariety} is a superspace
isomorphic to $\uspec A$ for some \textit{ affine}
superalgebra $A$ i. e. a finitely generated superalgebra such
that $A/J_A$ has no nilpotents.
We will call
$A$ the \textit{ coordinate ring} of the supervariety.

\end{definition}

\medskip

\begin{proposition}
A superspace $S$ is a superscheme if and only if it is locally
isomorphic to $\uspec A$ for some superalgebra $A$, i. e.
for all $x \in |S|$, there exists $U_x \subset |S|$ open such
that $(U_x, \cO_S|_{U_x}) \cong \uspec A$. (Clearly $A$
depends on $U_x$).
\end{proposition}
\medskip

\begin{proof}
Since $S$ is a superscheme, by definition
$S'=(|S|, \cO_{S,0})$ is an ordinary scheme, that is, it admits an open
cover $S' = \cup V_i$ so that $V_i \cong \uspec A_{i,0}$
where $A_{i,0}$ is a commutative algebra. 
Let  $x \in |S|$ and let
$U_i=(|V_i|, \cO_S|_{V_i})$, such that $x \in |U_i|$.

The $U_i$ can be chosen so that
there exists a $A_{i,0}$-module $A_{i,1}$ such that
$\cO_{S,1}|_{U_i} \cong \cO_{A_{i,1}}$, where
$\cO_{A_{i,1}}$ denotes the sheaf induced by the $A_{i,0}$-module
$A_{i,1}$. 
So we have that:
$$
\cO_S|_{U_i}=\cO_{S,0}|_{U_i} \oplus \cO_{S,1}|_{U_i}=
\cO_{A_{i,0}} \oplus \cO_{ A_{i,1}}=\cO_{A_i}
$$

Since these are sheaves of superalgebras, $A_i$ is
also a superalgebra, in fact 
$A_i=\cO_S|_{U_i}(U_i)$. Hence $U_i \cong \uspec A_i$. 
The other direction is clear. 
\end{proof}

\medskip

Given a superscheme $X=(|X|,\cO_X)$,
the scheme $(|X|,\cO_X^r)$ is called
the \textit{ reduced scheme} associated
to $X$.
Notice that the reduced scheme associated to a given superscheme
may not be reduced, i.e. $\cO_X^r(U)$, $U$ open in $X$,
can contain nilpotents.

\medskip

As in the classical setting we can define closed
subsuperschemes.

\medskip

\begin{definition}
A \textit{ closed subsuperscheme} $Y$ of a given superscheme $X$
is such that $|Y| \subset |X|$ and $\cO_Y = \cO_X/\cI$ for
a quasi-coherent sheaf of ideals in $\cO_X$.

Notice that if $X= \uspec A$, closed subschemes are in one
to one correspondence with ideals in $A$ as it happens in
the ordinary case.
\end{definition}

\medskip

\begin{example}

1. \textit{ Affine superspace $\bA^{m|n}$}. \label{affinespace}

Consider the polynomial superalgebra $k[x_1 \dots x_m, \xi_1 \dots \xi_n]$
over an algebraically closed field $k$
where $x_1 \dots x_m$ are even indeterminates and $\xi_1 \dots \xi_n$
are odd indeterminates (see chapter \ref{sgchap1}).
We will call $\uspec k[x_1 \dots x_m, \xi_1 \dots \xi_n]$ the affine
superspace of superdimension $m|n$ and we denote it by
$\bA^{m|n}$.
$$
k[\bA^{m|n}]=k[x_1 \dots x_m, \xi_1 \dots \xi_n].
$$

As a topological space
$\spec k[x_1 \dots x_m, \xi_1 \dots \xi_n]_0$ will consists of
the even maximal ideals
$$
(x_i-a_i, \xi_j\xi_k),
\qquad i=1 \dots m, \quad j,k=1 \dots n
$$
and the even prime ideals
$$
(p_1 \dots p_r, \xi_j\xi_k),
\qquad i=1 \dots m, \quad j,k=1 \dots n
$$
where $(p_1 \dots p_r)$ is a prime ideal in $k[x_1 \dots x_m]$.

The structural sheaf of $\bA^{m|n}$ will have stalk at the
point $\bp \in \spec k[\bA^{m|n}]_0$:
$$
k[\bA^{m|n}]_{\bp}=\{{f \over g} \quad | \quad  f \in k[\bA^{m|n}],\quad
g\in k[\bA^{m|n}]_0, \quad  g\notin \bp\}.
$$
\end{example}

\medskip

2. \textit{ Supervariety over the sphere $S^2$}.

Consider the polynomial
superalgebra generated over an algebraically closed field $k$
$k[x_1,x_2, x_3,\xi_1, \xi_2, \xi_3],$ and the ideal
$$
\cI=(x_1^2+ x_2^2+ x_3^2-1,\; x_1\cdot \xi_1+ x_2\cdot \xi_2+ x_3 \cdot \xi_3).
$$
Let $k[X]=k[x_1, x_2, x_3,\xi_1, \xi_2, \xi_3]/\cI$
and $X=\uspec k[X]$. $X$ is a 
supervariety whose reduced variety $X^r$ is the
sphere $S^2$. A maximal ideal in $k[X]_0$ is given by:
$$
\bm=(x_i-a_i,\, \xi_i\xi_j)\quad \hbox{with}
\;\; i,j=1,2,3, \;\; a_i\in k\;\; \hbox{and}\;\;
a_1^2+a_2^2+a_3^2=1.
$$

\medskip

\begin{observation} \label{equiv}

Let $\sschemesaff$ denote the category of affine superschemes.

The functor
$$
\begin{array}{ccc}
F:\salg^o & \lra & \sschemesaff  \\
A & \mapsto & \uspec A
\end{array}
$$
gives an equivalence between the category of
superalgebras and the category of affine superschemes.
The inverse functor is given by:
$$
\begin{array}{ccc}
G:\sschemesaff & \lra & \salg  \\
\uspec A & \mapsto & \cO_A(A_0)\cong A.
\end{array}
$$
\end{observation}

We need to specify both $F(\phi)$, for $\phi:A \lra B$ and 
$G(f)$, for $f:\uspec A \lra \uspec B$ and show that they
realize a bijection:
$$
\Hom_\salg (A,B)\cong \Hom_\sschemesaff(\uspec B, \uspec A)
$$
Let $\phi:A \lra B$ and $\phi_0=\phi|_{A_0}$. 
We want to build a morphism $f=F(\phi): \uspec B \lra \uspec A$.
We have immediately $|f|: \spec B_0 \lra  \spec A_0$
defined as $|f|(\bp)=\phi_0^{-1}(\bp)$.
We also have a map
$$
f^*_U: \cO_{A}(U) \lra \cO_{B}(|f|^{-1}(U)), \qquad U
\hbox{ open in } \spec A_0,
$$
defined as in the classical
case. That is if $a \in \cO_{A}(U)$, $f_U^*(a): \bp
\mapsto \phi_{\bp}(a(|f|(\bp)))$ and $\phi_\bp: A_{\phi_0^{-1}(\bp)} \lra
B_{\bp}$, $\bp \in \spec B_0$.

Vice-versa if we have a map $f: \uspec B \lra  \uspec A$, since
global sections of the structural sheaves coincide with the
rings $B$ and $A$ respectively, we obtain
immediately a map from $A$ to $B$:
$$
G(f)=f^*_{\spec A_0}: \cO_{A}(\spec A_0) \cong A \lra
\cO_{B}(\spec B_0) \cong B.
$$

\medskip

When we restrict the functor $F$ to the category of
affine superalgebras, it gives an equivalence
of categories between affine superalgebras and affine 
supervarieties. 

\medskip

We now would like to give an example of a non affine superscheme
which is of particular importance: the projective superspace.

\medskip

\begin{example}
\textit{ Projective superspace}

Let $S= k[x_0, \dots x_m, \xi_1 \dots \xi_m]$.
$S=\s-0 \oplus S_1$ is a $\Z/2\Z$ and 
$\Z$ graded algebra and the two gradings are
compatible. Define the topological space $\proj S_0$ as the
set of $\Z$-homogeneous non irrelevant primes in $S_0$,
with the Zariski topology. 
$\proj S_0$ is covered by open affine $U_{i}$ consisting
of those primes non containing $(x_i)$. 
As in the classical setting we have that
$$
U_i = \spec k[x_0 \dots \hat x_i \dots x_m, \xi_1 \dots \xi_n]_0,
\quad i=1 \dots m.
$$
So we can define the sheaves 
$$
\cO_{U_i}=\cO_{k[x_0 \dots \hat x_i \dots x_m, \xi_1 \dots \xi_n]}
$$ 
corresponding to
these open affine subsets. One can check that these sheaves glue
to give a sheaf $\cO_S$ on all $\proj S$. 
So we define \textit{ projective superspace} $\bP^{m|n}$, 
as the superscheme $(\proj S_0, \cO_S)$.

\medskip

The same construction can be easily repeated for a generic 
$\Z$-graded superalgebra.

\end{example}

\bigskip
\section{The functor of points}

\medskip

As in $C^{\infty}$ geometry, we employ the functor of points
approach from algebraic geometry to better handle to 
nilpotent elements and to bring back geometric intuition.

\medskip

\begin{definition}
For a superscheme $X$,  the
\textit{ functor of points} of $X$ is a representable functor
$$
h_X:\sschemes^o \rightarrow \sets,
\qquad
h_X(Y)= \Hom_{\sschemes}(Y,X)
$$
$h_X(Y)$ are called the \textit{ $Y$-points} of the superscheme $X$.
\end{definition}

In the previous chapters, we have used the same notation to denote both
a supergeometric object, say a supermanifold, and its functor of points.
In this chapter, however, we want to make a distinction, since we
will also deal with non representable functors.

\medskip

As in the ordinary setting, the functor of points of
a superscheme $h_X$ is determined by looking at its restriction to the
affine superschemes $h_X^a$, that is looking at the functor
$$
h_X^a:\salg \lra \sets, \qquad h_X^a(A)=\Hom_{\sschemes}(\uspec A,X).
$$

This is proven in the same way as the ordinary
 case. In fact a morphism $\phi \in \Hom(Y,X)$ is
determined by its restrictions to the open affine subschemes
that form an open cover of $Y$.

\medskip

When the superscheme $X$ is affine,
i.e. $X=\uspec R$, $h_X^a$ is representable.
In fact by Observation \ref{equiv}:
$$
h_X^a(A)=\Hom_{\sschemes}(\uspec A,\uspec R)=
\Hom_{\salg}(R,A).
$$

\begin{observation}
Since we have the equivalence of categories between
affine superschemes and superalgebras, we can
define an affine superscheme
equivalently as a representable functor
$$
F:\salg \lra \sets, \qquad F(B)=\Hom_{\salg}(A,B).
$$
\end{observation}

\begin{remark}
To simplify notation we drop the suffix $a$
in $h_X^a$, the context will make clear whether we are
considering $h_X$ or its restriction to affine superschemes.
Moreover, whenever want the restriction of $h_X$
to affine superschemes we will not use a different functor name for
$h_X(A)$, $h_X: \salg \lra \sets$ and for $h_X(\uspec A)$
$h_X: \sschemesaff \lra \sets$.
\end{remark}

\medskip

\begin{observation}
Let $X^0$ be an affine variety over an algebraically closed
field $k$. Consider
an affine supervariety $X$ whose reduced part coincides with
$X^0$. Then one can immediately check that
the $k$-points of $X$ correspond to
the points of the affine variety $X^0$.
\end{observation}

\medskip


\begin{examples}

1. \textit{ Affine superspace revisited}.

Let $A \in \salg$ and let $V=V_0 \oplus V_1$ be a
free supermodule (over $k$). Let $\smod$ denote
the category of $k$-modules.
Define
$$
V(A)=(A \otimes V)_0= A_0 \otimes V_0 \oplus A_1 \otimes V_1.
$$
In general this functor is not representable. However, if $V$ is finite
dimensional we have:
$$
(A \otimes V)_0 \cong 
\Hom_{\smod}(V^*,A) \cong \Hom_{\salg}(\Sym(V^*),A)
$$
where $\Sym(V^*)$ denotes the symmetric algebra over the dual space $V^*$.
Recall that
$V^*$ is the set of linear maps $V \lra k$ not necessarily preserving
the parity, and $\Sym(V^*)=\Sym(V^*_0) \otimes \wedge V^*_1$,
where $\wedge V^*_1$ denotes the exterior algebra over the
ordinary space $V_1$. 

Let's fix a basis for $V$ and let $\dim V=p|q$.
The functor $V$ is represented by:
$$
k[V]=k[x_1 \dots x_p, \xi_1 \dots \xi_q]
$$
where $x_i$ and $\xi_j$ are respectively even and odd indeterminates.

Hence the functor $V$ is the functor of points of the
affine supervariety $\bA^{m|n}$ introduced in Example \ref{affinespace}.

\medskip

We also want to remark that the functor $D_V$ defined as:
$$
D_V(A)=_{\defi}\Hom_{\smod}(V,A)
$$
is representable for any $V$ (not necessarily finite dimensional),
and it is represented by the superalgebra $\Sym(V)$.
Clearly $V=D_V$ when $V$ is finite dimensional.
\medskip

2. \textit{ Supermatrices}.

Let $A \in \salg$.
Define $\rM_{m|n}(A)$ as the set of endomorphisms of the $A$-supermodule
$A^{m|n}$. Choosing coordinates we can write:
$$
\rM_{m|n}(A)=
\left\{ \begin{pmatrix} a & \alpha \\ \beta & b \end{pmatrix} 
\right\}
$$
where $a$ and $b$ are $m \times m$, $n \times n$ blocks of even elements
and $\forall$, $\al$, $\be$ $m \times n$, $n \times m$ blocks of odd elements.

\medskip

This is the functor of points
of an affine supervariety represented by the commutative
superalgebra: $k[\rM(m|n)]=$ $k[x_{ij},\xi_{kl}]$
where $x_{ij}$'s
and $\xi_{kl}$'s are respectively even and odd variables
with
$1 \leq i,j \leq m$ or  $m+1 \leq i,j \leq m+n$,
$1 \leq k \leq m$, $m+1 \leq l \leq m+n$ or
$m+1 \leq k \leq m+n$, $1 \leq l \leq m$.

\medskip

Notice that $\rM_{m|n} \cong h_{\bA^{m^2+n^2|2mn}}$.
\end{examples}

\bigskip

\section{A representability criterion}

\medskip
We now want to single out among all the functors
$F: \salg \lra \sets$
those that are the functor of points of superschemes.

\medskip

We first need
the definition of local functor and open subfunctor.

\medskip
Let $A$ be a superalgebra. Given $f \in A_0$, let $A_f$ denote:
$$
A_f=_{\defi} \{a/f^n | a \in A\}.
$$

The sets $U_f=\spec (A_{f})_0$ are
open sets in the topological space
$X=\spec A_0$. In fact recall that since $A_0$ is an ordinary
commutative algebra, 
by definition the open sets
in the Zariski topology of $\spec A_0$ are:
$$
U_I=\{\bp \in \spec A_0| I \not\subset \bp \}
$$
for all the ideals $I$ in $A_0$.

\medskip

We now want to define the notion of open subfunctor of a functor
$F$. If we assume $F$ is the functor of points
of superscheme $X$ an open subfunctor could simply be defined as the
functor of points of an open superscheme $U \subset X$.
However because we are
precisely interested in a characterization of those $F$ that come from
superschemes, we have to carefully extend this notion.

\medskip

\begin{definition} Let $U$ be a subfunctor of a functor
$F:\salg \lra \sets$ (this means that we have
a natural transformation $U \lra F$ such that 
$U(A) \lra F(A)$ is injective for all $A$).
We say that \textit{ $U$ is an open subfunctor of $F$} 
if for all $A \in \salg$ given
a natural transformations $f:h_{\uspec A} \lra F$, the
subfunctor $f^{-1}(U)$ is equal to $h_V$, for some open $V$ 
in $\uspec A$ where
$$
f^{-1}(U)(R)=_{\defi}f_R^{-1}(U(R)), \qquad
f_R:h_\uspec A(R) \lra  F(R).
$$

We say $U$ is an \textit{ open affine subfunctor} of $F$
if it is open and representable.

\end{definition}

\medskip

\begin{observation} \label{opencover}
Let $X=(|X|,\cO_X)$ be a superscheme
and $U \subset X$ open affine in $X$.
Then $h_U$ is an open affine subfunctor of $h_X$.

By Yoneda's lemma $f:h_{\uspec A} \lra h_X$ corresponds
to a map $f':\uspec A \lra X$. Let $V=f'^{-1}(U)$ open in $\uspec A$.
We claim
$$
f_R^{-1}(h_U(R))=h_V(R).
$$
Let $\phi \in h_{\uspec A}(R)$, then $f_R(\phi)=f' \cdot \phi \in h_X(R)$.

Hence if $f_R(\phi) \in h_U(R)$ immediately:
$$
\phi: \uspec R \lra V=f^{-1}(U) \lra \uspec A.
$$
So $f_R(\phi) \in h_U(R)$ if and only if $\phi \in h_V(R)$.

\end{observation}

\medskip

We want to define the notion of an open cover of a functor.

\medskip

\begin{definition} Let $F:\salg \lra \sets$ be a functor.
$F$ is covered by the open subfunctors  $(U_i)_{i \in I}$, if
and only if for any affine superscheme $\uspec A$ 
and map $f: h_{\uspec A} \lra F$ we have that
the fibered product
$h_{\uspec A} \times_F U_i \cong h_{V_i}$ and $(V_i)_{i \in I}$ is an
open cover of $\uspec A$. (For the definition of fibered product
see the Appendix \ref{sgappendix}).
\end{definition}

Notice that by the very definition of open subfunctor the
functor $h_{\uspec A} \times_F U_i$ is always representable.
In fact it is equal to $f^{-1}(U)$ which is by definition
the functor of points of an open and affine $V_i$ in ${\uspec A}$.

\medskip

Before going to our main result we need the notion of local 
functor.

\begin{definition}
A functor
$$
F: \salg \lra \sets
$$
is called \textit{ local} or
\textit{ sheaf in the Zariski topology}, if
for each $A \in \salg$,
there exists $f_i\in A_0$, $i \in I$, $(f_i, i \in I)=(1)$,
such that for every collection of $\al_i \in
F(A_{f_i})$ which map to the same element in $F(A_{f_if_j})$, then
there  exists a unique $\alpha \in F(A)$ mapping to each $\al_i$.
\end{definition}

\medskip

\begin{proposition} \label{local}
The functor of points $h_X$  of a superscheme $X$
is local.
\end{proposition} 

\begin{proof} We briefly sketch the proof since it is the same as
in the ordinary case.
Let the notation as in the previous definition.
Consider a collection of maps $\al_i \in h_X(A_{f_i})$ which map to
the same element in $h_X(A_{f_if_j})$. 
Each $\al_i$ consists of two maps: $|\al_i|: \spec {A_{f_i}}_0 \lra |X|$
and a family of 
$\al_{i,U}^*: \cO_X(U) \lra \cO_{A_{f_i}}(|\al|^{-1}(U))$.
The fact the $|\al_i|$ glue together is clear.
The gluing of the $\al_i^*$'s to give $\al:\uspec A \lra X$
depends on the fact that $\cO_X$, $\cO_{A_{f}}$
are sheaves.  
\end{proof}

\medskip

We are ready to state the result that characterizes among all
the functors from $\salg$ to $\sets$ those which are the functors of points
of superschemes.

\medskip

\begin{theorem} \label{representability}
{A functor
$$
F: \salg \lra \sets
$$
is the functor of points of a superscheme $X$, i. e. 
$F=h_X$ if and only if

1. $F$ is local.

2. $F$ admits a cover by affine open subfunctors.
}
\end{theorem}

\begin{proof} Again the proof of this result is 
similar to that in the ordinary case.
We include a sketch of it for lack of an
appropriate reference.
We first observe that if $h_X$ is the functor of points of a
superscheme, by $\ref{local}$ it is local and by $\ref{opencover}$
it admits a cover by open affine subfunctors.


Let's now assume to have $F$ satisfying the properties (1) and (2) of
\ref{representability}.
We need to construct a superscheme $X=(|X|, \cO_X)$
such that $h_X=F$.
The construction of the topological space $|X|$ is the same
as in the ordinary case. Let's sketch it.

Let $\{h_{X_\al}\}_{\al \in \rA}$
be the affine open subfunctors that cover $F$.
Define $h_{X_{\al\be}}=h_{X_\al} \times_F h_{X_\be}$.
($X_{\al\be}$ will correspond
to the intersection of the two open affine $X_\al$ and $X_\be$ in the
superscheme $X$).
Notice that $h_{X_\al} \times_F h_{X_\be}$ is representable.

We have the commutative diagram:
\begin{equation*}
\begin{CD}
h_{X_{\al\be}}=h_{X_\al} \times_F h_{X_\be}  @> j_{\be,\al} >>  h_{X_\be} \\
@VV{j_{\al,\be}}V @VV{i_\be}V\\
h_{X_\al} @> i_\al >>F
\end{CD}
\label{cd}
\end{equation*}

As a set we define:
$$
|X|=_{\defi} \coprod_{\al} |X_\al| / \sim
$$
where $\sim$ is the following relation:
$$
\forall x_\al \in |X_\al|, x_\be \in |X_\be|,
x_\al \sim x_{\be}  \Longleftrightarrow  
\exists x_{\al\be} \in |X_{\al\be}|,  
j_{\al,\be}(x_{\al\be})=x_\al, j_{\be,\al}(x_{\al\be})=x_{\be}.
$$
This is an equivalence relation. $|X|$ is a topological space and
$\pi_\al:|X_\al| \hookrightarrow |X|$ is an injective map.

\medskip

We now need to define a sheaf of superalgebras $\cO_X$, 
by using the sheaves in the open affine $X_\al$ and ``gluing''.

\medskip

Let $U$ be open in $|X|$ and let $U_\al=\pi_\al^{-1}(U)$. Define:
$$
\cO_X(U)=_{\defi}\{(f_\al) \in \coprod_{\alpha \in I} \cO_{X_\al}(U_\al)|
\quad j_{\be,\gamma}^*(f_\be)=j_{\gamma,\be}^*(f_\gamma),
\forall \be, \gamma \in I \}.
$$
The condition $j_{\be,\gamma}^*(f_\be)=j_{\gamma,\be}^*(f_\gamma)$ simply
states that to be an element of  $\cO_X(U)$, the collection 
$\{f_\al\}$ must be
such that $f_\be$ and $f_\gamma$ agree on the intersection of $X_\be$
and $X_\gamma$ for any $\be$ and $\gamma$.

\medskip

One can check that $\cO_X$ is a sheaf of superalgebras.

\medskip

We now
need to show $h_X \cong F$. We are looking for a
functorial bijection
$h_X(A)=\Hom_{\sschemes}
(\uspec A , X )= F(A)$,
for all $ A \in \salg$. 
It is here that we use the hypothesis
of $F$ being local.

To simplify the notation let $T=\uspec A$. We also write
$h_X(T)$ instead of $h_X(A)$. 
So we want to show $h_X(T) \cong F(T)$.
\medskip

We first construct a natural transformation
$\rho_T:F( T ) \lra h_X( T )$.

\medskip

Let $t \in F( T)=\Hom(h_ T,F)$, by Yoneda's lemma.
Consider the diagram:
$$
\begin{CD}
h_{ T_\al}=_{\defi}h_{X_\al} \times_F h_{ T}
@>  >> h_{ T} \\
@VV{ t_\al}V @VV{ t}V\\
h_{X_\al} @> i_\al >> F.
\end{CD}
\label{cd}
$$

Notice that $\{T_\al\}$ form an open affine cover of $T$. 
Since by Yoneda's lemma:
$\Hom(h_{ T_\al},h_{X_\al}) \cong \Hom( T_\al,X_\al)$
we obtain a map: $ t_\al: T_\al \lra X_\al \subset X$. One can check that
the maps $ t_\al$ glue together to give a map $ t':T \lra X$, hence
$ t' \in h_X( T)$. So we define $\rho_ T( t)= t'$. 

\medskip

Next we construct another natural transformation
$\s_T:h_X( T) \lra F( T)$, which turns out to be the inverse of $\rho$.

\medskip
Assume we have $f \in h_X( T )$ i.e.
$f: T \lra X$. Let $ T_\al=f^{-1}(X_\al )$.
We immediately obtain maps $ g_\al :T_\al \lra X_\al \subset F$.
By Yoneda's lemma, $g_\al$ corresponds to a map 
$g_\al:h_{T_\al} \lra h_{X_\al}$. 
Since $F$ is local,
the maps $i_\al \cdot g_\al$ glue together to give a map 
$ g: h_T \lra F$, i. e. an element $g \in F(T)$.  
Define $\s(f)= g$.

\medskip
One can directly check that 
$\rho$ and $\sigma$ are inverse to each other and that the
given correspondence is functorial.
\end{proof}

\medskip

This theorem has an important corollary.

\medskip

\begin{corollary} Fibered products exist in the category of
superschemes. The fibered product $X \times_  Z Y$, for
superschemes $X$, $Y$, $  Z$ with morphisms $f:X \lra Z$,
$g:Y \lra Z$
is the superscheme whose
functor of points is $h_X \times_{h_  Z} h_Y$.
\end{corollary}

\begin{proof}  The proof follows the classical proof, and full
details can be found for example in \cite{dg} I \S 1, 5.1. For completeness
we will briefly sketch the argument. 
Let $F=X \times_  Z Y$. We want to show $F$ is representable. 
One can check that $F$ is local. We then want to show it
can be covered by open affine subfunctors. Let $\{Z_i\}$
be a cover by affine open subsuperschemes of $Z$. 
Define $X_i=X \times_Z Z_i=f^{-1}(Z_i)$ and
$Y_i=Y \times_Z Z_i=g^{-1}(Z_i)$. Let $X_{i\al}$ and
$Y_{j\be}$ open affine covers of $X_i$ and $Y_j$ respectively.
One can check $X_{i\al} \times_{Z_i} Y_{i\be}$ form an affine open
cover of $F$. Hence $F$ is representable.
\end{proof}

\medskip

\begin{remark} One could also
prove directly the existence of fibered product in
the category of superschemes. This is done
exactly as in the classical case, see for
example  Theorem 3.3
in chapter II of \cite{Harts}.
\end{remark}

\bigskip

\begin{remark}
Theorem \ref{representability} can be stated also in the 
$C^\infty$ category:

\medskip 

\noindent
\textit{Let $F$ be a functor $F: \hbox{{\text (smfld)}} \lra \sets$, 
such that
when restricted to the category of manifolds is representable.}

\noindent  \textit{Then the functor $F$ is representable if and only if:}

\noindent \textit{1. $F$ is local, i. e. it has the sheaf property.}

\noindent 2. \textit{$F$ is covered by  open supermanifold functors.}

\medskip

where an open supermanifold functors is a subfunctor $U$ of $F$
such that for all $f:h_X \lra F$, $f^{-1}(U)=h_V$ where $V$ is
a submanifold of $X$ (here $h_X$ denotes the functor of points
of the supermanifold $X$).

\medskip

The proof of this result in the $C^\infty$ category
is essentially the same as the one seen in the algebraic category.

\end{remark}

\bigskip


\section{The Grassmannian superscheme}

In this section we
want to discuss the grassmannian of the $r|s$-dimensional
superspaces inside a super vector space of dimension $m|n$, $r<m$, $s<n$.
We will show that it is a superscheme using the
Theorem \ref{representability}. 
This is a particularly important example since it is the first 
non trivial example of a non affine superscheme.

\medskip

Consider the functor $Gr: \salg \lra \sets$, where for any
superalgebra $A$, 
$Gr(A)$ is the set of projective $A$-submodules of
rank $r|s$ of $A^{m|n}$ (for the definition of the rank of a
projective $A$-module see the Appendix \ref{sgappendix}).

Equivalently $Gr(A)$ can also be defined as:
$$
Gr(A)=\{\al: A^{n|m} \lra L \quad | \quad \al \hbox{ surjective, 
$L$ projective $A$-module of rank $r|s$}\}
$$
(modulo equivalence).

We need also to specify $Gr$ on morphisms $\psi: A \lra B$.

Given  a morphism $\psi:A\rightarrow B$ of superalgebras, we can give to
$B$ the structure of right $A$-module by setting
$$
a\cdot b=\psi(a)b,\qquad a\in A,\;b\in B.
$$ 
Also, given an  $A$-module
$L$, we can construct the $B$-module $L\otimes_AB$. So given
$\psi$ and the element of $Gr(A)$, $f:A^{m|n}\rightarrow L$, we
have an element of $Gr(B)$,
$$
Gr(\psi)(f):B^{m|n}=A^{m|n}\otimes_A B\rightarrow L\otimes_AB.
$$

We want to show that $Gr$ is the
functor of points of a superscheme.

\medskip

We will start by showing it admits a cover of open affine
subfunctors. 
Consider the multiindex $I=(i_1, \dots, i_r| \mu_1,
\dots, \mu_s)$ and the map $\phi_I: A^{r|s} \lra A^{m|n}$ where
$\phi_I(x_1,\dots x_r|\xi_1, \dots \xi_s)$ is the
$m|n$-uple with $x_1, \dots x_r$ occupying the position 
$i_1,\dots, i_r$ and $\xi_1, \dots \xi_s$ occupying the position 
$\mu_1, \dots, \mu_s$ and the other positions are occupied by
zero.
For example, let $m=n=2$ and $r=s=1$. Then
$\phi_{1|2}(x,\xi)=(x,0|0, \xi)$.

\medskip

Now define the subfunctors $v_I$ of $Gr$ as follows.
The $v_I(A)$ are the maps $\al: A^{m|n} \lra L$ such that
$\al \cdot \phi_I$ is invertible. 

\medskip

We want to show that the $v_I$ are open affine subfunctors
of $Gr$. 
The condition that $v_I$ is an open subfunctor is equivalent
to asking that $f^{-1}(v_I)$ is open for any map 
$f: \spec A \lra Gr$.

By Yoneda's lemma, a map $f:\spec A \lra Gr$ corresponds to
a point $f$ in $Gr(A)$. So we are asking if there exists an
open subscheme $V_I$ in $\spec A$, such that
$$
h_{V_I}(B)=\{\psi:A \lra B | Gr(\psi)(f) \in v_I(B)\} 
\subset h_{\uspec A}(B)
$$

To show $V_I$ is open, consider the matrix 
$$
Z=(f(e_{i_1}) \dots f(e_{i_r}), f(\ep_{\iota_1}) \dots f(e_{\iota_s}))
$$
and define $b_A(f)$ the product of the determinants of the two
even diagonal blocks of $Z$.

If  $b_A(f)$ is invertible, then any map
$\psi: A \lra B$ is forced to send $b_A(f)$ into an
invertible element in $B$, hence all maps $Gr(\psi)(f)$ are in $v_I(B)$.
Hence $V_I=\uspec A$.

\medskip

If $b_A(f)$ is zero, then no map can send $b_A(f)$ into
an invertible element, so ${V_I}$ is 
empty.

\medskip

The only non trivial case is when $b_A(f)$ is non-zero and
not invertible. In this case since
$b_A(f)$ is sent to an invertible element in $B$ by $\psi$ we
have a one to one correspondence between such maps
$\psi$ and $\psi': A[b_A(f)^{-1}] \lra B$.
So we have obtained that $h_{V_I} \cong \uspec A[b_A(f)^{-1}]$
which is open in  $\uspec A$.

\medskip

It remains to show that these subfunctors cover $Gr$.

Given $f \in Gr(A)$, that
is a function from $h_{\uspec A}  \lra Gr$, we have that since
$f$ is surjective, there exists at least an index $I$ so that
$b_A(f)$ is invertible, hence $f \in v_I(A)$ for this $I$. 
The above argument shows that we obtain a cover of 
any $\uspec A$ by taking ${v_I} \times_{Gr} h_{\uspec A}$.
\medskip

Finally we want to show that $Gr$ is local. This is immediate once
we identify $Gr(A)$ with coherent sheaves with locally constant
of rank $r|s$. 

$$
Gr(A) \cong \{\cF \subset \cO_A^{m|n}\,/\, \cF \hbox{ is a subsheaf,
of locally costant rank } r|s \}
$$
where $\cO_A^{m|n}=k^{m|n} \otimes \cO_A$.

\medskip

By its very definition this functor is local. 

\medskip

This identification is possible since by the Appendix \ref{projective}
we prove that a projective module $M$ is locally free
and the correspondence between coherent sheaves and finitely
generated modules in the supersetting.

\medskip

So we have shown that $Gr$ is the functor of points of
a superscheme that we will call the \textit{ supergrassmannian} of
$r|s$ subspaces into a $m|n$ dimensional space.

\bigskip

\section{The infinitesimal theory}
\bigskip
In this section we discuss the infinitesimal theory
of superschemes.
We define the notion of tangent
 space to a superscheme and to a supervariety at a point of
the underlying topological space.
We then use these definitions in explicit computations. 
\medskip


\medskip

Let $k$ be a field.

\medskip


\medskip

\begin{definition} Let $ X =(X, \cO_X)$ be a superscheme
(supervariety). We say that
$ X $ is \textit{ algebraic} if it admits an open affine finite
cover $\{ X _i\}_{i \in I}$ 
such that $\cO_X(X_i)$ is a finitely generated
superalgebra for each $X_i$.
\end{definition}
\medskip

Unless otherwise specified all superschemes are assumed to be algebraic.

\medskip

Given a superscheme $ X =(|X|,\cO_X)$ each
point of $x$ in the topological space $|X|$ belongs
to an open affine subsuperscheme $\uspec A$, $x \cong \bp
\in \spec A_0$,
so that $\cO_{X,x} \cong A_{\bp}$. 
Recall that $A_\bp$ is the
localization of the $A_0$-module $A$ into the prime ideal 
$\bp \subset A_0$ and that
$$
A_{\bp}=\{ {f \over g} | g \in A_0-\bp\}.
$$

The local ring $A_\bp$ contains the maximal ideal
$\bp_{x}$ generated by the maximal ideal in the local ring
$(A_\bp)_0$ and the generators of $(A_\bp)_1$ as $A_0$-module.

\medskip

We want to define the notion of a rational
point of a scheme. We will then define
the tangent space to a scheme in a rational
point. 


\medskip

\begin{definition} Let $ X =(|X|,\cO_X)$ be a superscheme.
A point $x \in |X|$ is said to be \textit{ rational} if
$\cO_{X,x}/\bp_x \cong k$.
\end{definition}

\medskip

\begin{remark}
As in the commutative case we have that
if $k$ is algebraically closed all closed points of $|X|$ are
rational. This is because the field $\cO_{X,x}/\bp_x$ is a finite algebraic
extension of $k$. 
(see \cite{am} 7.9 for more details).

(Recall that a point $x \in |X|$ is  \textit{ closed} if it corresponds to
a maximal ideal in $\spec A_0$, where
$(\spec A_0, \cO_A) \subset  X $ is any affine 
open neighbourhood of $x$).
\end{remark}

\medskip

It is important not to confuse the points of the
underlying topological space $|X|$ with the elements obtained via
 the functor of points, $h_X(A)$ for a generic $A \in \salg$. These
 are called \textit{ $A$-points} of the superscheme $X$. The next
 observation clarifies the relationship between the points of $X$
and $h_ X$ the functor of points of $X$.

\medskip

\begin{observation} 
There is a bijection between the rational points of a superscheme $ X $
and the set of its $k$-points $h_ X (k)$.
In fact, an element $(|f|,f^*) \in h_ X (k)$,
$|f|:\spec k \lra |X|$, $f^*:\cO_{X,x} \lra k$, determines
immediately a point $x=|f|(0)$, which is rational. 
\end{observation}
\medskip

\begin{definition} Let $A$ be a superalgebra and $M$
an $A$-module. Let $D:A \lra M$ be an additive map
with the property $D(a)=0$, $\forall a \in k$.
We say that $D$ is an \textit{ super derivation} if:
$$
D(fg)=D(f)g+(-1)^{p(D)p(f)}fD(g), \qquad f,g \in A
$$
where $p$ as always denotes the parity.
\end{definition}


\medskip

\begin{definition} Let $ X =(|X|, \cO_ X )$
be a superscheme and $x$ a  rational point in $|X|$.
We define \textit{ tangent space of $ X $ at $x$}:
$$
T_x X =\Der(\cO_{X,x}, k)
$$
where $k$ is viewed as $\cO_{X,x}$-module via the identification
$k \cong \cO_{X,x}/\bm_{X,x}$, where $\bm_{X,x}$ is the 
maximal ideal in $\cO_{X,x}$.
\end{definition}

\medskip

The next proposition gives an equivalent definition for the
tangent space.

\medskip

\begin{proposition} Let $ X $ be a superscheme, then: \end{proposition}
$$  
T_x X =\Der(\cO_{X,x}, k) \cong \uHom_{\smod}(\bm_{X,x}/\bm^2_{X,x},k).
$$

\medskip

Note that $\bm_{X,x}/\bm^2_{X,x}$ is a $\cO_{X,x}$-supermodule which is
annihilated by $\bm_{X,x}$, hence it is a $k=\cO_{X,x}/\bm_{X,x}$-supermodule
i.e. a super vector space.

\medskip

\begin{proof} 
Let $D \in \Der(\cO_{X,x},k)$. Since $D$ is zero on $k$
and $\cO_{X,x}=k \oplus \bm_{X,x}$
we have that $D$ is determined by 
its restriction to $\bm_{X,x}$, $D|_{\bm_{X,x}}$. 
Moreover since $\bm_{X,x}$ acts as
zero on $k \cong \cO_{X,x}/\bm_{X,x}$ one can check that
$$
\begin{array}{ccc}
\psi:\Der(\cO_{X,x},k) & \lra & \uHom_{\smod}(\bm_{X,x}/\bm^2_{X,x},k)
\\
D & \mapsto & D|_{\bm_{X,x}}
\end{array}
$$
is well defined.

Now we construct the inverse. Let $\al: \bm_{X,x} \lra k$,
$\al(\bm^2_{X,x})=0$. Define
$$  D_\al: \cO_{X,x}=k \oplus \bm_{X,x} \lra k, \qquad
D_\al (a,f)=\al(f).
$$
This is a well defined superderivation.

Moreover one can check that the
map $\al \mapsto D_\al$ is $\psi^{-1}$. 
\end{proof}

\medskip

The next proposition provides a characterization of
the tangent space, that is useful for explicit calculations.

\medskip
\begin{proposition} \label{tangentspace}
Let $ X =(|X|,\cO_X)$ be a supervariety
$x \in |X|$ a rational closed point. Let $U$ be an affine neighbourhood
of $x$, $\bm_x \subset k[U]$ the maximal ideal corresponding to $x$.
Then
\end{proposition}
$$
T_x X  \cong \uHom_{\smod}(\bm_{X,x}/\bm^2_{X,x},k) \cong
\uHom_{\smod}(\bm_{x}/\bm^2_{x},k).
$$

\begin{proof}
The proof is the same as in the ordinary case
and is based on the fact that localization commutes with exact sequences.
\end{proof}

\medskip
Let's compute explicitly the tangent space in an example.
\medskip

\begin{example}\label{extangent}
Consider the affine supervariety represented by the coordinate ring:
$$
\C[x,y,\xi,\eta]/(x\xi+y\eta).
$$
Since $\C$ is algebraically closed, all closed points are rational.
Consider the closed
point $P=(1,1,0,0) \cong \bm_P=(x-1,y-1,\xi,\eta) \subset
\C[x,y,\xi,\eta]/(x\xi+y\eta)$ (we identify $(x_0,y_0,0,0)$ with
maximal ideals in the ring of the supervariety, as we do in the
commutative case).
By Proposition \ref{tangentspace},  
the tangent space at $P$ is given by all the functions
$\al:\bm_P \lra k$, $\al(\bm_P^2)=0$.

A generic $f \in \bm_P$ lifts to the family of $f=f_1+f_2(x\xi+y\eta)
\in \C[\bA^{1|1}]=\C[x,y,\xi,\eta]$ with
$f_1(1,1,0,0)=0$ and where $f_2$ is any function in $\C[\bA^{1|1}]=
\C[x,y,\xi,\eta]$.
Thus $f$ can be formally expanded in power series around $P$
(see \cite{VSV2} for more details).
$$
\begin{array}{c}
f={\del f_1 \over \del x}(P) (x-1)+ 
{\del f_1 \over \del y}(P) (y-1)  +
({\del f_1 \over \del \xi}(P)+f_2(P))\xi+ \\ \\ 
 ({\del f_1 \over \del \eta}(P)+f_2(P))\eta+ 
 \hbox{ higher order terms}.
\end{array}
$$
Define:
$$
X={\del f_1 \over \del x}(P), \quad
Y={\del f_1 \over \del y}(P), \quad
\Xi={\del f_1 \over \del \xi}(P), \quad
\Eta={\del f_1 \over \del \eta}(P).
$$
These are coordinates for the supervector space $\btm_P/\btm_P^2$,
$\btm_P=(x-1,y-1,\xi,\eta) \subset \C[x,y,\xi,\eta]$.
A basis for the dual space $(\btm_p/ \btm_p^2)^*$
consists of sending the coefficient
of one of the $x-1, y-1, \xi, \eta$ to a non zero element and the
others to zero. This gives relations that
allow us to elimine the parameter $f_2(P)$. We get equations:
$$
\quad \Xi+f_2(P)=0, \quad \Eta+f_2(P)=0.
$$
Eliminating the parameter we get the equation for the tangent space:
$$
\Xi+\Eta=0.
$$

So we have described the tangent space $(\bm_P/\bm_P^2)^*$ as
a subspace of $(\btm_P/\btm_P^2)^*$, the tangent space to the affine
superspace $\bA^{m|n}$.
\end{example}
\medskip

There is yet another way to compute the tangent space, in the
case $X$ is an affine supervariety. Before we examine this
construction, we must understand first 
the notion of differential of a function and of a
morphism.

\medskip

\begin{definition} Let $ X =(|X|, \cO_X)$ be a superscheme,
$x$ a rational point.

Consider the projections:
$$
\pi:\cO_{X,x} \lra \cO_{X,x} /\bm_{X,x} \cong k, \qquad
p:\bm_{X,x} \lra \bm_{X,x}/\bm^2_{X,x}
$$
Let $f \in \cO_{X,x}$, 
we define \textit{value of $f$ at $x$}:
$$
f(x) =_{\defi} \pi(f).
$$
We also define \textit{differential of $f$ at $x$}:
$$
(df)_x=_{\defi} p(f-f(x)).
$$

We now want to define value and differential of a section
in a point.

If $U$ is an open neighbourhood of $x$ and $f \in \cO_X(U)$
we define \textit{ value of $f$ at $x$} to be:
$$
f(x)=_{\defi} \pi(\phi(f)), \qquad \phi: \cO_X(U) \lra \cO_{X,x}.
$$ 
We define \textit{ differential of $f$ at $x$}
$$
(df)_x=_{\defi} (d\phi(f))_x.
$$
\medskip

For example, if $P=(x_1^0 \dots x_m^0, 0 \dots 0)$ is a
closed rational point of the affine superspace $\bA^{m|n}$ with coordinate ring
$k[x_1 \dots x_m, \xi_1 \dots \xi_n]$ a basis of $\bm_P/\bm^2_P$ is
$\{x-x_i^0, \xi_j\}_{\{i=1 \dots m,j=1 \dots n\}}$. Hence:
$$
(dx_i)_P=x-x_i^0,
\qquad (d\xi_j)_P=\xi_j, \qquad i=1 \dots m,\quad j=1 \dots n
$$
\medskip
Let $(|\al|,\al^*): X  \lra  Y$ be a morphism of 
superschemes and $x$ a rational
point, $|\al|(x)$ also rational, $\al$
induces a morphism $d\al_x:T_x  X  \lra T_{|\al|(x)}  Y$, by:
$$
\begin{array}{c}
d\al_x(D)f=D(\al_x^*(f)), \quad D \in T_x  X =\Der(\cO_{X,x},k),\quad
\al^*_x:\cO_{Y,|\al|x} \lra \cO_{X,x}
\end{array}
$$
\end{definition}

It is simple to check that if $(|\al|,\al^*)$ is an immersion
i.e. it identifies $ X $ with a subscheme of $ Y$, then $d\al_x$ is
injective. Hence if $ X $ is a subsupervariety of $\bA^{m|n}$ it makes sense
to ask for equations that determine the tangent superspace to $ X $ as
a linear subsuperspace of $T_x\bA^{m|n} \cong k^{m|n}$.

\medskip

\begin{proposition} \label{tangent}
Let $ X $ be a subvariety of $\bA^{m|n}$
and let $x$ be a rational closed point of $ X $. Then
$$
T_x X =\{v \in k^{m|n}| (df)_x(v)=0, \forall f \in I\}
$$
where $I$ is the ideal defining $X$ in $k[\bA^{m|n}]$.
\end{proposition}

\begin{proof} The immersion $\al: X  \subset \bA^{m|n}$ corresponds to
a surjective morphism $\phi:k[\bA^{m|n}] \lra k[ X ]$, hence $k[ X ] \cong
k[\bA^{m|n}]/I$.
Let $\bm_x$ and $\btm_x$ denote respectively the maximal ideal
associated to $x$ in $X$ and $\spec k[\bA^{m|n}]$ respectively.
$\phi$ induces a surjective linear map $\psi$ between superspaces:
$$
\psi: \btm_x/\btm_x^2 \lra \bm_x/\bm_x^2.
$$
Let's recall the following simple fact of linear algebra.
\medskip
If $a:V_1 \lra V_2$ is a surjective linear map between finite
dimensional vector spaces $V_1$, $V_2$ and $b:V_2^* \subset V_1^*$
is the injective linear map induced by $a$ on the dual vector spaces
then $s \in Im(b)$ if and only if $s|_{ker(a)}=0$.
\medskip
We apply this to the differential
$$
(d\al)_x:T_x( X )=(\bm_x/\bm_x^2)^* \lra
T_{\al(x)}(\bA^{m|n})=(\btm_x/\btm_x^2)^*
$$
and we see that
$$
T_x( X )=\{v \in T_{\al(x)}(\bA^{m|n}) |
v(ker(\psi))=0\}.
$$
Observe that $ker(\psi)=\{(df)_x| f \in I\}$.
By identifying $\bA^{m|n}=k^{m|n}$ with its double dual
$(k^{m|n})^{**}$ we obtain the result.
\end{proof}

\medskip

\begin{remark} In the notation of the previous proposition,
if $I=(f_1 \dots f_r)$ one can check that:
$$
T_x X =\{v \in k^{m|n}| (df_i)_x(v)=0, \forall i=1 \dots r \}.
$$
\end{remark}

\medskip

Let's revisit Example \ref{extangent} and see 
how the calculation is made using
Proposition \ref{tangent}.

\medskip

\begin{example} Consider again the supervariety represented by:  
$$
\C[x,y,\xi,\eta]/(x\xi+y\eta).
$$
We want to compute the tangent space at $P=(1,1,0,0)=(x_0,y_0,\xi_0,\eta_0)$.
$$
\begin{array}{cc}
d(x\xi+y\eta)_P& =x_0(d\xi)_P+\xi_0 (dx)_P+y_0(d\eta)_P+\eta_0 (dy)_P \\
&=(d\xi)_P+(d\eta)_P
\cong(0,0,1,1).
\end{array}
$$
Hence by \ref{tangent}
the tangent space is the subspace of $k^{2|2}$ given by the equation:
$$
\xi+\eta=0.
$$
 
\end{example}

\chapter{Algebraic Supergroups and Lie Superalgebras}\label{sgchap6}

In this section we introduce the notion of supergroup scheme,
and of its Lie superalgebra. 

\medskip

Let $k$ be a noetherian ring.

\medskip

All superschemes are assumed to be algebraic.

\bigskip

\section{Supergroup functors and algebraic supergroups}

A supergroup scheme is a superscheme whose functor of points is group
valued. In order to study supergroup
schemes we need first to understand 
the weaker notion of supergroup functor.

\medskip

\begin{definition}  A \textit{ supergroup functor} is a group
valued functor:
$$
 G :\salg \lra \sets
$$
\end{definition}

\begin{remark}
Saying that $G$ is group valued
is equivalent to have the following natural transformations:

1. Multiplication $\mu:  G  \times  G  \lra  G $, such that
$\mu \circ (\mu \times \id)=(\mu \times \id) \circ \mu$, i. e.
$$
\begin{CD}
G \times G \times G @> \mu \times \id >> G \times G \\
@V{\id \times \mu}VV @VV \mu V \\
  G \times G @> \mu >> G  \\
\end{CD}
$$
2. Unit $e:e_k \lra  G $, where $e_k:\salg \lra \sets$,
$e_k(A)=1_A$, such that $\mu \circ (\id \otimes e)=
\mu \circ (e \times \id)$, i. e.
$$
\begin{array}{ccccc}
G \times e_k & \stackrel{\id \times e} \lra &
G \times G & \stackrel{e \times id} \longleftarrow & e_k \times G \\
             & _{}\searrow & _{\mu}\Big\downarrow & \swarrow  _{}& \\
             &       & G & & \\
\end{array}
$$

3. Inverse 
$i: G  \lra  G $, such that $\mu \circ (\id \times i)=e \circ \id$,
i. e. 
$$
\begin{CD}
G @> (\id, i)>> G \times G \\
@V{}VV @VV \mu V \\
  e_k @>e>> G  \\
\end{CD}
$$

\end{remark}

\medskip

If $ G $ is the functor of points of a superscheme $X$, we say that
$ X$ is a \textit{ supergroup scheme}.
An \textit{ affine supergroup scheme} is a supergroup
scheme which is an affine superscheme. To make the terminology easier
we will drop the word ``scheme'' when speaking of supergroup
schemes.

\medskip

\begin{observation}
The functor of points of an affine supergroup 
$ G $ is a representable functor. It is represented by the superalgebra
$k[ G ]$. This superalgebra has a Hopf superalgebra structure,
so we identify the category of affine supergroups 
with the category of commutative
Hopf superalgebras.

\end{observation}

\medskip

\begin{observation} 
If $k$ is a field, we may interpret the unit $e$ 
of a supergroup $G=(|G|,\cO_G)$
as a rational point
of $G$ that we will denote with $1_G$. 
In fact $e:e_k \lra G$, $e=(|e|,e^*)$.
Define $1_G=|e|(|e_k|)$. This is a rational point,
in fact $\cO_{G,1_G}/m_{1_G} \cong k$.
Moreover by the very definition of $e$, $1_G$ has
the property of a unit for the group $|G|$.
\end{observation}

\medskip 

\begin{example}

1. \textit{ Supermatrices $\rM_{m|n}$}. It is immediate to verify that the
supermatrices discussed in chapter 1 are an
affine supergroup where $\mu$ is interpreted as the usual 
matrix addition.

\medskip

2. \textit{ The general linear supergroup $\rGL_{m|n}$}.

Let $A \in \salg$.
Define ${\rGL_{m|n}}(A)$ as $\rGL(A^{m|n})$
(see chapter 1) the set of automorphisms of the $A$-supermodule
$A^{m|n}$. Choosing coordinates we can write
$$
\rGL_{m|n}(A)=\left\{ \begin{pmatrix} a & \alpha \\ \beta & b \end{pmatrix} 
\right\}
$$
where $a$ and $b$ are $m \times m$, $n \times n$ blocks of even elements
and $\forall$, $\al$, $\be$ $m \times n$, $n \times m$ blocks of odd elements
and $a$ and $b$ are invertible matrices.

\medskip

It is not hard to see that
this is the functor of points of
an affine supergroup $\rGL(m|n)$ represented by the Hopf
superalgebra
$$
k[\rGL(m|n)]=k[x_{ij},\xi_{kl}][T]/(T \Ber-1),
$$
where $x_{ij}$'s
and $\xi_{kl}$'s are respectively even and odd variables
with
$1 \leq i,j \leq m$ or  $m+1 \leq i,j \leq m+n$,
$1 \leq k \leq m$, $m+1 \leq l \leq m+n$ or
$m+1 \leq k \leq m+n$, $1 \leq l \leq m$
and $\Ber$ denotes the Berezinian.

In general if $V$ is a super vector space we define the
functor $\rGL(V)$ as $\rGL(V)(A)=\rGL(V(A))$, the invertible
transformations of $V(A)$ preserving parity.
\end{example}

3. \textit{ The special linear group $\rSL_{m|n}$}.

For a superalgebra $A$, 
let's define $\rSL_{m|n}(A)$ to be the subset of $\rGL_{m|n}(A)$
consisting of matrices with Berezinian equal to 1. This is the
functor of points of an affine supergroup and it is represented
by the Hopf superalgebra:
$$
k[\rSL(m|n)]=k[x_{ij},\xi_{kl}]/(\Ber-1).
$$

Similarly one can construct the functor of points and
the representing Hopf superalgebras for all the classical algebraic
supergroups (see \cite{DM, VSV2}).

\bigskip

\section{Lie superalgebras}

\medskip

Assume  2 and 3 are not zero divisors in $k$.

\medskip

In this section we define functorially the notion of 
Lie superalgebra. 

Our definition is only apparently different from the one we
have introduced in chapter \ref{sgchap1}, which is the one
mostly used in the literature. 

\medskip

Let $\cO_k:\salg \lra \sets$ be the functor represented by $k[x]$.
$\cO_k$
corresponds to an ordinary algebraic variety,
namely the affine
line. For a superalgebra $A$ we have that $\cO_k(A)=A_0$.

\medskip

\begin{definition} \label{liefunctor}
Let $\fg$ be a free $k$-module. 
We say that the group valued functor
$$
L_{\fg}:\salg \lra \sets, \qquad L_{\fg}(A)=(A \otimes \fg)_0
$$
is a \textit{ Lie superalgebra} if
there is a $\cO_k$-linear natural transformation
$$
[\;,\;]:L_{\fg} \times L_{\fg} \lra L_{\fg}
$$
that satisfies commutative
diagrams corresponding to the antisymmetric property and the
Jacobi identity. For each superalgebra $A$, the bracket $[\;,\;]_A$ defines a
Lie algebra structure on the $A$-module
$L_{\fg}(A)$, hence the functor $L_{\fg}$ is \textit{ Lie
algebra} valued. We will drop the suffix $A$ from the bracket
and the natural transformations to ease the notation.
\end{definition}

\medskip
\begin{remark}
In general, a Lie superalgebra is not representable.
However if $\fg$ is finite dimensional,
$$
L_{\fg}(A)=(A \otimes \fg )_0=
\Hom_{\smod}(\fg^*,A)=\Hom_{\salg}(\sym(\fg^*),A)
$$
where $\smod$ denotes the category of supermodules and
$\sym(\fg^*)$ the symmetric algebra over $\fg^*$. In this special
case $L_{\fg}$ is representable and it
is an affine superscheme represented by the
superalgebra $\sym(\fg^*)$.
\end{remark}
\medskip

We want to see that the usual notion
of Lie superalgebra, as defined by Kac (see \cite{ka}) among many
others, is equivalent to this functorial definition.

Recall that in chapter \ref{sgchap1} we gave the following definition of
Lie superalgebra.

\medskip

Let $k$ be a field, char$(k) \neq 2,3$.

\medskip

\begin{definition} \label{superliealgebra}
Let $ \fg $ be a super vector space. We say that
a bilinear map 
$[,]:  \fg  \times  \fg  \lra  \fg$
is a \textit{ superbracket} if  $\forall x,y,z \in  \fg $:

\noindent
a) $[x,y]=(-1)^{p(x)p(y)}[y,x]$

\noindent
b) $[x,[y,z]]+(-1)^{p(x)p(y)+p(x)p(z)}[y,[z,x]]+(-1)^{p(x)p(z)+p(y)p(z)}
[z,[x,y]]$.

The super vector space
$ \fg $ satisfying the above properties
is commonly defined as Lie superalgebra
in the literature.
\end{definition}

\medskip

\begin{observation}
The definition \ref{liefunctor} and \ref{superliealgebra} are
equivalent. In other words given a Lie algebra valued
functor $L_{\fg}:\salg \lra \sets$ we can build a super vector space
$\fg$ with a superbracket and vice-versa. Let's see these constructions
in detail.

If we have a Lie superalgebra $L_{\fg}$
there is always, by definition, a super vector space
$\fg$ associated to it.
The superbracket on $\fg$ is given following the \textit{ even rules}.
Let's see in detail what it amounts to in this case
(for a complete treatment of even rules see pg 57 \cite{DM}).
\medskip
Given $v, w \in \fg$, since the Lie bracket on $L_\fg(A)$ is
$A$-linear we can define the element $\{v,w\} \in \fg$ as:
$$
[a \otimes v, b \otimes w]=(-1)^{p(b)p(v)}ab \otimes \{v,w\} \in
(A \otimes \fg)_0 \in L_\fg(A).
$$
Clearly the bracket $\{v,w\} \in \fg$ does not depend on $a,b \in A$.
It is straightforward to verify that it is a superbracket. Let's see,
for example, the antisymmetry property.
Observe first that if $a \otimes v \in (A \otimes \fg)_0$
must be $p(v)=p(a)$, since $(\fg \otimes A)_0=$
$A_0 \otimes \fg_0 \oplus A_1 \otimes \fg_1 $.
So we can write:
$$
[a \otimes v, b \otimes w]=(-1)^{p(b)p(v)}\{v,w\} \otimes ab=
(-1)^{p(w)p(v)}ab \otimes \{v,w\}. 
$$
On the other hand:
$$
\begin{array}{c}
[b \otimes w, a \otimes v]=(-1)^{p(a)p(w)}ba \otimes \{w,v\} =
(-1)^{p(a)p(w)+p(a)p(b)} ab \otimes \{w,v\}= \\
\\
(-1)^{2p(w)p(v)}ab \otimes \{w,v\} =\{w,v\}.
\end{array}
$$
By comparing the two expressions we get the antisymmetry of the superbracket.
For the super Jacobi identity the calculation is the same.

\medskip
A similar calculation also shows that given 
a supervector space with a super bracket one obtains a Lie
superalgebra.

\medskip

Hence a Lie superalgebra $L_\fg$ according to
Definition \ref{superliealgebra} 
is equivalent to a supervector space $\fg$
with a superbracket. With an abuse of language we will
refer to both $\fg$ and $L_{\fg}$
as  ``Lie superalgebra''.
\end{observation}

\medskip
\begin{remark}
Given a supervector space $\fg$ one may also define a Lie superalgebra
to be the representable functor $D_\fg:\salg \lra \sets$ so that
$$
D_\fg(A)=\Hom_\smod (\fg^*,A)=
\Hom_\salg (\sym(\fg^*),A)
$$
with a $\cO_k$ linear natural transformation $[,]:D_\fg \times D_\fg
\lra D_\fg$ satisfying the commutative diagrams corresponding
to antisymmetry and Jacobi identity.
When $\fg$ is finite dimensional this definition coincides with the
previous one, however we have preferred the one given in
\ref{superliealgebra} since its
immediate equivalence with the definition mostly used in the literature.
\end{remark}

The purpose of the next two sections is to naturally associate 
a Lie superalgebra $\Lie(G)$ to a supergroup $G$.

\bigskip

\section{$\Lie(G)$ as tangent superspace to a supergroup scheme}

For the rest of this chapter, let $k$ be a field, char$(k) \neq 2,3$.

\medskip

Let $G$ be a supergroup functor.


\medskip


\medskip

Let $A$ be a commutative superalgebra and let
$A(\ep)=_{\defi}A[\ep]/(\ep^2)$ be the algebra of dual numbers ($\ep$
here is taken as an \textit{ even} indeterminate). We have that
$A(\ep)=A \oplus \ep A$ and there are two natural
morphisms: 
$$
\begin{array}{c}
i: A \rightarrow
A(\ep), \qquad i(1)=1 \\ 
p: A(\ep) \rightarrow A, \qquad
p(1)=1, \quad p(\ep)=0, \quad p \cdot i= \id_A
\end{array}
$$
\medskip

\begin{definition} \label{supergroup}
Consider the homomorphism
$ G (p): G (A(\ep)) \longrightarrow  G (A)$. For each $ G $ there is a
supergroup functor, 
$$
\Lie( G ):\salg \lra \sets, \qquad
\Lie( G )(A)=_{\defi}\ker( G (p)).
$$
\end{definition}

If $G$ is a supergroup scheme, we denote $\Lie(h_G)$ by $\Lie(G)$.
\medskip

\begin{example} \label{gl}

1. \textit{ The super general linear algebra.} 

We want to determine the functor $\Lie(\rGL_{m|n})$.
Consider the map:
$$
\begin{array}{cccc}
\rGL_{m|n}(p):  & \rGL_{m|n}(A(\ep)) &  \lra  & \rGL_{m|n}(A) \\ \\
& \begin{pmatrix} p+\ep p' & q+\ep q' \\ r+\ep r' & s+ \ep s' \end{pmatrix}
& \mapsto &
\begin{pmatrix} p & q \\ r & s \end{pmatrix}
\end{array}
$$
with $p,p',s,s'$ having entries in $A_0 $ and $q,q',r,r'$ having
entries in $A_1$; 
the blocks  $p$ and  $s$ are invertible matrices. One can
see immediately that
$$
\Lie(\rGL_{m|n})(A)=\ker(\rGL_{m|n}(p))= \left\{ \begin{pmatrix}
I_m+\ep p' & \ep q' \\ \ep r' & I_n+ \ep s' \end{pmatrix}\right\}
$$
where $I_n$ is a $n \times n$ identity matrix.
The functor $\Lie(\rGL_{m|n})$ is clearly group valued and can be
identified with the (additive) group functor $M_{m|n}$ defined as:
$$
M_{m|n}(A)=\Hom_\smod(M(m|n)^*,A)=\Hom_\salg(\Sym(M(m|n)^*),A)
$$
where $M(m|n)$ is the supervector space
$$
M(m|n)=\left\{ \begin{pmatrix}P & Q \\ R & S
\end{pmatrix} \right\} \cong k^{m^2+n^2|2mn}
$$
where $P$, $Q$, $R$, $S$ are respectively
$m \times m$, $m \times n$, $n \times m$, $n \times n$
matrices with entries in $k$.
An element $X \in M(m|n)$ is even if $Q=R=0$, odd if $P=S=0$.

Notice that $M(m|n)$ is a Lie superalgebra with superbracket:
$$
[X,Y]=XY-(-1)^{p(X)p(Y)}YX
$$
So $\Lie(\rGL_{m|n})$ is a Lie superalgebra. In the next section
we will see that in general we can give a Lie superalgebra
structure to $\Lie(G)$ for any group scheme $G$.
\medskip
\end{example}


\medskip

2. \textit{ The special linear superalgebra. }

A similar computation shows that 
$$
\Lie(\rSL_{m|n})(A)=\left\{ W=\begin{pmatrix}
I_m+\ep p' & \ep q' \\ \ep r' & I_n+ \ep s' \end{pmatrix}
 \quad | \quad \Ber(W)=1 \right\}.
$$
The condition on the Berezinian is equivalent to:
$$
\det(1-\ep s') \det(1+\ep p')=1
$$
which gives:
$$
\mathrm{tr}(p')-\mathrm{tr}(s')=0.
$$
Hence
$$
\Lie(\rSL_{m|n})(A)=\{ X \in M_{m|n}(A) \quad | \quad 
\mathrm{Tr}(X)=0\}.
$$

Similar calculations can be done also for the other classical
supergroups.

\medskip

Let's now assume that $G$ is a supergroup scheme.

We now want to show that $\Lie(G):\salg \lra \sets$ 
is a representable functor and its representing superscheme is
identified with the tangent space at the identity of the supergroup
$G$.

\medskip

\begin{definition} \label{firstnbhd}
Let $X=(|X|,\cO_X)$ be a superscheme, $x \in |X|$. We define
\textit{ the first neighbourhood of $X$},
$X_x$, to be the superscheme
$\uspec \cO_{X,x}/m_{X,x}^2$.
The topological space $|X_x|$ consists of the one point $m_{X,x}$
which is the maximal ideal in $\cO_{X,x}$.
\end{definition}

\medskip

\begin{observation} There exists a natural map
$f:X_x \lra X$. In fact we can write immediately
$$
\begin{array}{ccc}
|f|: \spec \cO_{X,x}/m_{X,x}^2=\{m_{X,x}\} & \lra & X \\
m_{X,x} & \mapsto & x
\end{array}
$$
$$
\begin{array}{ccc}
f^*_{U}:  \cO_{X}(U) & \lra & \cO_{X_x}(|f|^{-1}(U))=
\cO_{X,x}/m_{X,x}^2 \\
\end{array}
$$
where $f^*_{U}$ is the composition of natural map 
from $\cO_{X}(U)$ to the direct
limit $\cO_{X,x}$ and the projection $\cO_{X,x} \lra \cO_{X,x}/m_{X,x}^2$.

\end{observation}

\medskip


We now want to make some observation on the identity element
of a supergroup $G$. By definition we have that the identity
is a map $e: \uspec k \lra G$. This corresponds to a map of
the functor of points: $h_{\uspec k} \lra h_G$ assigning to
the only map $1_A \in h_{\uspec k}(A)$ a map that we will
denote $1_{G(A)} \in h_G(A)= \Hom(\uspec A, G)$. 
The topological space map $|1_{G(A)}|$ sends all the maximal
ideals in $\spec A$ to $1_G \in |G|$.
The sheaf map $\cO_G \lra k$ is the evaluation at $1_G$
that is $\cO_G(U) \lra \cO_{G,1_G} \lra \cO_{G,1_G}/m_{G,1_G} \cong k$
(the identity is a rational point).
Hence it is immediate to verify that $1_{G(A)}$ \textit{ factors
through} $G_1$ (the first neighbourhood at the identity $1_G$), i. e.
$1_{G(A)}: \uspec A \lra G_1 \lra G$. This fact will be
crucial in the proof of the next theorem.

\medskip

\begin{theorem} \label{tangentspace}  Let  $ G $ be an
algebraic supergroup.
Then 
$$
\Lie( G )(A)=\Hom_\smod({m_{G,1_G}/m_{G,1_G}^2},A)=
(A \otimes T_{1_G}( G ))_0
$$

where $T_{1_G}( G )$ denotes the tangent space in the rational point $1_G$.
\end{theorem}

\begin{proof}
Let $d:m_{1_G}/m_{1_G}^2 \lra A$ be a linear map.
Let $d'$ be the map:
$$
\begin{array}{ccccc}
\cO_{G,1_G}/m_{1_G}^2 \cong
k \oplus m_{1_G}/m_{1_G}^2 & \lra & A(\ep)
\\  \\
(s,t) & \mapsto & s+d(t)\ep.
\end{array}
$$
So we have $d' \in h_{G_1}(A(\ep))$ since $G_1$ is a superscheme
represented by $\cO_{G,1_G}/m_{1_G}^2$.

This shows that we have a correspondence between
$h_{G_1}(A)$ and the elements of $\Hom_\smod(m_{1_G}/m_{1_G}^2,A)$.
Let $\phi:G_1 \lra G$ be the
map described in \ref{firstnbhd}. By Yoneda's lemma
$\phi$ induces $\phi_{A(\ep)}: h_{G_1}(A(\ep)) \lra h_G(A(\ep))$,
hence we have a map 
$$
\begin{array}{ccc}
\psi: \Hom_\smod(m_{1_G}/m_{1_G}^2,A) & \lra & h_G(A(\ep)) \\
d & \mapsto & \phi_{A(\ep)}(d'). \\
\end{array}
$$
The following commutative diagram shows that
$\psi(d) \in ker(h_{G}(p))=\Lie(G)(A)$. 
$$
\begin{array}{ccc}
h_{G_1}(A(\ep)) & \stackrel{h_{G_1}(p)} \lra & h_{G_1}(A)\\
d'' & \mapsto & 1_{G_1(A)} \\
\downarrow & & \downarrow \\
h_{G}(A(\ep)) & \stackrel{h_G(p)} \lra & h_G(A) \\
\psi(d) & \mapsto & 1_{G(A)}.
\end{array}
$$

\medskip
We now want to build an inverse for $\psi$. Let $z \in ker(h_{G}(p))$
i.e. $h_G(p)z=1_{G(A)}$, where:
$$
\begin{array}{ccc}
h_G(p):h_G(A(\ep))=\Hom(\uspec A(\ep), G) & \lra &
h_G(A)=\Hom(\uspec A, G) \\ \\
\end{array}
$$
Here $z$ factors via $G_1$ and this is because $1_{G(A)}$
splits via $G_1$
(recall $p:A(\ep) \lra A$ and it induces
$p^\#:\uspec A \lra \uspec A(\ep)$).
Since $z$ factors via $G_1$, that is $z:\uspec A(\ep) \lra G_1 \lra G$,
this provides immediately a map $\uspec A(\ep) \lra G_1$
corresponding to an element in $\Hom_\smod (m_{1_G}/m_{1_G}^2, A)$.
\end{proof}

\bigskip

\section{The Lie superalgebra of a supergroup scheme}


We now want to show that, for any supergroup scheme $ G $,
$\Lie( G )$ is a Lie superalgebra.

\medskip

To ease the notation, throughout this section $G$ will denote
the functor of points of a supergroup scheme.

\medskip

\begin{observation}
$\Lie( G )$ has a structure of $\cO_k$-module. In fact
let $u_a:A(\ep) \lra A(\ep)$ be the endomorphism,
$u_a(1)=1$, $u_a(\ep)=a\ep$, for $a \in A_0$. $\Lie( G )$ admits a
$\cO_k$-module structure, i.e. there is a natural transformation
$\cO_k \times \Lie( G ) \lra \Lie( G )$, such that for any superalgebra $A$
$$
(a,x) \mapsto ax=_{\defi}\Lie( G )(u_a)x, \qquad
a\in \cO_k(A), \quad x\in \Lie( G )(A).
$$
For subgroups of
$\rGL(m|n)(A)$, $ax$ corresponds to the multiplication of the
matrix $x$ by the
even scalar $a$.
\end{observation}

\medskip

We now want to define a natural transformation $[,]:\Lie( G ) \times \Lie( G )
\lra \Lie( G )$ which has the properties of a superbracket.

\medskip

Let $\rGL(\Lie( G ))(A)$ be  the (multiplicative)
group of linear automorphisms and
$\End(\Lie( G ))(A)$ be the (additive) group of linear
endomorphisms of
$\Lie( G )(A)$. The natural $\cO_k$-module structure of
$\Lie( G )$ gives two group functors (one multiplicative the other additive)
$$
\rGL(\Lie( G )):\salg \lra \sets,
\qquad
\End(\Lie( G )):\salg \lra \sets.
$$
One can  check
$$
\Lie(\rGL(\Lie(G)))=\End(\Lie(G)).
$$

\medskip

\begin{definition}
The \textit{ adjoint action} $\Ad$ of $ G $ on
$\Lie( G )$ is defined as the natural transformation
$$
\begin{array}{c}
\Ad:  G    \lra  \rGL(\Lie( G )) \\ \\
\Ad(g)(x)=G (i)(g)x  G (i)(g)^{-1}, \quad g \in G(A),
\quad x \in \Lie(G)(A).
\end{array}
$$

The \textit{ adjoint action} $\ad$ of $\Lie( G )$ on $\Lie( G )$ is defined as
$$
\ad=_{\defi}\Lie(\Ad):\Lie( G ) \lra
\Lie(\rGL(\Lie( G )))=\End(\Lie( G )).
$$

On $\Lie(G)$ we are ready to define a \textit{bracket}:
$$
[x,y]=_{\defi} \ad(x)y, \qquad x,y \in \Lie(G)(A).
$$
\end{definition}

\medskip

\begin{observation}\label{conjugation}
One can check that 
$$
\Ad(g)=\Lie(c(g))
$$
where $c(g):G(A) \lra G(A)$, $c(g)(x)=gxg^{-1}$.
\end{observation}

\medskip

Our goal is now to prove that $[,]$ is a Lie bracket.

In the next example we work out the bracket for $\rGL_{m|n}$.
This example will be crucial for the next propositions.

\medskip

\begin{example} \label{glbracket}
We want to see that in the case of $\rGL_{m|n}$, the Lie bracket
$[\;,\;]$ coincides with the bracket
defined in Example \ref{gl}. We have:
$$
\begin{array}{cccc}
\Ad:& \rGL(A) & \lra & \rGL(\Lie(\rGL_{m|n}))(A)=
\rGL(\rM_{m|n}(A)) \cr\cr
& g & \mapsto & \Ad(g),
\end{array}
$$
Since $ G (i):\rGL_{m|n}(A) \lra \rGL_{m|n}(A(\ep))$ is an inclusion,
if we identify $\rGL_{m|n}(A)$ with its image we can write:
$$
\begin{array}{c}
\Ad(g)x=gxg^{-1}, \qquad x\in \rM_{m|n}(A).
\end{array}
$$
By definition we have: $\Lie(\rGL(\rM_{m|n}))(A)=\{ 1+\ep \be
\quad | \quad \be \in \rGL(\rM_{m|n})(A) \}$
So we have, for $a,b \in \rM_{m|n}(A) \cong \Lie(\rGL_{m|n})(A)=$
$\{1+\ep a\;|\;a \in \rM_{m|n}(A)\}$:
$$
\ad(1+\ep a)b=(1+\ep a)b(1-\ep a)=b+(ab-ba)\ep=b+\ep [a,b].
$$
Hence $\ad(1+\ep a)=\id+\ep \be(a)$, with $\be(a)= [a,\;]$.

\end{example}

\medskip

It is important to observe that in $G(A(\ep))$ it is customary
to write the product of two elements $x$ and $y$ as $xy$. However
as elements of $Lie(G)(A)$, their product is written as $x+y$
(hence the unit is $0$ and the inverse of $x$ is $-x$).
In order to be able to switch between these two way of writing it is
useful to introduce the notation $e^{\ep x}$.

\medskip

\begin{definition}
Let $\phi: A(\ep) \lra B$ be a superalgebra morphism such that
$\phi(\ep)=\al$. We define $e^{\al x}=G(\phi)(x)$. 
\end{definition}

\medskip

The following properties are immediate:

\medskip

1. $e^{\ep x}=x$,

2. $e^{\al(x+y)}=e^{\al x} e^{\al y}$,

3. $e^{\al(ax)}=e^{a\al x}$,

4. If $f:G \lra H$, $f(e^{\al x})=e^{\al Lie(f) x}$.

\medskip

\begin{observation} \label{Ad}
Observe that Property 4 above and the Example \ref{glbracket}
give us:
$$
\Ad(e^{\ep x})y=y+\ep[x,y]=(\id+\ep \ad(x))y
$$ 
\end{observation}

\begin{lemma} \label{exp}
Let the notation be as above and let $\ep$, $\ep'$ be two
elements with square zero.
$$
e^{\ep x} e^{\ep' y} e^{-\ep x}e^{-\ep' y}=e^{\ep \ep'[x,y]}
$$
\end{lemma}

\begin{proof} By the Property 4 and Observation \ref{Ad} 
we have that:
$$
e^{\ep x} e^{\ep' y} e^{-\ep x}=e^{\ep' \Ad(e^{\ep x}) y}
=e^{\ep'(y+\ep[x,y])}.
$$
So we have
$$
e^{\ep x} e^{\ep' y} e^{-\ep x}=e^{\ep' y}e^{\ep \ep'[x,y]}=
e^{\ep \ep'[x,y]}e^{\ep' y}
$$
which gives the result.

\end{proof}

\medskip

\begin{proposition} The bracket $[,]$ is  antisymmetric.
\end{proposition}

\begin{proof} From the previous proposition we have that:
$$
e^{\ep \ep'[x,y]}=
e^{\ep \ep'[-y,x]}
$$
from which we get the result.
\end{proof}
\medskip

\begin{proposition}
Let $\rho: G \lra \rGL(V)$ be a morphism of supergroup
functors.
Then $\Lie(\rho): \Lie(G) \lra \Lie(V)$ is a Lie
superalgebras morphism.
\end{proposition}

\begin{proof} Using the notation introduced previously we have that
by Observation \ref{Ad} and Property 4:
$$
\rho(e^{\ep x})=e^{\ep \Lie(\rho) x} = \id+ \ep \Lie(\rho) x
$$
Using Proposition \ref{exp} we have:
$$
\rho(e^{\ep \ep' [x,y]})=\rho(e^{\ep x})
\rho(e^{\ep' y})\rho(e^{-\ep x})\rho(e^{\ep' y}).
$$
Hence using Property 4:
$$
\id + \ep \ep' \Lie(\rho)[x,y]=(\id+\ep \Lie(\rho)x)
(\id+\ep' \Lie(\rho)y)(\id-\ep \Lie(\rho)x)(\id-\ep \Lie(\rho)y),
$$
which immediately gives:
$$
\Lie(\rho)[x,y]=[\Lie(\rho)(x), \Lie(\rho)(y)].
$$
\end{proof}

\begin{proposition} The bracket $[,]$ satisfies the Jacobi identity.
\end{proposition}

\begin{proof} In the previous proposition take $\rho=\Ad$.
Then we have:
$$
[\ad(x), \ad(y)]=\ad([x,y]), \qquad \forall x,y \in \Lie(G)(A)
$$
which gives us immediately the Jacobi identity.
\end{proof}
\medskip

\begin{corollary} The natural transformation
$[,]:\Lie( G ) \times \Lie( G ) \lra \Lie( G )$ defined as
$$
[x,y]=_{\defi}\ad(x)y, \qquad x,y \in \Lie( G )(A).
$$
is a Lie bracket for all $A$.
\end{corollary}

\begin{proof} Immediate from previous propositions.
\end{proof}

\medskip

\bigskip

\section{Affine algebraic supergroups. Linear representations.}
\medskip

We now want to restrict our attention to the case of the supergroup scheme
$ G $ to be an affine algebraic group.

\medskip

Let's recall few facts from chapter 1.
Let $A$ and $B$ be superalgebras.
A morphism of algebras $f:A \lra B$ (not necessarily in $\Hom_\salg (A,B)$)
is called
\textit{ even} if it preserves parity, 
\textit{ odd} if it ``reverses'' the parity
i.e. sends even elements in odd elements. Clearly any morphism of
algebras can be
written as sum of an even and an odd one. Recall that
$\Hom_{\salg} (A,B)$ consists only of the even maps.
The set of all morphisms between $A$ and $B$ is called \textit{ inner
Hom}, it is denoted with $\uHom_{\salg}(A,B)$ and it can
be made an object of $\salg$.
Its even part is $\Hom_\salg (A,B)$.

\medskip

\begin{definition} Let $ G $ be an affine algebraic supergroup,
$k[ G ]$ its Hopf superalgebra. We define the additive map
$D:k[ G ] \lra k[ G ]$
a \textit{ left invariant super derivation} if it satisfies
the following properties.

\noindent
1. $D$ is $k$-linear i. e. $D(a)=0$, for all $a \in k$,

\noindent
2. $D$ satisfies the Leibniz identity,
$
D(fg)=D(f)g+(-1)^{p(D)p(f)}fD(g),
$

\noindent
3. $\Delta \circ D=(id \otimes D) \circ \Delta$, where $\Delta$
denotes the comultiplication in $k[ G ]$.
\end{definition}

\medskip

\begin{observation}
The set $L( G )$ of left invariant
derivations of $k[ G ]$ is a Lie superalgebra with bracket:
$$
[D_1, D_2]=_{\defi} D_1 D_2- (-1)^{p(D_1)p(D_2)} D_2D_1.
$$
\end{observation}

\medskip

\begin{theorem}
{ Let $ G $ be an affine supergroup scheme. Then we have natural
bijections among the sets:

\noindent
a) $L( G )$ left invariant derivations in $\Der(k[ G ], k[ G ])$,

\noindent
b) $\Der(k[ G ],k)$,

\noindent
c) $\Lie( G )$.
}
\end{theorem}

\begin{proof}. Let's examine the correspondence between (a) and (b).
We want to construct a map $\phi:\Der(k[ G ],k) \lra L( G )$.
Let $d \in \Der(k[ G ],k)$. Define $\phi(d)= (\id \otimes d) \D$. Then
$\phi(d) \in \Der(k[ G ], k[ G ])$, moreover it is left invariant as one can
readily check. Vice-versa if $D \in L( G )$ define
$\psi(D)= D \circ \ep$ ($\ep$ is the
counit in the Hopf algebra $k[ G ]$). One can check that
$\psi$ is the inverse of $\phi$.
We now want a correspondence between (b) and (c).
By Theorem \ref{tangentspace} we have that
$\Lie( G ) = \Hom_\smod(T_{1_G}( G )^*,A)=Der(\cO_{ G ,1_G},k)$.
Observe that as in the commutative case:
$$
Der(\cO_{ G ,1_G},k)=Der(k[ G ],k),
$$
that is, the derivation on the localization of the ring $k[ G ]$
is determined by the derivation on the ring itself.
\end{proof}

\medskip

We want to show that, as in the classical case, every affine algebraic
supergroup $ G $ can be embedded into some $\rGL(m|n)$.

\begin{definition} Let $f: X \lra  Y $ be a superscheme morphism.
We say that $f$ is a \textit{ closed immersion} if
the topological map
$|f|:|X| \lra |Y|$ is a homeomorphism of the topological space
$|X|$ onto its image in $|Y|$ and the 
sheaf map $f^*:\cO_Y \lra f_*\cO_X$ is
a surjective morphism of sheaves of superalgebras.
\end{definition}

This means that we may identify $X$  with a closed subscheme of $Y$,
so its sheaf is identified with $\cO_Y/\cI$, for some quasicoherent
sheaf of ideals $\cI$.
If both $X$ and $Y$ are affine superschemes we have immediately that $f$
is a closed immersion if and only if $k[X] \cong k[Y]/I$
for some ideal $I$.

\medskip

We now need to introduce the notion of linear representation of
a supergroup.

\medskip

\begin{definition} Let $V$ be a super vector space
and $ G $ an algebraic supergroup (not necessarily affine). We define
\textit{ linear representation} a natural transformation $\rho$
$$
\rho:  h_G  \lra \End(V),
$$
where $\End(V)$ is the functor
$$
\End(V): \salg \lra \sets, \qquad
\End(V)(A)=\End(A \otimes V).
$$
Here $\End(A \otimes V)$ denotes the endomorphisms of the $A$-module
$A \otimes V$ preserving the parity.
We will also say that $ G $ \textit{ acts} on $V$.
\end{definition}

\medskip

Assume now that $ G $ is an affine algebraic supergroup, $k[G]$
the Hopf superalgebra representing it, $\Delta$ and $\ep$ its comultiplication
and counit respectively.

\medskip

\begin{definition} Let $V$ be a supervector space. We say
that $V$ is a \textit{ right $ G $-comodule} if there exists a linear map:
$$
\Delta_V:V \lra V \otimes k[G]
$$
called a \textit{ comodule map} with the properties:

\noindent
1) $(\Delta_V \otimes \id_G) \Delta_V=(\id_V \otimes \Delta)\Delta_V$

\noindent
2) $(\id_V \otimes \ep)\Delta_V=\id_V$.

where $\id_G:k[G] \lra k[G]$ is the identity map.
\end{definition}

\medskip

One can also define a left $ G $-comodule in the obvious way.

\medskip

\begin{observation} The two notions of $ G $ acting on $V$
and $V$ being a (right) $G$-comodule are essentially equivalent.
In fact, given a representation $\rho: G  \lra \End(V)$, it defines
immediately a comodule map:
$$
\Delta_V(v)=\rho_{k[G]}(\id_G)v, \qquad 
\id_G \in  h_G(k[G])=\Hom_\salg (k[G],k[G])
$$
where we are using the natural identification (for $A=k[G]$)
$$
\End(V)(A) \cong \Hom_{\smod}(V, V \otimes A ).
$$

Vice-versa if we have a comodule map $\Delta_V$ we can define a
representation in the following way:
$$
\begin{array}{ccc}
\rho_A:h_G(A) & \lra & \End(V)(A) \cong \Hom_{\smod}(V, V \otimes A )\\ \\
g & \mapsto &  v \mapsto (\id \otimes g)(\Delta_V(v))
\end{array}
$$
where $g \in h_G(A)=\Hom_\salg(k[G],A)$.

\end{observation}

\medskip

Let's see this correspondence in a special, but important case.

\medskip

\begin{example} Let's consider the
natural action of  $\rGL_{m|n}$  on $k^{m|n}$:
$$
\begin{array}{ccc}
\rho_A:\rGL_{m|n}(A) & \lra & \End(k^{m|n})(A) \\ \\
g=(g_{ij}) & \mapsto & e_j \mapsto \sum e_i \otimes g_{ij}
\end{array}
$$
where $\{e_j\}$ is the canonical homogeneous basis for the framed
supervector space $k^{m|n}$.
We identify the morphism $g \in \rGL_{m|n}(A) =$ $\Hom_\salg(k[G],A)$
with the matrix with entries $g_{ij}=g(x_{ij})$, where $x_{ij}$'s are
the generators of $k[\rGL(m|n)]$.

This corresponds to the comodule map
$$
\begin{array}{ccc}
\Delta_{k^{m|n}}: k^{m|n} & \lra & k^{m|n}  \otimes k[\rGL(m|n)]\\ \\
e_j & \mapsto & \sum e_j \otimes x_{ij}
\end{array}
$$
where $x_{ij}$ are the generators of the algebra $k[\rGL(m|n)]$.

Vice-versa, given the comodule map as above:
$e_j  \mapsto \sum e_j \otimes x_{ij}$
it corresponds to the representation:
$$
\begin{array}{ccccc}
{\rho}_A & : \rGL_{m|n}(A) & \lra & \End(k^{m|n})(A) &  \\
& g=(g_{ij}) & \mapsto & e_j \mapsto
(\id \otimes g)(\sum e_i \otimes x_{ij})= &
\sum e_i \otimes g_{ij}.
\end{array}
$$
\end{example}

\medskip

\begin{definition} Let $ G $ act on the superspace $V$,
via a representation $\rho$ corresponding to the comodule map
$\Delta_V$. We say
that the subspace $W \subset V$ is \textit{$ G $-stable} if
$\Delta_V(W) \subset W \otimes V$.
Equivalently $W$ is $ G $-stable if $\rho_A(g)(W \otimes A)
\subset W \otimes A$.
\end{definition}

\medskip

\begin{definition} The \textit{ right regular representation}
of the affine algebraic group $ G $ is the representation of $ G $ in
the (infinite dimensional) super vector space $k[G]$
corresponding to the comodule map:
$$
\Delta:k[G] \lra k[G] \otimes k[G].
$$
\end{definition}

\medskip

\begin{proposition} \label{finite}
Let $\rho$ be a linear representation
of an affine algebraic supergroup $G$. Then each finite dimensional
supersubspace of $V$ generates a finite dimensional stable subspace of
$V$.
\end{proposition}

\begin{proof} It is the same as in the commutative case. Let's sketch it.
It is enough to prove for one element $x \in V$.
Let $\Delta_V:V \lra V \otimes k[G]$ be the comodule structure associated
to the representation $\rho$. Let
$$
\Delta_V(x)=\sum_i x_i \otimes a_i
$$
where $\{a_i\}$ is a basis for $k[G]$.

We claim that $\Span_k \{x_i\}$ is a $ G $-stable subspace.

By definition of comodule we have:
$$
(\Delta_V \otimes \id_G)(\Delta_V(x))=(\id_V \otimes \Delta)(\Delta_V(x)),
$$
that is
$$
\sum_j \Delta_V(x_j) \otimes a_j= \sum_i x_i \otimes \Delta(a_i)=
\sum_{i,j} x_i \otimes b_{ij} \otimes a_j.
$$
Hence
$$
\Delta_V(x_j)=\sum_{i} x_i \otimes b_{ij}.
$$
The finite dimensional stable subspace is given by the span
of the $x_i$'s.
\end{proof}

\medskip

\begin{theorem} Let $ G $ be an affine supergroup variety. Then
there exists a closed embedding:
$$
 G  \subset \rGL(m|n)
$$
for suitable $m$ and $n$.
\end{theorem}

\begin{proof} We need to find a
surjective superalgebra morphism
$k[GL(m|n)] \lra k[G]$ for suitable $m$ and $n$.
Let $k[G]=k[f_1 \dots f_n]$, where $f_i$ are homogeneous and
chosen so that $W =\Span\{f_1 \dots f_n\}$ is $G$-stable, according
to the right regular representation. This choice is possible because
of Proposition \ref{finite}.
We have:
$$
\Delta_{k[G]}(f_i)= \sum_j f_j \otimes a_{ij}.
$$
Define the morphims:
$$
\begin{array}{ccc}
k[GL(m|n)] & \lra & k[G] \\
x_{ij} & \mapsto & a_{ij}
\end{array}
$$
where $x_{ij}$ are the generators for $k[GL(m|n)]$. This is the
required surjective algebra morphism.
In fact, since $k[G]$ is both a right and left $ G $-comodule we have:
$$
f_i=(\ep \otimes \id)\Delta(f_i)=(\ep \otimes \id)(\sum_j f_j \otimes a_{ij})=
\sum_j \ep(f_j) \otimes a_{ij}
$$
which proves the surjectivity.
\end{proof}

\medskip

\begin{corollary} $ G $ is an affine supergroup scheme
if and only if it is a closed subgroup of $\rGL(m|n)$.
\end{corollary}


\newpage

\appendix  \label{sgappendix}

\chapter{Appendix} \label{sgappendix}

\section{Categories}

We want to make a brief summary of 
formal properties and definitions relative to categories.
For more details one can see for example \cite{lang}.

\begin{definition}

A \textit{category} $\mathcal{C}$ consists of a collection of objects, 
denoted by $Ob(\mathcal{C})$, and sets of \textit{morphisms} 
between objects.
For all pairs $A,B \in Ob(\mathcal{C})$, we denote the 
set of morphisms from $A$ to $B$ by
$\text{Hom}_{\mathcal{C}}(A,B)$ so that
for all $A,B,C \in \mathcal{C}$, there exists 
an association
\[
\Hom_{\mathcal{C}}(B,C) \times \Hom_{\mathcal{C}}(A,B) 
\lra \Hom_{\mathcal{C}}(A,B)
\] 
called the 
``composition law"  ($(f,g) \rightarrow f \circ g$) which satisfies the properties\\
\noindent (i) the law ``$\circ$" is associative,\\
\noindent (ii) for all $A,B \in Ob(\mathcal{C})$, there exists 
$id_A \in \Hom_{\mathcal{C}}(A,A)$ so that we get 
$f \circ id_A = f$ for all $f \in \Hom_{\mathcal{C}}(A,B)$ and 
$id_A \circ g = g$ for all $g \in \Hom_{\mathcal{C}}(B,A)$,\\
\indent (iii) $\Hom_{\mathcal C} (A,B)$ and  $\Hom_{\mathcal C} (A',B')$
are disjoint unless $A=A'$, $B=B'$ in which case they are equal.
\end{definition}

Once the category is understood, it is conventional to write $A \in \mathcal{C}$ 
instead of $A\in Ob(\mathcal{C})$ for objects.  
We may also suppress the ``$\mathcal{C}$" from $\Hom_{\mathcal{C}}$ and just write $\Hom$
whenever there is no danger of confusion.

Essentially a category is a collection of objects which share some basic
structure, along with maps between objects which preserve that structure.

\begin{example}
Let $\mathcal{G}$ denote the category of groups.  
Any object $G \in \mathcal{G}$ is a group, and for any two groups $G, H \in
Ob(\mathcal{G})$, the set $\Hom_{\mathcal{G}}(G,H)$ is the set of 
group homomorphisms from $G$ to $H$.
\end{example}





\begin{definition}
A category $\mathcal{C}'$ is a \textit{subcategory} of category 
$\mathcal{C}$ if $Ob(\mathcal{C}' )\subset Ob(\mathcal{C})$ and if 
for all $A,B \in\mathcal{C}'$, $\Hom_{\mathcal{C}'} (A,B) \subset 
\Hom_{\mathcal{C}} (A,B)$ so that the composition law ``$\circ$" on 
$\mathcal{C}'$ is induced by that on $\mathcal{C}$.
\end{definition}

\begin{example}
The category $\mathcal{A}$ of \textit{abelian groups} is a 
subcategory of the category of groups $\mathcal{G}$.
\end{example}

\begin{definition}
Let $\mathcal{C}_1$ and $\mathcal{C}_2$ be two categories.  
Then a \textit{covariant [contravariant] functor} 
$F: \mathcal{C}_1 \longrightarrow \mathcal{C}_2$ consists of\\
\noindent (1) a map $F: Ob(\mathcal{C}_1) \longrightarrow Ob(\mathcal{C}_2)$ and\\
\noindent (2) a map (denoted by the same $F$) $F: \Hom_{\mathcal{C}_1} (A,B) 
\longrightarrow \Hom_{\mathcal{C}_2}(F(A),F(B))$ 
[$F: \Hom_{\mathcal{C}_1} (A,B) \longrightarrow
\Hom_{\mathcal{C}_2}(F(B),F(A))$] 
so that\\
\indent (i) $F(id_A) = id_{F(A)}$ and\\
\indent (ii) $F(f \circ g) = F(f) \circ F(g)$ [$F(f \circ g) = F(g) \circ F(f)$]\\
\noindent for all $A,B \in Ob(\mathcal{C}_1)$.
\end{definition}

When we say ``functor'' we mean covariant functor. 
A contravariant functor
$F: \mathcal{C}_1 \lra \mathcal{C}_2$ is the same as a covariant functor
from $\mathcal{C}^o_1 \lra \mathcal{C}_2$ where $\mathcal{C}_1^o$ denotes
the \textit{opposite} category i. e. the category where all morphism arrows
are reversed.

\begin{definition}
Let $F_1, F_2$ be two functors from $\mathcal{C}_1$ to $\mathcal{C}_2$.  
We say that there is a natural transformation of functors 
$\varphi: F_1 \longrightarrow F_2$ if for all 
$A \in \mathcal{C}_1$ there is a set of morphisms 
$\varphi_A: F_1(A) \longrightarrow F_2(A)$ so that for any 
$f \in \Hom_{\mathcal{C}_1}(A,B)$ ($B \in \mathcal{C}_1$), 
the following diagram commutes:
\begin{equation}
\begin{array}{lcr}
F_1(A) & \stackrel{\varphi_A}{\longrightarrow} & F_2(A)\\
 _{F_1(f)} \downarrow & & \downarrow_{F_2(f)}\\
F_1(B) & \stackrel{\varphi_B}{\longrightarrow} & F_2(B).
\end{array}
\end{equation}
\end{definition}

We say that the family of functions $\varphi_A$ is \textit{functorial} in $A$.  

The notion of equivalence of categories is important since it allows
to identify two categories which are apparently different.

\begin{definition}

We say that two categories $\mathcal{C}_1$ and $\mathcal{C}_2$ are
{\it equivalent} if there exists two functors $F:\mathcal{C}_1 \lra \mathcal{C}_2$ and
$G:\mathcal{C}_1 \lra \mathcal{C}_2$ such that $FG \cong id_{\mathcal{C}_2}$, 
$GF \cong id_{\mathcal{C}_1}$
(where $id_{\mathcal{C}}$ denotes the identity functor of a given category,
defined in the obvious way).
\end{definition}

Next we want to formally define what it means for a functor to be 
\textit{representable}.

\begin{definition}
Let $F$ be a functor from the category $\mathcal{C}$ to 
the category of sets $\mathcal{S}$.  We say that $F$ 
is \textit{representable by $X \in \mathcal{C}$} if for all $A \in \mathcal{C}$,
$$
\begin{array}{c}
F(A) = \Hom_{\mathcal{C}}(X,A), \\
F(f): F(A) \lra F(B), \quad F(f)(\al)=_{def} f \cdot \al\qquad 
\end{array}
$$
for all $f: A \lra B$.

\end{definition}

We end our small exposition of categories by the constructing the fibered product which is
be important in chapter \ref{sgchap3}.

\begin{definition}

Given functors $A$, $B$, $C$, from a category $\mathcal{C}$ to the category of $\sets$, and given natural transformations $f:A \lra C$, $g:B \lra C$ the \textit{fibered product} $A \times_C B$ is the universal object making the following diagram commute:
\begin{equation*}
\begin{array}{ccc}
A \times_C B & \longrightarrow & B \\
\downarrow &  & \downarrow g\\
A & \stackrel{f} \longrightarrow & C.
\end{array}
\label{cd}
\end{equation*}
One can see that:
\[
(A \times_C B)(R)=A(R) \times_{C(R)} B(R)=
\{(a,b) \in A \times B| f(a)=g(b)\}
\]
\end{definition}

One can see that if $g$ is injective, that is $g_R:B(R) \subset C(R)$
we have that $(A \times_C B)(R)=f^{-1}(B(R))$.

The language of categories allows us to make (and prove) some sweeping 
generalizations about geometric objects without too much ``forceful"
computation.  In particular, it also allows us to generalize the notion of a ``point" 
to a $T$-point; this allows us to make more intuitive calculations with 
supergeometric objects.  The main categories we discuss in this exposition 
are the categories of $C^{\infty}$-supermanifolds, super Lie groups 
(a subcategory of $C^{\infty}$-supermanifolds), superschemes, and 
super algebraic groups.

\section{SuperNakayama's Lemma and Projective Modules}
\label{projective}

Let $A$ be a commutative superalgebra.

\begin{definition}
A projective $A$-module $M$ is a direct summand of $A^{m|n}$. In other words it is a  projective module in the classical sense respecting the grading: $M_0 \subset A^{m|n}_0$, $M_1 \subset A^{m|n}_1$.
\end{definition}

\begin{observation}
As in the classical setting being projective is equivalent to the exactness of the functor $\Hom(M, \_ )$.
\end{observation}

We want to show that a projective $A$-module has the property of being locally
free, that is its localization $L_p$ into primes $p$ of $A_0$ is free as
$A_p$-module. This result allows to define the \textit{rank} of a
projective module as it happens in the ordinary case.

We start with a generalization of the Nakayama's lemma.

\begin{lemma}
\label{snakayama} (super Nakayama's Lemma)
Let $A$ be a local supercommutative ring with maximal
homogeneous ideal $\mathfrak{m}$. Let $E$ be a finitely generated
module for the ungraded ring $A$.\\

\noindent (i) If $\mathfrak{m}E=E$, then $E=0$; more
generally, if $H$ is a submodule of $E$ such that $E=\frak mE+H$,
then $E=H$.\\

\noindent (ii) Let $(v_i)_{1\le i\le p}$ be a basis
for the $k$-vector space $E/\mathfrak{m}E$ where $k=A/\mathfrak{m}$. Let
$e_i\in E$ be above $v_i$. Then the $e_i$ generate $E$. If $E$ is
a supermodule for the super ring $A$, and $v_i$ are homogeneous
elements of the super vector space $E/\mathfrak{m}E$, we can choose the
$e_i$ to be homogeneous too (and hence of the same parity as
the $v_i$).\\

\noindent (iii) Suppose $E$ is projective, i.e. there is a $A$-module $F$ such that $E\oplus
F=A^N$ where $A^N$ is the free module for the ungraded ring $A$
of rank $N$. Then $E$ (and hence $F$) is free, and the
$e_i$ above form a basis for $E$.
\end{lemma}

\begin{proof}
The proofs are easy extensions of the ones in the
commutative case. We begin the proof of (i) with the following
observation: if $B$ is a commutative local ring with $\mathfrak{n}$ a
maximal ideal, then a square matrix $R$ over $B$ is invertible if
and only if it is invertible modulo $\mathfrak{n}$ over the
field $B/\mathfrak{n}$.  In fact if this is so, $\text{det} (R)\notin \mathfrak{n}$
and so is a unit of $B$. This said, let $u_i$, $(1< i < N)$
generate $E$. If $E=\mathfrak{m}E$, we can find $m_{ij}\in \mathfrak{m}$ so
that $u_i=\sum _j m_{ij} e_j$ for all $i$. Hence, if $L$ is the
matrix with entries $\delta _{ij}-m_{ij}$, then
\[
L\begin{pmatrix} u_1\cr u_2\cr\vdots \cr u_N\cr \end{pmatrix} = 0.
\]
It is now enough to prove that $L$ has a left inverse.  Then multiplying
the above from the left by $P$, we get $u_i=0$ for all $i$ and so
$E=0$. It is even true that $L$ is invertible. To
prove this, let us consider $B=A/J$ where $J$ is the ideal generated
by $A_1$. Since $J\subset \mathfrak{m}$ we have
\[
A\longrightarrow B=A/J\longrightarrow k=A/\mathfrak{m}.
\]
Let $L_B$ (resp. $L_k$) be the reduction of $L$ modulo $J$ (respectively modulo
$\mathfrak{m}$). Then $B$ is local, and its maximal ideal is $\mathfrak{m}/J$ where $L_k$ is the reduction of $L_B$ mod $\frak m/J$. But $B$ is
commutative and $L_k=I$, and so $L_B$ is invertible. But then $L$
is invertible.  If more generally we have $E=H+\mathfrak
{m}E$, then $E/H=\mathfrak{m}(E/H)$ and so $E/H=0$, which is to say that $E=H$.

To prove (ii), let $H$ be the submodule of $E$ generated by the $e_i$.
Then $E=\mathfrak{m}E+H$ and so $E=H$.

We now prove (iii). Clearly $F$ is also finitely generated. We
have $k^N=A^N/\mathfrak{m}^N=E/\mathfrak{m}E \oplus F/\mathfrak{m}F$. Let 
$(w_j)$
be a basis of $F/\mathfrak{m}F$ and let $f_j$ be elements of $F$ above
$w_j$. Then by (ii), the $e_i, f_j$ form a basis of $A^N$ while the
$e_i$ (respectively the $f_j$) generate $E$ (resp. $F$). Now there are exactly
$N$ of the $e_i, f_j$, and so if $X$ denotes the $N\times N$
matrix with columns $e_1, \dots , f_1, \dots $, then for some
$N\times N$ matrix $Y$ over $A$ we have $XY=I$.  Hence $X_BY_B=I$
where the suffix ``$B$" denotes reduction modulo $B$.  However, $B$ is
commutative and so $Y_BX_B=I$.  Thus $X$ has a left inverse over
$A$, which must be $Y$ so that $YX=I$. If there is a linear
relation among the $e_i$ and the $f_j$, and if $x$ is the column vector whose
components are the coefficients of this relation, then $Xx=0$; but
then $x=YXx=0$.  In particular $E$ is a free module with basis
$(e_i)$.
\end{proof}

We now wish to give a characterization of projective modules.

\begin{theorem}
Let $M$ be a finitely generated $A$-module, $A$ finitely
generated over $A_0$ where $A_0$ noetherian. Then

\noindent
i) $M$ is projective if and only if $M_p$ is free for
all $p$ primes in $A_0$ and

\noindent
ii) $M$ is projective if and only if $M[f_i^{-1}]$ free
for all $f_i$'s such that $(f_1 \dots f_r)=A_0$.
\end{theorem}

\begin{proof} (i) If $M$ is projective, by part (iii) of Nakayama's Lemma, we have that $M_{p}$ is free since it is a module over the supercommutative ring $A_p$.




















Now assume that $M_p$ is free for all primes $p \in A_0$.
Recall that
\[
\Hom_{A[U^{-1}]}(M[U^{-1}], N[U^{-1}])=\Hom_A(M,N)[U^{-1}]
\]
for $U$ a multiplicatively closed set in $A_0$.
Recall also that given $A_0$-modules $N$, $N'$, $N''$, we have that
$0 \lra N' \lra N \lra N''$ is exact if and only if
$0 \lra N'_p \lra N_p \lra N''_p$ is exact for all prime $p$ in $A_0$.
So given an exact sequence
$0 \lra N' \lra N \lra N''$, since $M_p$ is free,
we obtain the exact sequence
\[
0 \lra \Hom(M_p,N_p') \lra \Hom(M_p,N_p) \Hom(M_p,N_p'') \lra 0
\]
for all the primes $p$.  Hence by the previous observation,
$$
0 \lra \Hom(M,N')_p \lra \Hom(M,N)_p \lra \Hom(M,N'')_p \lra 0,
$$
and
$$
0 \lra \Hom(M,N') \lra \Hom(M,N) \lra \Hom(M,N'') \lra 0.
$$
Hence $M$ is projective.

\noindent (ii) That $M_p$ is free for all primes $p$ is equivalent
to $M[f_i^{-1}]$ being free for $(f_1 \dots f_r)=A_0$ is a standard fact
of commutative algebra and can be found in \cite{eisenbud} p. 623 for example.
\end{proof}

\begin{remark}
As in the ordinary setting we have a correspondence between 
projective $A$-modules and locally free sheaves on $\spec A_0$. In this 
correspondence, given a projective $A$-module $M$, we view
$M$ as an $A_0$-module and build a sheaf of modules $\cO_M$ on $A_0$.  The
global sections of this sheaf are isomorphic to $M$ itself, and locally, i.e. on the open sets 
$U_{f_i}=\{p \in \spec A_0| (f_i) \not\subset p\}$, $f_i \in A_0$,
\[
\cO_M(U_{f_i})=M[f_i^{-1}].
\]
More details on this construction can be found, for example, in \cite{Harts} chapter 2.
\end{remark}


\bibliography {bib/network,bib/naming}    

\bibliographystyle {uclathes}

\bibliography {references}

\begin{thebibliography}{aaa}

\bibitem{am} M. F. Atyah, I. G. McDonald
\textit{ Introduction to commutative algebra},
Perseus books , 1969. 

\bibitem{Berezin}
F.~A.~Berezin.  \textit{Introduction to superanalysis}.  
D.~Reidel Publishing Company, Holland, 1987.

\bibitem{bmpz} Yu. A. Bahturin,
A. A. Mikhalev, V. M. Petrogradsky,  M. V. Zaicev
\textit{Infinite dimensional Lie Superalgebras}
De Gruyter Expositions in Mathematics, 7, Walter de Gruyter, 1992.

\bibitem{bbh} C. Bartocci, U. Bruzzo, D. Hernandez-Ruiperez.
\textit{The Geometry of Supermanifolds}.
MAIA 71, Kluwer Academic Publishers, 1991.

\bibitem{bo} A. Borel, \textit{Linear algebraic groups}. Springer
Verlag, New York, 1991.

\bibitem{eh} D. Eisenbud and J. Harris, \textit{The geometry of
schemes.} Springer Verlag, New York, 2000.

\bibitem{DM}
P.~Deligne, J.~W.~Morgan.  Notes on supersymmetry following Bernstein.
\textit{Quantum fields and strings; a course for mathematicians, Vol.~1} 
(Princeton, NJ, 1996/1997), 41-96, Amer. Math. Soc., Providence, RI, 1999.


\bibitem{dewitt}
B. DeWitt. \textit{Supermanifolds}. Cambridge University Press,
London, 1984.

\bibitem{dg} M. Demazure, P. Gabriel, \textit{Groupes
Alg\'ebriques, Tome 1.} Mason$\&$Cie, \'editeur. North-Holland
Publishing Company, The Netherlands, 1970. 


\bibitem{eisenbud} D. Eisenbud, \textit{ Commutative algebra
with a view towards algebraic geometry}, GTM Springer Verlag, 1991.

\bibitem{ferrara} 
Daniel Z. Freedman, P. van Nieuwenhuizen , S. Ferrara  
\textit{Progress Toward A Theory Of Supergravity},
Phys.Rev.D 13, 3214-3218, 1976.

\bibitem{Groth}
A.~Grothendieck, M.~Demazure, M.~Artin [et al.].  
\textit{Sch\'{e}mas en groupes; s\'{e}minaire de g\'{e}om\'{e}trie
  alg\'{e}brique du Bois Marie}, 1962/64, SGA 3. Springer-Verlag, 
Berlin-New York, 1970.

\bibitem{Harts}
R.~Hartshorne.   \textit{Algebraic geometry}. 
Graduate Text In Mathematics.  
Springer-Verlag, New York, 1977.

\bibitem{ka} V. G. Kac, \textit{Lie superalgebras} Adv. in Math. 
{26} (1977) 8-26.

\bibitem{Kostant}
B.~Kostant.  Graded manifolds, graded Lie theory, and prequantization. 
\textit{Differential geometrical methods in mathematical physics} 
(Proc. Sympos., Univ. Bonn, Bonn, (1975), pp. 177--306. 
\textit{Lecture Notes  in Math.}, Vol. 570, Springer, Berlin, 1977.

\bibitem{lang}
S. Lang. \textit{ Algebra}, Addison Wesley, 1984.

\bibitem{Leites}
D.~A.~Leites.  \textit{Introduction to the theory of supermanifolds}.  
Russian Math. Surveys \textbf{35}:~1 (1980), 1-64.


\bibitem{Manin}
Y.~I.~Manin.   \textit{Gauge field theory and complex geometry}; 
translated by N. Koblitz and J.R. King.  Springer-Verlag, Berlin-New York, 1988.



\bibitem{Mumford}
D.~Mumford.   \textit{The red book of varieties and schemes}.  
Springer-Verlag, Berlin-New York, 1988.



\bibitem{ss} A. Salam, J. Strathdee.
\textit{Super-gauge transformations}, Nucl. Phys. B, 
\textbf{76}, (1974), 477-482.


\bibitem{Tuyn}
G. M. ~Tuynman.   \textit{Supermanifolds and supergroups}.
MAIA 570, Kluwer Academic Publishers, 2004.


\bibitem{VSV1}
V.~S.~Varadarajan.  \textit{Lie groups, Lie algebras, and their
  representations}.  Graduate Text in Mathematics.  Springer-Verlag, New York, 1984.



\bibitem{VSV2}
V.~S.~Varadarajan.  \textit{Supersymmetry for mathematicians: an
  introduction}.  Courant Lecture Notes.  Courant Lecture Notes Series, New York, 2004.


\bibitem{Water}
W.~C.~Waterhouse.  \textit{Introduction to affine group schemes.}  
Graduate Texts in Mathematics, 66. Springer-Verlag, New York-Berlin, 1979.


\bibitem{wz} J. Wess, B. Zumino. \textit{ Supergauge transformations
in four dimensions}. Nucl. Phys. B, \textbf{70}, (1974), 39-50.


\end{thebibliography}

\bibliographystyle {plain}





\end{document}